\pdfoutput=1
\RequirePackage{ifpdf}
\ifpdf 
\documentclass[pdftex]{sigma}
\else
\documentclass{sigma}
\fi

\usepackage{mathrsfs}

\numberwithin{equation}{section}
\newtheorem{Theorem}{Theorem}[section]
\newtheorem{Corollary}[Theorem]{Corollary}
\newtheorem{Lemma}[Theorem]{Lemma}
\newtheorem{Proposition}[Theorem]{Proposition}
 { \theoremstyle{definition}
\newtheorem{Definition}[Theorem]{Definition}
\newtheorem{Example}[Theorem]{Example}
\newtheorem{Remark}[Theorem]{Remark} }

\begin{document}

\allowdisplaybreaks

\renewcommand{\thefootnote}{$\star$}

\renewcommand{\PaperNumber}{009}

\FirstPageHeading

\ShortArticleName{The Heisenberg Relation~-- Mathematical Formulations}

\ArticleName{The Heisenberg Relation~-- Mathematical\\
Formulations\footnote{This paper is a~contribution to the Special Issue on Noncommutative Geometry and
Quantum Groups in honor of Marc A.~Rief\/fel. The full collection is available at \href{http://www.emis.de/journals/SIGMA/Rieffel.html}{http://www.emis.de/journals/SIGMA/Rieffel.html}}}

\Author{Richard V.~KADISON~$^\dag$ and Zhe LIU~$^\ddag$}

\AuthorNameForHeading{R.V.~Kadison and Z.~Liu}

\Address{$^\dag$~Department of Mathematics, University of Pennsylvania, USA}
\EmailD{\href{mailto:kadison@math.upenn.edu}{kadison@math.upenn.edu}}
\URLaddressD{\url{http://www.math.upenn.edu/~kadison/}}

\Address{$^\ddag$~Department of Mathematics, University of Central Florida, USA}
\EmailD{\href{mailto:zhe.liu@ucf.edu}{zhe.liu@ucf.edu}} \URLaddressD{\url{http://math.cos.ucf.edu/~zheliu/}}

\ArticleDates{Received July 26, 2013, in f\/inal form January 18, 2014; Published online January 25, 2014}

\Abstract{We study some of the possibilities for formulating the Heisenberg relation of quantum mechanics
in mathematical terms.
In particular, we examine the framework discussed by Murray and von Neumann, the family (algebra) of
operators af\/f\/iliated with a~f\/inite factor (of inf\/inite linear dimension).}


\Classification{47L60; 47L90; 46L57; 81S99}
\begin{flushright}
\begin{minipage}{110mm} \it
Dedicated, with affection, to Marc Rieffel on the occasion of his 75th birthday by his
proud $($much older$)$ ``mathematical father'' and his admiring $($much younger$)$ ``mathematical sister''.
\end{minipage}
\end{flushright}

\renewcommand{\thefootnote}{\arabic{footnote}}
\setcounter{footnote}{0}

\section{Introduction}

The celebrated Heisenberg relation,
\begin{gather*}
QP-PQ=i\hbar I,
\end{gather*}
where $\hbar=\frac{h}{2\pi}$ and $h$ is Planck's experimentally determined quantum of action
${\approx}6.625\times10^{-27}$~erg~sec, is one of the basic relations (perhaps, the {\it most} basic
relation) of quantum mecha\-nics.
Its form announces, even before a~tentative mathematical framework for quantum mecha\-nics has been
specif\/ied, that the mathematics of quantum mechanics {\it must be non-commutative}.
By contrast, the mathematics suitable for classical mechanics is the commutative mathematics of algebras of
real-valued functions on spaces (manifolds).

Our program in this article is to study specif\/ic mathematical formulations that have the attributes
necessary to accommodate the calculations of quantum mechanics, more particularly, to accommodate the
Heisenberg relation.
We begin that study by noting the inadequacy of two natural candidates.
We turn, after that, to a~description of the ``classic'' representation and conclude with a~model,
suggested by von Neumann, especially suited to calculations with unbounded operators.
Von Neumann had hoped that this model might resolve the mathematical problems that the founders of quantum
mechanics were having with those calculations.
In connection with this hope and the Heisenberg relation, we answer a~question that had puzzled a~number of
us.

\looseness=-1
There is a~very def\/inite expository (review) aspect to this article, with, nevertheless, much new
material and some new results.
It is intended that it can be read completely (even by non-experts) without reference to other resources.
As a~result, there is a~substantial amount of material copied from~\cite{KRI,KRII,KRIII,KRIV} and~\cite{Z}.
We have the reader's convenience very much in mind.

\section{Some physical background}

At the beginning of the twentieth century, the study of dynamical systems comprised of ``large'' material
bodies~-- those not af\/fected (signif\/icantly) by the act of observing them~-- was fully developed.
Physicists could be justly proud of what they had accomplished.
Many thought that physics was ``complete''; all that was needed were clever algorithms that would make it
possible to f\/inish the calculations resulting from applying the techniques of classical mechanics to
``small'' systems that resembled macroscopic mechanical systems.
However, there were surprising dif\/f\/iculties: predictions made on the basis of classical mechanical
models of ``small'' systems (those related to and measurable in the elementary units, of which, it was
believed, all matter is composed) were at wide variance with the data observed.

As workable formulae were developed, formulae that yielded numerical results more in line with experimental
data, one pattern stood out again and again: some process that classical mechanical computation would
predict should occur ``continuously'' seemed to occur in discrete steps.
The f\/irst instance of this, the one that may be identif\/ied as the origin of quantum mecha\-nics, was
Planck's meticulous development of his radiation formula (1900) and the introduction of his physical
``quantum of action'', which in the prechosen units of those times was experimentally determined to be
(approximately) ${\approx}6.625\times10^{-27}$ erg sec.
The formula
\begin{gather*}
\frac{8\pi hc\lambda^{-5}}{e^{hc/k\lambda T}-1}{\rm d}\lambda
\end{gather*}
expresses the energy per unit volume inside a~cavity with ref\/lecting walls associated with the wave
lengths lying between $\lambda$ and $\lambda+{\rm d}\lambda$ of a~full (black body) radiator at (absolute)
temperature~$T$.
Try as he could to explain his formula purely in terms of Maxwell's classical electromagnetic theory,
Planck could not rid himself of the need for the assumption that energy emitted and absorbed by one of the
basic units of the radiating system, a~linear (harmonic) oscillator of frequency~$\nu$ occurred as integral
multiples of~$h\nu$.

If Planck assumed that there is no smallest non-zero amount of energy that each unit (oscillator) could
emit, then he was led to the classical Rayleigh--Jeans radiation formula $8\pi kT\lambda^{-4}{\rm d}\lambda$,
which is a~good approximation to Planck's formula for larger wave lengths.
However, as $\lambda \rightarrow 0$, the energy associated with the oscillator of this wave length tends to
$\infty$; and the total energy per unit volume in the cavity due to this part of the spectrum is, according
to the Rayleigh--Jeans formula,
\begin{gather*}
8\pi kT\int_0^r\lambda^{-4}{\rm d}\lambda=\lim_{\lambda\to0}(8\pi kT/3)\big(\lambda^{-3}-r^{-3}\big),
\end{gather*}
which also tends to $\infty$ as $\lambda\rightarrow 0$.
This breakdown of the formula, being associated with short wave lengths, those corresponding to the
ultraviolet end of the spectrum, was termed ``ultraviolet catastrophe''.

Planck was forced to assume that there was the ``quantum of action'', $h$ erg sec, for his formula to agree
with experimental measurements at the high frequency as well as the low frequency spectrum of the radiation.
It is no small irony that, for some years, Planck was deeply disappointed by this shocking break with the
principles of classical mechanics embodied in his revolutionary discovery.
It was not yet clear to him that this discovery was to become the fundamental physical feature of ``small''
systems.
Not many years later, others began to make the assumption of ``quanta'' in dif\/ferent physical processes
involving small systems.

Any discussion of the inadequacy of classical mechanics for explaining the phenomena noted in experiments
on systems at the subatomic scale must include some examination of the startling evidence of the dual,
corpuscular (material ``particle'') and wave, nature of light.
Maxwell (1865) had shown that the speed of electromagnetic-wave propagation in a~vacuum is a~constant
(representing the ratio of the electromagnetic unit of charge to the electrostatic unit, abcoulomb/esu)
$3\times 10^{10}$ cm/sec (approximately).
He noted that this is (close to) the speed of light and concluded that light is a~form of electromagnetic
wave.
Measurements indicate that the wave length of the visible spectrum lies between $0.000040$ cm (${=}4000\times
10^{-8}$~cm ${=}4000$~{\AA}ngstr{\o}m) at the violet end and $0.000076$ cm (${=}7600$~\AA) at the red end.
Above this, to 0.03 cm is the infrared spectrum.
Below the violet is the ultraviolet spectrum extending to 130~\AA, and below this, the $X$-rays from 100~\AA\ to 0.1~\AA, and the $\gamma$-rays from 0.1~\AA\ to 0.005~\AA.

Another type of wave behavior exhibited by light is the phenomenon of {\it polarization}.
A pair of thin plates cut from a~crystal of tourmaline allows no light to pass through it if one is held
behind the other with their optical axes perpendicular.
As one of the plates is rotated through a~$90^\circ$ angle, more and more of the light passes through, the
maximum occurring when the axes are parallel~-- an indication of light behaving as a~transverse wave.

The phenomenon of (wave) {\it interference} provides additional evidence of the transverse-wave character
of light.
Two waves of the same frequency and amplitude are superimposed.
If they are in phase, they ``reenforce'' one another.
If they are in phase opposition, they cancel.

Further evidence of the wave nature of light is found in the phenomenon of {\it diffraction}~-- the
modif\/ication waves undergo in passing the edges of opaque bodies or through narrow slits in which the
wave direction appears to bend producing fringes of reenforcement and cancelation (light and dark).
A {\it diffraction grating}, which consists, in essence, of a~transparent plate on which parallel, evenly
spaced, opaque lines are scribed~-- several thousand to the centimeter~-- uses interference and
dif\/fraction to measure wave length.
A brief, simple, geometric examination of what is happening to the light, considered as transverse waves,
during interference and dif\/fraction shows how this measurement can be made.

In 1912, von Laue proposed that crystals might be used as natural dif\/fraction gratings for dif\/fracting
high frequency $X$-rays.
The spaces between the atomic planes would act as slits of the grating for dif\/fracting $X$-rays.
The nature of~$X$-rays was far from clear at that point.
Von Laue was convinced by experiments a~year earlier of C.G.~Barkla that $X$-rays are electromagnetic
waves, but, waves of {\it very} short wave lengths.
In order for interference ef\/fects to produce fringe patterns with a~grating, the distance between the
``slits'' cannot be much larger than the wave lengths involved.
Presumably the atoms in the crystals von Laue envisioned are close and symmetrically spaced.
The spaces between the atomic planes would act as slits of the grating.
As it turns out, the distance between neighboring atoms in a~crystal is about 1~\AA.
Von Laue suggested to W.~Friedrich and F.~Knipping that they examine his ideas and calculations
experimentally.
They did, and conf\/irmed his conjectures.

Einstein's 1905 description of the ``photo electric ef\/fect'' is one of the basic known instances of these
early ``ad hoc'' quantum assumptions.
J.J.~Thompson and Lenard noted that ultraviolet light falling on metals causes an emission of electrons.
Varying the intensity of the light does not change the velocity, but does change the number of electrons
emitted.
Einstein used Planck's assumption that energy in radiation is emitted and absorbed in quanta of size
$h\nu$, where $\nu$ is the frequency.
Einstein pictures the light as waves in which energy is distributed discretely over the wave front in
quanta (called {\it photons}) with energy $h\nu$ and momenta $h/\lambda$.
His {\it photo electric equation}
\begin{gather*}
\frac{1}{2}mv_m^2=h\nu-a
\end{gather*}
expresses the maximum kinetic energy of an emitted electron when the frequency of the incident radiation is
$\nu$ and $a$ is the energy required to remove one of the lightly bound electrons ($a$ varies with the
metal).
The photo electric ef\/fect is an indication of the corpuscular (material-particle) nature of light.

Perhaps, the most dramatic instance of the early appearances of the {\it ad hoc} quantum assumptions was
Niels Bohr's 1913 explanation, in theoretical terms, of the lines in the visible portion of the spectrum of
hydrogen, the ``Balmer series''.
Their wave lengths are: 6563 \AA\ (red), 4861~\AA\ (blue), 4380~\AA, 4102~\AA, and 3921~\AA\ (at the ultra
violet end).
Bohr uses Rutherford's ``planetary'' model of the atom as negatively charged electrons moving with uniform
angular velocities in circular orbits around a~central nucleus containing positively charged protons under
an attractive Coulomb force.
In the case of hydrogen, there is one electron and one proton with charges $-\epsilon$ and $+\epsilon$; so,
the attractive force is $-\frac{\epsilon^2}{r^2}$ (between the electron and the proton).
If~$\omega$ is the uniform angular velocity of the electron, its linear velocity is~$r\omega$ (tangential
to the orbit), where $r$ is the radius of its orbit, and its linear acceleration is~$r\omega^2$ directed
``inward'', along ``its radius''.
The moment of inertia $I$ of the electron about the nucleus is $mr^2$ (the measure of its tendency to
resist change in its rotational motion, as mass is the measure of its tendency to resist change in its
linear motion).
The ``angular momentum'' of the electron is~$I\omega$, where $m$ is the (rest) mass of the electron
($9.11\times 10^{-28}$~gm).
Bohr's single quantum assumption is that the angular momentum in its stable orbits should be an {\it
integral} multiple of~$\hbar$ (${=}\frac h{2\pi}$).
That is, $mr^2\omega=k\hbar$, with $k$ an integer, for those $r$ in which the electron occupies a~possible
orbit.

At this point, it is worth moving ahead ten years in the chronological development, to note de Broglie's
1923 synthesis of the increasing evidence of the dual nature of waves and particles; he introduces ``matter
waves''.
De Broglie hypothesized that particles of small mass $m$ moving with (linear) speed~$v$ would exhibit
a~wave like character with wave length~$h/mv$.
Compare this with Einstein's assumption of momentum $h/\lambda$ (${=}mv$).
So, for perspective, an electron moving at $c/3$ would have wave length:
\begin{gather*}
\frac{h}{mv}=\frac{6.625\times10^{-27~}{\rm erg~sec}}{9.11\times10^{-28}~{\rm gm}\times10^{10}~{\rm cm/sec}}
=\frac{66.25\times10^{-10}~{\rm dyne~cm}}{9.11~{\rm gm~cm/sec^2}}\approx0.0727~\text{\AA}. 
\end{gather*}
Returning to the Bohr atom, Bohr's quantum assumption, $mr^2\omega=k\hbar$, can be rewritten as
\begin{gather*}
2\pi r=k\frac{h}{mr\omega}.
\end{gather*}
Combining this with de Broglie's principle, and noting that $r\omega$ is the linear speed of the electron
(directed tangentially to its orbit), $h/mr\omega$ is its wave length, when it is viewed as a~wave.
Of course, $2\pi r$ is the length of its (stable) orbit.
It is intuitively satisfying that the stable orbits are those with radii such that they accommodate an {\it
integral} number of ``complete'' wave cycles, a~``standing wave-train''.

Considering, again, the hydrogen atom, and choosing units so that no constant of proportionality is needed
(the charge $\epsilon$ is in esu~-- electrostatic units), we have that $mr\omega^2=\frac{\epsilon^2}{r^2}$.
From Bohr's quantum assumption, $mr^2\omega=\frac{kh}{2\pi}$.
Thus
\begin{gather*}
m^2r^4\omega^2=\frac{k^2h^2}{4\pi^2}=mr\big(mr^3\omega^2\big)=mr\epsilon^2
\qquad
\text{and}
\qquad
r=\frac{k^2h^2}{4\pi^2m\epsilon^2}.
\end{gather*}
The values 1, 2, 3 of~$k$ give possible values of~$r$ for the stable states.

The kinetic energy of the electron in the orbit corresponding to $r$ is $\frac12 mv^2$ (${=}\frac12
mr^2\omega^2$).
The potential energy of the electron in this Coulomb f\/ield can be taken as the work done in bringing it
from $\infty$ to its orbit of radius~$r$. That is, its potential energy is
\begin{gather*}
-\frac{\epsilon^2}r=-\frac{\epsilon^2}x\bigg|_\infty^r=\int_\infty^r\frac{\epsilon^2}{x^2}{\rm d}x.
\end{gather*}
The total energy is
\begin{gather*}
\frac12mr^2\omega^2-\frac{\epsilon^2}r
=\frac1{2mr^2}\big(m^2r^4\omega^2\big)-\frac{\epsilon^2}r=\frac{\epsilon^2}
{2r}-\frac{\epsilon^2}r=-\frac{\epsilon^2}{2r}=-\frac{2\pi^2m\epsilon^4}{k^2h^2}.
\end{gather*}
The dif\/ferences in energy levels will be
\begin{gather*}
\frac{2\pi^2m\epsilon^4}{h^2}\left(\frac1{k^2}-\frac1{l^2}\right)=h\nu=\frac{hc}{\lambda}.
\end{gather*}
The wave number (that is, the number of waves per centimeter) is given by
\begin{gather*}
w=\frac1\lambda=\frac{2\pi^2m\epsilon^4}{h^3c}\left(\frac1{k^2}-\frac1{l^2}\right)
\\
\phantom{w}
=\frac{2\pi^2\times9.11\times10^{-28}~{\rm gm}\times(4.8025\times10^{-10}~{\rm esu})^4}
{(6.625\times10^{-27}~{\rm erg~sec})^3\times2.99776\times10^{10}~{\rm cm/sec}}(=109,739.53/{\rm cm}
)\times\left(\frac1{k^2}-\frac1{l^2}\right).
\end{gather*}
If we substitute 2 for $k$ and then: 3 for $l$, we f\/ind that $\lambda=6561$~\AA, 4 for $l$ gives
$4860$~\AA, 5 for $l$ gives $4339$~\AA, 6 for $l$ gives $4101$~\AA, 7 for $l$ gives $3969$~\AA.

Comparing these wave lengths with those noted before (from spectroscopy), we see startling agreement,
especially when we note that the physical constants that we use are approximations derived from experiments.
Of course, poor choices for the rest mass and charge of the electron would produce unacceptable values for
wave lengths in the hydrogen spectrum.
It was a~happy circumstance that reasonably accurate values of the mass and charge of the electron were
available after the 1912 Millikan ``oil drop'' experiment.

This striking evidence of the ef\/f\/icacy of uniting the principles of Newtonian (Hamiltonian) mechanics
and ``{\it ad hoc}'' quantum assumptions makes clear the importance of f\/inding a~ma\-the\-matical model
capable of holding, comfortably, within its structure both classical mechanical principles and
a~mathematics that permits the formulation of those ``{\it ad hoc}'' quantum assumptions.
We study the proposal (largely Dirac's) for such a~mathematical structure in the section that follows.

\section{Quantum mechanics~-- a~mathematical model}

\looseness=-1
In Dirac's treatment of physical systems~\cite{D}, there are two basic constituents: the family of
\textit{observables} and the family of \textit{states} in which the system can be found.
In classical (Newtonian--Hamiltonian) mechanics, the observables are algebraic combinations of the
(canonical) coordinates and (conjugate) momenta.
Each state is described by an assignment of numbers to these observables~-- the values certain to be found
by measuring the observables in the given state.
The totality of numbers associated with a~given observable is its spectrum.
In this view of classical statics, the observables are represented as \textit{functions} on the space of
states~-- they form an algebra, necessarily \textit{commutative}, relative to pointwise operations.
The experiments involving atomic and sub-atomic phenomena made it clear that this Newtonian view of
mechanics would not suf\/f\/ice for their basic theory.
Speculation on the meaning of these experimental results eventually led to the conclusion that the only
physically meaningful description of a~state was in terms of an assignment of probability measures to the
spectra of the observables (a measurement of the observable with the system in a~given state will produce
a~value in a~given portion of the spectrum with a~specif\/ic probability).
Moreover, it was necessary to assume that a~state that assigns a~def\/inite value to one observable assigns
a~dispersed measure to the spectrum of some other observable~-- the amount of dispersion involving the
experimentally reappearing Planck's constant.
So, in quantum mechanics, it is not possible to describe states in which a~particle has both a~def\/inite
position and a~def\/inite momentum.
The more precise the position, the less precise the momentum.
This is the celebrated \textit{Heisenberg uncertainty principle}~\cite{H}.
It entails the \textit{non-commutativity} of the algebra of observables.

The search for a~mathematical model that could mirror the structural features of this system and in which
computations in accord with experimental results could be made produced the self-adjoint operators
(possibly unbounded) on a~Hilbert space as the observables and the unit vectors (up to a~complex multiple
of modulus 1) as corresponding to the states~\cite{KI}.
If~$A$ is an observable and $x$ corresponds to a~state of interest, $\langle Ax,x\rangle$, the inner product of
the two vectors~$Ax$ and~$x$, is the real number we get by taking the average of many measurements of~$A$
with the system in the state corresponding to~$x$.
Each such measurement yields a~real number in the spectrum of~$A$.
The probability that that measurement will lie in a~given subset of the spectrum is the measure of that
set, using the probability measure that the state assigns to $A$.
The ``expectation'' of the observable~$A$ in the state corresponding to~$x$ is~$\langle Ax,x\rangle$.

With this part of the model in place, Dirac assigns a~self-adjoint operator $H$ as the energy observable
and, by analogy with classical mechanics, assumes that it will ``generate'' the dynamics, the
time-evolution of the system.
This time-evolution can be described in two ways, either as the states evolving in time, the
``Schr\"{o}dinger picture'' of quantum mechanics, or the observables evolving in time, the ``Heisenberg
picture'' of quantum mechanics.
The prescription for each of these pictures is given in terms of the one-parameter unitary group $t\to
U_t$, where $t\in \mathbb{R}$, the additive group of real numbers, and~$U_t$ is the unitary operator
$\exp(itH)$, formed by applying the spectral-theoretic, function-calculus to the self-adjoint operator~$H$,
the Hamiltonian of our system.
If the initial state of our system corresponds to the unit vector~$x$, then at time~$t$, the system will
have evolved to the state corresponding to the unit vector $U_tx$.
If the observable corresponds to the self-adjoint operator~$A$ at time~0, at time~$t$, it will have evolved
to $U_t^*AU_t$ $({=}\alpha_t(A))$, where, as can be seen easily, $t\to\alpha_t$ is a~one-parameter group of
automorphisms of the ``algebra'' (perhaps, ``Jordan algebra'') of observables.
In any event, the numbers we hope to measure are $\langle AU_tx,U_tx\rangle$, the expectation of the observable
$A$ in the state (corresponding to) $U_tx$, as $t$ varies, and/or $\langle(U_t^*AU_t)x,x\rangle$, the
expectation of the observable $\alpha_t(A)$ in the state $x$, as $t$ varies.
Of course, the two varying expectations are the same, which explains why Heisenberg's ``matrix mechanics''
and Schr\"{o}dinger's ``wave mechanics'' gave the same results.
(In Schr\"{o}dinger's picture, $x$ is a~vector in the Hilbert space viewed as $L_2(\mathbb{R}^3)$, so that
$x$ is a~function, the ``wave function'' of the state, evolving in time as~$U_tx$, while in Heisenberg's
picture, the ``matrix'' coordinates of the operator $A$ evolves in time as~$\alpha_t(A)$.)

\looseness=1
The development of modern quantum mechanics in the mid-1920s was an important motivation for the great
interest in the study of operator algebras in general and von Neumann algebras in particular.
In 
 \cite{VNI} 
 von Neumann def\/ines a~class of algebras of bounded operators on a~Hilbert space that have acquired
the name ``von Neumann algebras''~\cite{Di} (Von Neumann refers to them as ``rings of operators''). Such algebras are \textit{self-adjoint},
\textit{strong-operator closed}, and contain the identity operator.
Von Neumann's article~\cite{VNI} opens up the subject of ``operator algebras'' (see also~\cite{MVI,MVII,MVIII,VNII}).

We use~\cite{KII,KRI,KRII,KRIII,KRIV} as our basic references for results in the theory of operator algebras as well as for much
of our notation and terminology.
Let ${\mathcal H}$ be a~Hilbert space over the complex numbers $\mathbb{C}$ and let $\langle\;,\;\rangle $
denote the (positive def\/inite) inner product on~${\mathcal H}$.
By def\/inition, ${\mathcal H}$~is \textit{complete} relative to the norm $\|\;\|$ def\/ined by the
equation $\|x\|=\langle x,x \rangle^{\frac{1}{2}}$ $(x\in {\mathcal H})$.
If~${\mathcal K}$ is another Hilbert space and $T$ is a~linear operator (or linear transformation) from
${\mathcal H}$ into~${\mathcal K}$, $T$~is continuous if and only if~$\sup \{\|Tx\|: x\in {\mathcal H},
\|x\| \leq 1\}<\infty$.
This supremum is referred to as the \textit{norm} or (\textit{operator}) \textit{bound} of~$T$.
Since continuity is equivalent to the existence of a~f\/inite bound, continuous linear operators are often
described as \textit{bounded} linear operators.
The family ${\mathcal B}({\mathcal H}, {\mathcal K})$ of all bounded linear operators from ${\mathcal H}$
into ${\mathcal K}$ is a~Banach space relative to the operator norm.
When ${\mathcal K}={\mathcal H}$, We write ${\mathcal B}({\mathcal H})$ in place of~${\mathcal B}({\mathcal
H},{\mathcal H})$.
In this case, ${\mathcal B}({\mathcal H})$ is a~Banach algebra with the operator $I$, the identity mapping
on ${\mathcal H}$, as a~unit element.

If~$T$ is in ${\mathcal B}({\mathcal H}, {\mathcal K})$, there is a~unique element $T^*$ of~${\mathcal
B}({\mathcal K}, {\mathcal H})$ such that $\langle Tx,y \rangle=\langle x,T^*y\rangle$
$(x\in{\mathcal H}$, $y\in{\mathcal K})$.
We refer to $T^*$ as the \textit{adjoint} of~$T$.
Moreover, $(aT+bS)^*=\bar{a}T^*+\bar{b}S^*$, $(T^*)^*=T$, $\|T^*T\|=\|T\|^2$, and $\|T\|=\|T^*\|$, whenever
$S,T\in {\mathcal B}({\mathcal H}, {\mathcal K}), and\;a,b\in \mathbb{C}$.
When ${\mathcal H}={\mathcal K}$, we have that $(TS)^*=S^*T^*$.
In this same case, we say that $T$ is \textit{self-adjoint} when $T=T^*$.
A~subset of~${\mathcal B}({\mathcal H})$ is said to be \textit{self-adjoint} if it contains $T^*$ when it
contains~$T$.

\looseness=-1
The metric on ${\mathcal B}({\mathcal H})$ that assigns $\|T-S\|$ as the distance between $T$ and $S$ gives
rise to the \textit{norm} or \textit{uniform topology} on ${\mathcal B}({\mathcal H})$.
There are topologies on ${\mathcal B}({\mathcal H})$ that are weaker than the norm topology.
The \textit{strong-operator topology} is the weakest topology on ${\mathcal B}({\mathcal H})$ such that the
mapping $T\rightarrow Tx$ is continuous for each vector $x$ in ${\mathcal H}$.
The \textit{weak-operator topology} on ${\mathcal B}({\mathcal H})$ is the weakest topology on ${\mathcal
B}({\mathcal H})$ such that the mapping $T\rightarrow \langle Tx,x\rangle$ is continuous for each vector~$x$ in~${\mathcal H}$.

The self-adjoint subalgebras of~${\mathcal B}({\mathcal H})$ containing $I$ that are closed in the norm
topology are known as $C^*$-algebras.
Each abelian $C^*$-algebra is isomorphic to the algebra $C(X)$ (under pointwise addition and
multiplication) of all complex-valued continuous functions on a~compact Hausdorf\/f space $X$.
Each $C(X)$ is isomorphic to some abelian $C^*$-algebra.
The identif\/ication of the family of abelian $C^*$-algebras with the family of function algebras $C(X)$
underlies the interpretation of the general study of~$C^*$-algebras as noncommutative (real) analysis.
This ``noncommutative'' view guides the research and provides a~large template for the motivation of the
subject.
When noncommutative analysis is the appropriate analysis, as in quantum theory~\cite{M,VNIII}, operator algebras provide
the mathematical framework.

\looseness=-1
Those self-adjoint operator algebras that are closed under the strong-operator topology are called von
Neumann algebras (each von Neumann algebra is a~$C^*$-algebra).
We describe some examples of commutative von Neumann algebras.
Suppose $(S, \mu)$ is a~$\sigma$-f\/inite measure space.
Let ${\mathcal H}$ be $L_2(S, \mu)$.
With $f$ an essentially bounded measurable function on $S$, we def\/ine $M_f(g)$ to be the product $f\cdot
g$ for each $g$ in ${\mathcal H}$.
The family $\mathcal{A}=\{M_f\}$ of these {\it multiplication operators} is an abelian von Neumann algebra
and it is referred to as {\it the multiplication algebra} of the measure space $(S, \mu)$.
Moreover, ${\cal A}$ is in no larger abelian subalgebra of~${\mathcal B}({\mathcal H})$.
We say ${\cal A}$ is a~maximal abelian (self-adjoint) subalgebra, a~{\it masa}.
Here are some specif\/ic examples arising from choosing explicit measure spaces.
Choose for $S$ a~f\/inite or countable number of points, say~$n$, each of which has a~positive measure
(each is an {\it atom}).
We write ``${\cal A}={\cal A}_n$'' in this case.
Another example is given by choosing, for $S$, $[0,1]$ with Lebesgue measure.
In this case, we write ``${\cal A}={\cal A}_{\rm c}$'' (``c'' stands for ``continuous'').
Finally, choose, for $S$, $[0,1]$ with Lebesgue measure plus a~f\/inite or countably inf\/inite number~$n$
of atoms.
We write ``${\cal A}={\cal A}_{\rm c}\oplus{\cal A}_n$'' in this case.

\begin{Theorem}
Each abelian von Neumann algebra on a~separable Hilbert space is isomorphic to one of~${\cal A}_n$, ${\cal
A}_{\rm c}$, or ${\cal A}_{\rm c}\oplus{\cal A}_n$.
Each maximal abelian von Neumann algebra on a~separable Hilbert space is unitarily equivalent to one of these.
\end{Theorem}

In the early chapters of~\cite{D}, Dirac is pointing out that Hilbert spaces and their orthonormal bases,
if chosen carefully, can be used to simplify calculations and for determinations of probabilities, for
example, f\/inding the frequencies of the spectral lines in the visible range of the hydrogen atom (the
Balmer series), that is, the spectrum of the operator corresponding to the energy ``observable'' of the
system, the Hamiltonian.
In mathematical terms, Dirac is noting that bases, carefully chosen, will simultaneously ``diagonalize''
self-adjoint operators in an abelian (or ``commuting'') family.
Notably, the masas we have just been describing.

The early experimental work that led to quantum mechanics made it clear that, when dealing with systems at
the atomic scale, where the measurement process interferes with what is being measured, we are forced to
model the physics of such systems at a~single instant of time, as an algebraic mathematical structure that
is not commutative.
Dirac thinks of his small, physical system as an algebraically structured family of ``observables''~--
elements of the system to be observed when studying the system, for example, the position of a~particle in
the system would be an observable $Q$ (a ``canonical coordinate'') and the (conjugate) momentum of that
particle as another observable $P$~-- and they are independent of time.
As the particle moves under the ``dynamics'' of the system, the position $Q$ and momentum $P$ become time
dependent.
By analogy with classical mechanics, Dirac refers to them, in this case, as ``dynamical variables''.
He recalls the Hamilton equation of motion for a~general dynamical variable that is a~function of the
canonical coordinates $\{q_r\}$ and their conjugate momenta $\{p_r\}$:
\begin{gather*}
\frac{{\rm d}q_r}{{\rm d}t}=\frac{\partial H}{\partial p_r},
\qquad
\frac{{\rm d}p_r}{{\rm d}t}=-\frac{\partial H}{\partial q_r},
\end{gather*}
where $H$ is the energy expressed as a~function of the~$q_r$ and~$p_r$ and, possibly, of~$t$.
This $H$ is the Hamiltonian of the system.
Hence, with~$v$ a~dynamical variable that is a~function of the~$q_r$ and~$p_r$, but not explicitly of~$t$,
\begin{gather*}
\frac{{\rm d}v}{{\rm d}t}=\sum_r\left(\frac{\partial v}{\partial q_r}\frac{{\rm d}q_r}{{\rm d}t}+\frac{\partial v}{\partial p_r}
\frac{{\rm d}p_r}{{\rm d}t}\right)
=\sum_r\left(\frac{\partial v}{\partial q_r}\frac{\partial H}{\partial p_r}
-\frac{\partial v}{\partial p_r}\frac{\partial H}{\partial q_r}\right)=[v,H],
\end{gather*}
where $[v,H]$ is the classical {\it Poisson bracket} of~$v$ and $H$.
Dirac is using Lagrange's idea of introducing canonical coordinates and their conjugate momenta, in terms
of which the dynamical variables of interest for a~given system may be expressed, even though those $q_r$
and $p_r$ may not be associated with actual particles in the system.
Noting the fundamental nature of the Poisson bracket in classical mechanics, and establishing its Lie
bracket properties, Dirac def\/ines a~quantum Poisson bracket $[u,v]$ by analogy with the classical bracket.
So, it must be ``real''.
Dirac then argues ``quasi'' mathematically, to show that $uv-vu$ must be $i\hbar[u,v]$, where the real
constant $\hbar$ has to be set by the basic quantum mechanical experiments (giving $\hbar=\frac{h}{2\pi}$,
with $h$ Planck's constant).
Again using classical analogy, the classical coordinates and their conjugate momenta have Poisson brackets
\begin{gather*}
[q_r,q_s]=[p_r,p_s]=0,
\qquad
[q_r,p_s]=\delta_{r,s},
\end{gather*}
where $\delta_{r,s}$ is the Kronecker delta, 1 when $r=s$ and 0 otherwise.
So, Dirac assumes that the quantum Poisson brackets of the position $Q$s and the momentum $P$s satisfy
these same relations.
In the case of one degree of freedom, that is, one $Q$ (and its conjugate momentum $P$), $QP-PQ=i\hbar I$,
the basic Heisenberg relation.
This relation encodes the non-commutativity needed to produce the so-called ``{\it ad hoc} quantum
assumptions'' made by the early workers in quantum physics.
At the same time, this relation gives us a~``numerical grip'' on ``uncertainty'' and ``indeterminacy'' in
quantum mechanics.
In addition, the Heisenberg relation makes it clear (regrettably) that quantum mechanics cannot be modeled
using f\/inite matrices alone.
The trace of~$QP-PQ$ is 0 when $Q$ and $P$ are such matrices, while the trace of~$i\hbar I$ is not 0 (no
matter how we normalize the trace).
It can be shown that the Heisenberg relation cannot be satisf\/ied even with bounded operators on an
inf\/inite-dimensional Hilbert space.
Unbounded operators are needed, even unavoidable for ``representing'' (that is, ``modeling'') the Heisenberg
relation mathematically.
This topic is studied in the following sections.

\section{Basics of unbounded operators on a~Hilbert space}

\looseness=-1
What follows is a~compendium of material drawn from Sections 2.7, 5.2, 5.6, and 6.1 of~\cite{KRI,KRII}:
material that we need in the succeeding sections gathered together here for the convenience of the reader.

\subsection{Definitions and facts}

Let $T$ be a~linear mapping, with domain ${\mathscr{D}}(T)$ a~linear submanifold (not necessarily closed),
of the Hilbert space ${\mathcal H}$ into the Hilbert space ${\mathcal K}$.
We associate a~\textit{graph} ${\mathscr{G}}(T)$ with $T$, where ${\mathscr{G}}(T)=\{(x,Tx):x\in
{\mathscr{D}}(T)\}$.
We say that $T$ is \textit{closed} when ${\mathscr{G}}(T)$ is closed.
The \textit{closed graph theorem} tells us that if~$T$ is def\/ined on all of~${\mathcal H}$, then
${\mathscr{G}}(T)$ is closed if and only if~$T$ is bounded.
The unbounded operators $T$ we consider will usually be \textit{densely defined}, that is,
${\mathscr{D}}(T)$ is dense in ${\mathcal H}$.
We say that $T_0$ \textit{extends} (or is \textit{an extension of}) $T$, and write $T\subseteq T_0$, when
${\mathscr{D}}(T)\subseteq {\mathscr{D}}(T_0)$ and $T_0x=Tx$ for each $x$ in ${\mathscr{D}}(T)$.
If~${\mathscr{G}}(T)^-$, the closure of the graph of~$T$ (a~linear subspace of~${\mathcal H}\bigoplus
{\mathcal K}$), is the graph of a~linear transformation $\overline{T}$, clearly $\overline{T}$ is the
``smallest'' closed extension of~$T$, we say that $T$ is \textit{preclosed} (or~\textit{closable}) and
refer to $\overline{T}$ as the \textit{closure} of~$T$.
From the point of view of calculations with an unbounded operator $T$, it is often much easier to study its
restriction $T|{\mathscr{D}}_0$ to a~dense linear manifold ${\mathscr{D}}_0$ in its domain
${\mathscr{D}}(T)$ than to study $T$ itself.
If~$T$ is closed and ${\mathscr{G}}(T|{\mathscr{D}}_0)^-={\mathscr{G}}(T)$, we say that ${\mathscr{D}}_0$
is a~\textit{core} for $T$.
Each dense linear manifold in ${\mathscr{G}}(T)$ corresponds to a~core for $T$.
\begin{Definition}
{\rm If~$T$ is a~linear transformation with ${\mathscr{D}}(T)$ dense in the Hilbert space ${\mathcal H}$
and range contained in the Hilbert space ${\mathcal K}$, we def\/ine a~mapping $T^*$, the \textit{adjoint}
of~$T$, as follows.
Its domain consists of those vectors $y$ in ${\mathcal K}$ such that, for some vector $z$ in ${\mathcal
H}$, $\langle x,z\rangle=\langle Tx,y\rangle$ for all $x$ in ${\mathscr{D}}(T)$.
For such $y$, $T^*y$ is $z$.
If~$T=T^*$, we say that $T$ is \textit{self-adjoint}.
(Note that the formal relation $\langle Tx, y\rangle=\langle x,T^*y\rangle$, familiar from the case of bounded
operators, remains valid in the present context only when $x\in {\mathscr{D}}(T)$ and $y\in
{\mathscr{D}}(T^*)$.)}
\end{Definition}
\begin{Remark}
\label{selfadjoint}
{\rm If~$T$ is densely def\/ined, then $T^*$ is a~closed linear operator.
If~$T_0$ is an extension of~$T$, then $T^*$ is an extension of~$T_0^*$.}
\end{Remark}
\begin{Theorem}\label{thm: preclosed}
If~$T$ is a~densely defined linear transformation from the Hilbert space ${\mathcal H}$ to the Hilbert
space ${\mathcal K}$, then
\begin{enumerate}\itemsep=0pt
\item[$(i)$] if~$T$ is preclosed, $(\overline{T})^*=T^*$;

\item[$(ii)$] $T$ is preclosed if and only if~${\mathscr{D}}(T^*)$ is dense in ${\mathcal K}$;

\item[$(iii)$] if~$T$ is preclosed, $T^{**}=\overline{T}$;

\item[$(iv)$] if~$T$ is closed, $T^*T+I$ is one-to-one with range ${\mathcal H}$ and positive inverse of bound
not exceeding~$1$.
\end{enumerate}
\end{Theorem}
\begin{Definition}
 We say that $T$ is \textit{symmetric} when ${\mathscr{D}}(T)$ is dense in ${\mathcal H}$ and $\langle Tx,
y\rangle=\langle x,Ty\rangle$ for all $x$ and $y$ in ${\mathscr{D}}(T)$.
Equivalently, $T$ is symmetric when $T\subseteq T^*$.
(Since $T^*$ is closed and ${\mathscr{G}}(T)\subseteq {\mathscr{G}}(T^*)$, in this case, $T$ is preclosed
if it is symmetric.
If~$T$ is self-adjoint, $T$ is both symmetric and closed.)

\end{Definition}
\begin{Remark}
\label{remark: maximal symmetric}
{\rm If~$A\subseteq T$ with $A$ self-adjoint and $T$ symmetric, then $A\subseteq T\subseteq T^*$, so that
$T^*\subseteq A^*=A\subseteq T\subseteq T^*$ and $A=T$.
It follows that $A$ has no proper symmetric extension.
That is, a~self-adjoint operator is \textit{maximal symmetric}.}
\end{Remark}
\begin{Proposition}
\label{prop:2.1.6}
If~$T$ is a~closed symmetric operator on the Hilbert space ${\mathcal H}$, the following assertions are
equivalent:
\begin{enumerate}\itemsep=0pt
\item[$(i)$] $T$ is self-adjoint;

\item[$(ii)$] $T^*\pm iI$ have $(0)$ as null space;

\item[$(iii)$] $T\pm iI$ have ${\mathcal H}$ as range;

\item[$(iv)$] $T\pm iI$ have ranges dense in ${\mathcal H}$.
\end{enumerate}
\end{Proposition}
\begin{Proposition}
\label{prop: RN}
If~$T$ is a~closed linear operator with domain dense in a~Hilbert space ${\mathcal H}$ and with range in
${\mathcal H}$, then
\begin{gather*}
R(T)=I-N\big(T^*\big),
\qquad\!
N(T)=I-R\big(T^*\big),
\qquad\!
R\big(T^*T\big)=R(T^*),
\qquad\!
N\big(T^*T\big)=N(T),
\end{gather*}
where $N(T)$ and $R(T)$ denote the projections whose ranges are, respectively, the null space of~$T$ and
the closure of the range of~$T$.
\end{Proposition}

\subsection{Spectral theory}

If~$A$ is a~bounded self-adjoint operator acting on a~Hilbert space ${\mathcal H}$ and $\mathscr{A}$ is an
abelian von Neumann algebra containing $A$, there is a~family $\{E_\lambda\}$ of projections in
$\mathscr{A}$ (indexed by ${\mathbb{R}}$), called the spectral resolution of~$A$, such that
\begin{enumerate}\itemsep=0pt
\item[$(i)$]   $E_\lambda=0$ if~$\lambda<-\|A\|$, and $E_\lambda=I$ if~$\|A\|\leq\lambda$;

\item[$(ii)$] $E_\lambda\leq E_{\lambda'}$ if~$\lambda\leq\lambda'$;

\item[$(iii)$] $E_\lambda=\wedge_{\lambda'>\lambda}E_{\lambda'}$;

\item[$(iv)$] $AE_\lambda\leq\lambda E_\lambda$ and $\lambda(I-E_\lambda)\leq A(I-E_\lambda)$ for each
$\lambda$;

\item[$(v)$] $A=\int_{-\|A\|}^{\|A\|}\lambda {\rm d} E_\lambda$ in the sense of norm convergence of
approximating Riemann sums; and $A$ is the norm limit of f\/inite linear combinations with coef\/f\/icients
in ${\rm sp} (A)$, the spectrum of~$A$, of orthogonal projections $E_{\lambda ^\prime} - E_\lambda$.
\end{enumerate}

\noindent
  $\{E_\lambda\}$ is said to be a~{\it resolution of the identity} if~$\{E_\lambda\}$ satisf\/ies
$(ii)$, $(iii)$, $\wedge_{\lambda\in {\mathbb{R}}}E_\lambda=0$ and $\vee_{\lambda\in {\mathbb{R}}}E_\lambda=I$.
With the abelian von Neumann algebra $\mathscr{A}$ isomorphic to $C(X)$ and $X$ an extremely disconnected
compact Hausdorf\/f space, if~$f$ and $e_\lambda$ in $C(X)$ correspond to $A$ and $E_\lambda$ in
$\mathscr{A}$, then~$e_\lambda$ is the characteristic function of the largest clopen subset $X_\lambda$ on
which $f$ takes values not exceeding $\lambda$.

The spectral theory described above can be extended to unbounded self-adjoint operators.
We associate an unbounded spectral resolution with each of them.
We begin with a~discussion that details the relation between unbounded self-adjoint operators and the
multiplication algebra of a~measure space.

If~$g$ is a~complex measurable function (f\/inite almost everywhere) on a~measure space~$(S, m)$, without
the restriction that it be essentially bounded~-- multiplication by $g$ will not yield an
everywhere-def\/ined operator on $L_2(S)$, for many of the products will not lie in~$L_2(S)$.
Enough functions $f$ will have product~$gf$ in $L_2(S)$, however, to form a~dense linear submanifold
${\mathscr{D}}$ of~$L_2(S)$ and constitute a~(dense) domain for an (unbounded) multiplication operator~$M_g$.
To see this, let $E_n$ be the (bounded) multiplication operator corresponding to the characteristic
function of the (measurable) set on which $|g|\leq n$.
Since $g$ is f\/inite almost everywhere, $\{E_n\}$ is an increasing sequence of projections with union~$I$.
The union ${\mathscr{D}}_0$ of the ranges of the~$E_n$ is a~dense linear manifold of~$L_2(S)$ contained in~${\mathscr{D}}$.
A measure-theoretic argument shows that $M_g$ is closed with ${\mathscr{D}}_0$ as a~core.
In fact, if~$\{f_n\}$ is a~sequence in ${\mathscr{D}}$ converging in $L_2(S)$ to $f$ and $\{gf_n\}$
converges in $L_2(S)$ to $h$, then, passing to subsequences, we may assume that~$\{f_n\}$ and~$\{gf_n\}$
converges almost everywhere to $f$ and $h$, respectively.
But, then, $\{gf_n\}$ converges almost everywhere to~$gf$, so that~$gf$ and~$h$ are equal almost everywhere.
Thus $gf\in L_2(S)$, $f\in {\mathscr{D}}$, $h=M_g(f)$, and~$M_g$ is closed.
With $f_0$ in ${\mathscr{D}}$, $E_nf_0$ converges to $f_0$ and $\{M_gE_nf_0\}=\{E_nM_gf_0\}$ converges to
$M_gf_0$.
Now $E_nf_0\in {\mathscr{D}}_0$, so that ${\mathscr{D}}_0$ is a~core for $M_g$.
Note that~$M_gE_n$ is bounded with norm not exceeding $n$.
One can show that~$M_g$ is an (unbounded) self-adjoint operator when $g$ is real-valued.
If~$M_g$ is unbounded, we cannot expect it to belong to the multiplication algebra~$\mathscr{A}$ of the
measure space $(S, m)$.
Nonetheless, there are various ways in which $M_g$ behaves as if it were in~$\mathscr{A}$~-- for example,
$M_g$ is unchanged when it is ``transformed'' by a~unitary operator~$U$ commuting with~$\mathscr{A}$.
In this case, $U\in \mathscr{A}$, so that $U=M_u$ where~$u$ is a~bounded measurable function on~$S$ with
modulus $1$ almost everywhere.
With $f$ in ${\mathscr{D}}(M_g)$, $guf\in L_2(S)$; while, if~$guh\in L_2(S)$, then $gh\in L_2(S)$ and $h\in
{\mathscr{D}}(M_g)$.
Thus $U$ transforms ${\mathscr{D}}(M_g)$ onto itself.
Moreover
\begin{gather*}
\big(U^*M_gU\big)(f)=\overline{u}guf=|u|^2gf=gf.
\end{gather*}
Thus $U^*M_gU=M_g$.
The fact that $M_g$ ``commutes'' with all unitary operators commuting with~$\mathscr{A}$ in conjunction
with the fact that each element of a~$C$*-algebra is a~f\/inite linear combination of unitary elements in
the algebra and the double commutant theorem (from which it follows that a~bounded operator that commutes
with all unitary operators commuting with~$\mathscr{A}$ lies in~$\mathscr{A}$) provides us with an
indication of the extent to which~$M_g$ ``belongs'' to~$\mathscr{A}$.
We formalize this property in the def\/inition that follows.
\begin{Definition}
\label{df affiliation}
 We say that a~closed densely def\/ined operator $T$ is \textit{affiliated} with a~von Neumann algebra
${\mathcal R}$ and write $T\eta{\mathcal R}$ when $U^*TU=T$ for each unitary operator $U$ commuting with~${\mathcal R}$.
(Note that the equality, $U^*TU=T$, is to be understood in the strict sense that $U^*TU$ and $T$ have the
same domain and formal equality holds for the transforms of vectors in that domain.
As far as the domains are concerned, the ef\/fect is that $U$ transforms ${\mathscr{D}}(T)$ onto itself.)
\end{Definition}

\begin{Remark}
\label{remark: check core}
  If~$T$ is a~closed densely def\/ined operator with core ${\mathscr{D}}_0$ and $U^*TUx=Tx$ for each~$x$
in~${\mathscr{D}}_0$ and each unitary operator $U$ commuting with a~von Neumann algebra~${\mathcal R}$,
then~$T\eta{\mathcal R}$.
\end{Remark}

\begin{Theorem}
If~$A$ is a~self-adjoint operator acting on a~Hilbert space ${\mathcal H}$, $A$ is affiliated with some
abelian von Neumann algebra $\mathscr{A}$.
There is a~resolution of the identity $\{E_\lambda\}$ in ${\mathscr{A}}$ such that $\cup_{n=1}^\infty
F_n({\mathcal H})$ is a~core for $A$, where $F_n=E_n-E_{-n}$, and $Ax=\int_{-n}^n\lambda {\rm d}E_\lambda x$ for
each~$x$ in~$F_n({\mathcal H})$ and all~$n$, in the sense of norm convergence of approximating Riemann sums.
\end{Theorem}
Since $A$ is self-adjoint, from Proposition~\ref{prop:2.1.6}, $A+iI$ and $A-iI$ have range ${\mathcal H}$
and null space~$(0)$; in addition, they have inverses, say $T_+$ and $T_-$, that are everywhere def\/ined
with bound not exceeding~1.
Let ${\mathscr{A}}$ be an abelian von Neumann algebra containing $I$, $T_+$ and $T_-$.
If~$U$ is a~unitary operator in ${\mathscr{A}}'$, for each $x$ in ${\mathscr{D}}(A)$,
$Ux=UT_+(A+iI)x=T_+U(A+iI)x$ so that $(A+iI)Ux=U(A+iI)x$; and $U^{-1}(A+iI)U=A+iI$.
Thus~$U^{-1}AU=A$ and~$A\eta{\mathscr{A}}$.
In particular, $A$ is af\/f\/iliated with the abelian von Neumann algebra generated by~$I$,~$T_+$ and~$T_-$.
Since ${\mathscr{A}}$ is abelian, ${\mathscr{A}}$ is isomorphic to~$C(X)$ with~$X$ an extremely
disconnected compact Hausdorf\/f space.
Let~$g_+$ and~$g_-$ be the functions in $C(X)$ corresponding to~$T_+$ and~$T_-$.
Let~$f_+$ and~$f_-$ be the functions def\/ined as the reciprocals of~$g_+$ and~$g_-$, respectively, at
those points where~$g_+$ and~$g_-$ do not vanish.
Then~$f_+$ and~$f_-$ are continuous where they are def\/ined on~$X$, as is the function~$f$ def\/ined by
$f=(f_++f_-)/2$.
In a~formal sense, $f$ is the function that corresponds to~$A$.
Let $X_\lambda$ be the largest clopen set on which $f$ takes values not exceeding~$\lambda$.
Let $e_\lambda$ be the characteristic function of~$X_\lambda$ and $E_\lambda$ be the projection in
${\mathscr{A}}$ corresponding to~$e_\lambda$.
In this case, $\{E_\lambda\}$ satisf\/ies $E_\lambda\leq E_{\lambda'}$ if~$\lambda\leq\lambda'$,
$E_\lambda=\wedge_{\lambda'>\lambda}E_{\lambda'}$, $\vee_\lambda E_\lambda=I$ and $\wedge_\lambda
E_\lambda=0$.
That is, we have constructed a~resolution of the identity $\{E_\lambda\}$.
This resolution is unbounded if~$f\notin C(X)$.
Let $F_n=E_n-E_{-n}$, the \textit{spectral projection} corresponding to the interval $[-n,n]$ for each
positive integer~$n$.
$AF_n$ is bounded and self-adjoint.
Moreover, $\cup_{n=1}^\infty F_n({\mathcal H})$ is a~core for~$A$.
From the spectral theory of bounded self-adjoint operators, $Ax=\int_{-n}^n\lambda {\rm d}E_\lambda x$, for each
$x$ in $F_n({\mathcal H})$ and all $n$.
If~$x\in{\mathscr{D}}(A)$, $\int_{-n}^n\lambda {\rm d}E_\lambda x=\int_{-n}^n\lambda {\rm d}E_\lambda
F_nx=AF_nx\rightarrow Ax$.
Interpreted as an improper integral, we write $Ax=\int_{-\infty}^\infty\lambda {\rm d}E_\lambda x$
$(x\in {\mathscr{D}}(A))$.

\subsection{Polar decomposition}

Each $T$ in ${\mathcal B}({\mathcal H})$ has a~unique decomposition as $VH$, the \textit{polar
decomposition} of~$T$, where $H=(TT^*)^{1/2}$ and $V$ maps the closure of the range of~$H$, denoted by
r(H), isometrically onto~$r(T)$ and maps the orthogonal complement of~$r(H)$ to $0$.
We say that $V$ is a~\textit{partial isometry} with \textit{initial space} $r(H)$ and \textit{final space}
$r(T)$.
If~$R(H)$ is the projection with range~$r(H)$ (the \textit{range projection} of~$H$), then $V^*V=R(H)$ and
$VV^*=R(T)$.
We note that the components~$V$ and~$H$ of this polar decomposition lie in the von Neumann algebra~${\mathcal R}$ when~$T$ does.
There is an extension of the polar decomposition to the case of a~closed densely def\/ined linear operator
from one Hilbert space to another.

\begin{Theorem}\label{thm: polar}
If~$T$ is a~closed densely defined linear transformation from one Hilbert space to another, there is
a~partial isometry $V$ with initial space the closure of the range of~$(T^*T)^{1/2}$ and final space the
closure of the range of~$T$ such that $T=V(T^*T)^{1/2}=(T^*T)^{1/2}V$.
Restricted to the closures of the ranges of~$T^*$ and $T$, respectively, $T^*T$ and $TT^*$ are unitarily
equivalent $($and~$V$ implements this equivalence$)$.
If~$T=WH$, where $H$ is a~positive operator and $W$ is a~partial isometry with initial space the closure of
the range of~$H$, then $H=(T^*T)^{1/2}$ and~$W=V$.
If~${\mathcal R}$ is a~von Neumann algebra, $T\eta{\mathcal R}$ if and only if~$V\in{\mathcal R}$ and
$(T^*T)^{1/2}\eta{\mathcal R}$.
\end{Theorem}

\section{Representations of the Heisenberg relation}

In this section, we study the Heisenberg relation: classes of elements with which it can't be realized,
a~classic example in which it can be realized with a~bounded and an unbounded operator (the argument drawn
from~\cite{KRIII})
and special information about extendability to self-adjoint operators.
The standard representation, involving a~multiplication operator (``position'') and dif\/ferentiation
(``momentum'') viewed as the inf\/initesimal generator of the one-parameter group of translations of the
additive group of the reals, appears in Section~\ref{s5.3}.
The account is precise and complete also with regard to domains and unbounded operator considerations.

\subsection{Bounded operators}

Heisenberg's encoding of the ad-hoc quantum rules in his commutation relation, $QP-PQ=i\hbar I$, where $Q$
and $P$ are the observables corresponding to the position and momentum (say, of a~particle in the system)
respectively, $I$ is the identity operator and $\hbar=\frac{h}{2\pi}$ with $h$ as Planck's constant,
embodies the characteristic \textit{indeterminacy} and \textit{uncertainty} of quantum theory.
The very essence of the relation is its introduction of non-commutativity between the particle's position
$Q$ and its corresponding conjugate momentum $P$.
This is the basis for the view of quantum physics as employing noncommutative mathematics, while classical
(Newtonian--Hamiltonian) physics involves just commutative mathematics.
If we look for mathematical structures that can accommodate this non-commutativity and permit the necessary
computations, families of matrices come quickly to mind.
Of course, we, and the early quantum physicists, can hope that the f\/inite matrices will suf\/f\/ice for
our computational work in quantum physics.
Unhappily, this is not the case, as the trace (functional) on the algebra of complex $n \times n$ matrices
makes clear to us.
The trace of the left side of the Heisenberg relation is $0$ for matrices~$P$ and~$Q$, while the trace of
the right side is $i\hbar$ (${\not=}0$).
That is to say, the Heisenberg relation cannot be satisf\/ied by f\/inite matrices.
Of course, the natural extension of this attempt is to wonder if inf\/inite-dimensional Hilbert spaces
might not ``support'' such a~representation with bounded operators.
Even this is not possible as we shall show.
\begin{Proposition}
If~$A$ and $B$ are elements of a~Banach algebra $\mathfrak{A}$ with unit $I$, then ${\rm sp}(AB)\cup\{0\}
={\rm sp}(BA)\cup\{0\}$.
\end{Proposition}
\begin{proof}
If~$\lambda\not=0$ and $\lambda\in{\rm sp}(AB)$, then $AB-\lambda I$ and, hence $(\lambda^{-1} A)B-I$ are not
invertible.
On the other hand, if~$\lambda\not \in {\rm sp}(BA)$, then $BA-\lambda I$ and, hence, $B(\lambda^{-1}A)-I$
are invertible.
Our task, then, is to show that $I-AB$ is invertible in $\mathfrak{A}$ if and only if~$I-BA$ is invertible
in $\mathfrak{A}$, for arbitrary elements $A$ and $B$ of~$\mathfrak{A}$.

Let us argue informally for the moment.
The following argument leads us to the correct formula for the inverse of~$I-BA$, and gives us a~proof that
holds in any ring with a~unit.
\begin{gather*}
(I-AB)^{-1}=\sum_{n=0}^\infty(AB)^n=I+AB+ABAB+\cdots
\end{gather*}
and
\begin{gather*}
B(I-AB)^{-1}A=BA+BABA+BABABA+\cdots=(I-BA)^{-1}-I.
\end{gather*}
Thus if~$I-AB$ has an inverse, we may hope that $B(I-AB)^{-1}A+I$ is an inverse to $I-BA$.
Multiplying, we have
\begin{gather*}
(I-BA)[B(I-AB)^{-1}A+I]=B(I-AB)^{-1}A+I-BAB(I-AB)^{-1}A-BA
\\
\qquad=B[(I-AB)^{-1}-AB(I-AB)^{-1}]A+I-BA=I,
\end{gather*}
and similarly for right multiplication by $I-BA$.
\end{proof}

Finally, ${\rm sp}(A+I)=\{1+a:a\in {\rm sp}(A)\}$, together with the proposition, yield the fact that the unit
element $I$ of a~Banach algebra is not the \textit{commutator} $AB-BA$ of two elements $A$ and $B$.
(If~$I=AB-BA$, then ${\rm sp}(AB)= 1+{\rm sp}(BA)$, which is not consistent with ${\rm sp}(AB)\cup\{0\} ={\rm
sp}(BA)\cup\{0\}$.) Therefore, in quantum theory, the commutation relations (in particular, the Heisenberg
relation) are not representable in terms of bounded operators.
(A.~Wintner~\cite{Win} proved the quantum result for bounded self-adjoint operators on a~Hilbert space.
H.~Wielandt~\cite{Wie} proved it for elements of a~Banach algebra by a~method dif\/ferent from what has
just been used.)

\subsection{With unbounded operators}

In Section~5.1, we showed that the Heisenberg relation is not representable in terms of elements of complex Banach
algebras with a~unit element.
Therefore, in our search for ways to represent the Heisenberg relation in some (algebraic) mathematical
structure, we can eliminate f\/inite matrices, bounded operators on an inf\/inite-dimensional Hilbert
space, and even elements of more general complex Banach algebras.
Is there anything left? It becomes clear that unbounded operators would be essential for dealing with the
non-commutativity that the Heisenberg relation carries.
The following example gives a~specif\/ic representation of the relation with one of the representing
operators bounded and the other unbounded.
\begin{Example}
\label{example: L2}
Let ${\mathcal H}$ be the Hilbert space $L_2$, corresponding to Lebesgue measure on the unit interval
$[0,1]$, and let ${\mathscr{D}}_0$ be the subspace consisting of all complex-valued functions $f$ that have
a~continuous derivative $f^\prime$ on $[0,1]$ and satisfy $f(0)=f(1)=0$.
Let $D_0$ be the operator with domain ${\mathscr{D}}_0$ and with range in ${\mathcal H}$ def\/ined by
$D_0f=f^\prime$.
We shall show that $iD_0$ is a~densely def\/ined symmetric operator and that
\begin{gather*}
(iD_0)M-M(iD_0)=iI|{\mathscr{D}}_0,
\end{gather*}
where $M$ is the bounded linear operator def\/ined by $(Mf)(s)=sf(s)$
$(f\in L_2$; $0\leq s\leq 1)$.
\end{Example}

\begin{proof}
Each element $f$ of~${\mathcal H}$ can be approximated (in $L_2$ norm) by a~continuous function~$f_1$.
In turn, $f_1$ can be approximated (in the uniform norm, hence in the~$L_2$ norm) by a~polynomial~$f_2$.
Finally, $f_2$ can be approximated (in $L_2$ norm) by an element $f_3$ of~${\mathscr{D}}_0$; indeed, it
suf\/f\/ices to take $f_3=gf_2$, where $g: [0,1]\rightarrow [0,1]$ is continuously dif\/ferentiable,
vanishes at the endpoint~$0$ and~$1$, and takes the value $1$ except at points very close to $0, 1$.

The preceding argument shows that ${\mathscr{D}}_0$ is dense in ${\mathcal H}$, so $D_0$ is a~densely
def\/ined linear operator.
When $f,g\in {\mathscr{D}}_0$, the function $\bar g$ has a~continuous derivative $\bar g^\prime$, and we
have
\begin{gather*}
\langle D_0f,g\rangle=\int_0^1f^\prime(s)\overline{g(s)}{\rm d}s=\Big[f(s)\overline{g(s)}
\Big]_0^1-\int_0^1f(s)\overline{g^\prime(s)}{\rm d}s
\\
\phantom{\langle D_0f,g\rangle}
=-\int_0^1f(s)\overline{g^\prime(s)}{\rm d}s=-\langle f,D_0g\rangle.
\end{gather*}
Thus $\langle iD_0f,g\rangle=\langle f,iD_0g\rangle$, for all $f$ and $g$ in ${\mathscr{D}}_0$; and $iD_0$ is symmetric.

When $f\in {\mathscr{D}}_0$, $Mf\in {\mathscr{D}}_0$ and
\begin{gather*}
(D_0Mf)(s)=\frac{{\rm d}}{{\rm d}s}\big(sf(s)\big)=f(s)+sf^\prime(s)=f(s)+(MD_0f)(s).
\end{gather*}
Thus $D_0Mf-MD_0f=f$
$(f\in {\mathscr{D}}_0)$.
\end{proof}

One can press this example further to show that $iD_0$ has a~self-adjoint extension.
\begin{Example}
Let ${\mathcal H}$, ${\mathscr{D}}_0$ and $D_0$ be def\/ined as in the preceding example, and let
${\mathcal H}_1=\{f_1\in{\mathcal H}: \langle f_1,u\rangle=0\}$, where $u$ is the unit vector in~${\mathcal H}$
def\/ined by $u(s)=1$ $(0\leq s\leq 1)$.
When $f\in{\mathcal H}$, def\/ine $Kf$ in ${\mathcal H}$ by
\begin{gather*}
(Kf)(s)=\int_0^sf(t){\rm d}t,
\qquad
 0\leq s\leq1 .
\end{gather*}
We shall show the following:
\begin{enumerate}\itemsep=0pt
\item[$(i)$] $K\in {\mathcal B}({\mathcal H})$, K has null space $\{0\}$ and ${\mathscr{D}}_0\subseteq
K({\mathcal H}_1)$.

\item[$(ii)$] The equation
\begin{gather*}
D_1Kf_1=f_1,
\qquad
f_1\in{\mathcal H}_1,
\end{gather*}
def\/ines a~closed linear operator $D_1$ with domain ${\mathscr{D}}_1=K({\mathcal H}_1)$, and $D_1$ is the
closure of~$D_0$.

\item[$(iii)$] The equation
\begin{gather*}
D_2(Kf+au)=fm,
\qquad
 f\in{\mathcal H},\quad a\in\mathbb{C},
\end{gather*}
def\/ines a~closed linear operator $D_2$, with domain ${\mathscr{D}}_2=\{Kf+au:f\in {\mathcal H},\;a\in
\mathbb{C}\}$, that extends~$D_1$.

\item[$(iv)$] Let ${\mathscr{D}}_3=\{Kf_1+au:f_1\in {\mathcal H}_1,\;a\in \mathbb{C}\}$, and let $D_3$ be the
restriction $D_2|{\mathscr{D}}_3$.
 $D_3$~is a~closed densely def\/ined operator and $D_1\subseteq D_3=-D_3^*\subseteq D_2$ so that $iD_3$ is
a~self-adjoint extension of~$iD_0$.
\end{enumerate}
\end{Example}
\begin{proof}
$(i)$ For any unit vector $y$ in ${\mathcal H}$,
\begin{gather*}
\|Ky\|^2=\int_0^1|(Ky)(s)|^2{\rm d}s=\int_0^1\left|\int_0^sy(t){\rm d}t\right|^2{\rm d}s
\leq\int_0^1\left(\int_0^s|y(t)|^2{\rm d}t\right){\rm d}s
\\
\phantom{\|Ky\|^2}
\leq\int_0^1\left(\int_0^1|y(t)|^2{\rm d}t\right){\rm d}s=\int_0^1\|y\|^2{\rm d}s=1.
\end{gather*}
Thus $K\in {\mathcal B}({\mathcal H})$.
If~$f\in {\mathcal H}$ and $Kf=0$, then $\int_0^sf(t){\rm d}t=0$ ($0\leq s\leq 1$), and $f=0$; so $K$ has null
space $\{0\}$.
If~$g\in {\mathscr{D}}_0$, then $g$ has a~continuous derivative $g'$ on $[0,1]$ and $g(0)=g(1)=0$.
Since $g'\in {\mathcal H}$ and
\begin{gather*}
\langle g',u\rangle=\int_0^1g'(s){\rm d}s=g(1)-g(0)=0,
\end{gather*}
it follows that $g'\in{\mathcal H}_1$.
Moreover
\begin{gather*}
(Kg')(s)=\int_0^sg'(t){\rm d}t=g(s)-g(0)=g(s),
\qquad
 0\leq s\leq1,
\end{gather*}
so $g=Kg'\in{\mathcal K}({\mathcal H}_1)$.
Thus ${\mathscr{D}}_0\subseteq K({\mathcal H}_1)$.

$(ii)$ From $(i)$, $K({\mathcal H}_1)$ is dense in ${\mathcal H}$ (since ${\mathscr{D}}_0$ is dense in
${\mathcal H}$).
Now $K$ is one-to-one, the equation $D_1Kf_1=f_1$ $(f_1\in {\mathcal H}_1)$ def\/ines a~linear operator
$D_1$ with dense domain ${\mathscr{D}}_1=(K({\mathcal H}_1))$.

If~$\{g_n\}$ is a~sequence in ${\mathscr{D}}_1$ such that $g_n\rightarrow g$ and $D_1g_n=f$, then
$g_n=Kf_n$ and $D_1g_n=f_n$ for some sequence $\{f_n\}$ in ${\mathcal H}_1$.
Since $f_n\rightarrow f$, ${\mathcal H}_1$ is closed, and $K$ is bounded, we have $f\in{\mathcal H}_1$ and
$Kf=\lim Kf_n=\lim g_n=g$.
Thus $g\in K({\mathcal H}_1)={\mathscr{D}}_1$, and $D_1g=f$; so $D_1$ is closed.

If~$g\in{\mathscr{D}}_0(\subseteq K({\mathcal H}_1))$, then $g=Kg'$ and $g'\in {\mathcal H}_1$.
Thus $g\in{\mathscr{D}}_1$, $D_1g=g'=D_0g$; so $D_0\subseteq D_1$.
Since $D_1$ is closed, $\overline{D}_0\subseteq D_1$.

To prove that $D_1\subseteq \overline{D}_0$, suppose that $g\in{\mathscr{D}}_1$ and $D_1g=f$.
Then $f\in {\mathcal H}_1$, and $Kf=g$.
There is a~sequence $\{h_n\}$ of continuous functions on $[0,1]$ such that $\|f-h_n\|\rightarrow 0$; and
$\langle h_n,u\rangle\rightarrow \langle f,u\rangle=0$.
With $f_n$ def\/ined as $h_n-\langle h_n,u\rangle u$, $f_n$ is continuous, $\langle f_n,u\rangle=0$, and
$\|f-f_n\|\rightarrow 0$.
Let $g_n=Kf_n$, so that $g_n\rightarrow Kf=g$.
Since 
\begin{gather*}
g_n(s)=\int_0^sf_n(t){\rm d}t,
\qquad
\int_0^1f_n(t){\rm d}t=\langle f_n,u\rangle=0,
\end{gather*}
it follows that $g_n$ has a~continuous derivative $f_n$, and satisf\/ies $g_n(0)=g_n(1)=0$.
Thus $g_n\in{\mathscr{D}}_0$, $g_n\rightarrow g$, and $D_0g_n=f_n\rightarrow f=D_1g$.
This shows that each point $(g,D_1g)$ in the graph of~$D_1$ is the limit of a~sequence $\{(g_n,D_0g_n)\}$
in the graph of~$D_0$; so $D_1\subseteq \overline{D}_0$.

$(iii)$ If~$f\in {\mathcal H}$, $a\in \mathbb{C}$ and $Kf+au=0$, then
\begin{gather*}
a+\int_0^sf(t){\rm d}t=au(s)+(Kf)(s)=0
\end{gather*}
for almost all $s$ in $[0,1]$ and hence, by continuity, for all $s$ in $[0,1]$.
With $s=0$, we obtain $a=0$; it follows that $f$ is a~null function.
So the equation
\begin{gather*}
D_2(Kf+au)=f,
\qquad
f\in{\mathcal H},\quad a\in\mathbb{C},
\end{gather*}
def\/ines a~linear operator $D_2$ with domain ${\mathscr{D}}_2=\{Kf+au:f\in {\mathcal H},\;a\in
\mathbb{C}\}$.
In addition, $D_1\subseteq D_2$.
In particular, $D_2$ is densely def\/ined.

If~$\{g_n\}$ is a~sequence in ${\mathscr{D}}_2$ such that $g_n\rightarrow g$ and $D_2g_n\rightarrow f$,
then $g_n=Kf_n+a_nu$, where $f_n\in{\mathcal H}$ and $a_n\in \mathbb{C}$; and $D_2g_n=f_n$.
Thus
\begin{gather*}
f_n\rightarrow f,
\qquad
Kf_n\rightarrow Kf,
\qquad
a_nu=g_n-Kf_n\rightarrow g-Kf,
\end{gather*}
and therefore $g-Kf=au$ for some scalar $a$.
Thus $g=Kf+au\in{\mathscr{D}}_2$, $D_2g=f$; and $D_2$ is closed.

$(iv)$ Since ${\mathscr{D}}_1\subseteq {\mathscr{D}}_3\subseteq {\mathscr{D}}_2$ and $D_1\subseteq D_2$, it
is evident that $D_3$ $({=}D_2|{\mathscr{D}}_3)$ is densely def\/ined and $D_1\subseteq D_3\subseteq D_2$.
We shall show that $D_3=-D_3^*$.
It follows that $D_3$ is closed and $iD_3$ is self-adjoint.

First, we note that if~$f_1\in{\mathcal H}_1$, $f\in{\mathcal H}$ and $a\in \mathbb{C}$, then
\begin{gather}
\langle Kf_1,f\rangle+\langle f_1,Kf+au\rangle=\langle Kf_1,f\rangle+\langle f_1,Kf\rangle
\nonumber
\\
\qquad=\int_0^1(Kf_1)(s)\overline{f(s)}{\rm d}s+\int_0^1f_1(t)\overline{(Kf)(t)}{\rm d}t
\nonumber
\\
\qquad=\int_0^1\overline{f(s)}\left(\int_0^sf_1(t){\rm d}t\right){\rm d}s+\int_0^1f_1(t)\left(\int_0^t\overline{f(s)}
{\rm d}s\right){\rm d}t
\nonumber
\\
\qquad=\int_0^1f_1(t)\left(\int_t^1\overline{f(s)}{\rm d}s\right){\rm d}t+\int_0^1f_1(t)\left(\int_0^t\overline{f(s)}
{\rm d}s\right){\rm d}t
\nonumber
\\
\qquad=\int_0^1f_1(t)\left(\int_0^1\overline{f(s)}{\rm d}s\right){\rm d}t=\langle f_1,u\rangle\langle u,f\rangle=0.
\label{equation: 5.1}
\end{gather}

Suppose that $g_1,g_2\in{\mathscr{D}}_3$, and let $g_j=Kf_j+a_ju$, where $f_1,f_2\in{\mathcal H}_1$ and
$a_1,a_2\in \mathbb{C}$.
Since $\langle f_j,u\rangle=0$, from~\eqref{equation: 5.1} we have
\begin{gather*}
\langle D_3g_1,g_2\rangle+\langle g_1,D_3g_2\rangle
=\langle f_1,Kf_2+a_2u\rangle+\langle Kf_1+a_1u,f_2\rangle
=\langle f_1,Kf_2\rangle+\langle Kf_1,f_2\rangle=0.
\end{gather*}
Thus $g_2\in{\mathscr{D}}(D_3^*)$, and $D_3^*g_2=-D_3g_2$; so $-D_3\subseteq D_3^*$.

It remains to show that ${\mathscr{D}}(D_3^*)\subseteq {\mathscr{D}}_3$.
Suppose that $g\in{\mathscr{D}}(D_3^*)$, and $D_3^*g=h$.
For any $f_1\in{\mathcal H}_1$ and $a\in \mathbb{C}$, $Kf_1+au\in{\mathscr{D}}_3$, and $D_3(Kf_1+au)=f_1$.
Thus
\begin{gather*}
\langle f_1,g\rangle=\langle D_3(Kf_1+au),g\rangle=\langle Kf_1+au,h\rangle.
\end{gather*}
By varying $a$, it follows that $\langle h,u\rangle=0$; so $h\in{\mathcal H}_1$, and
$\langle f_1,g\rangle=\langle Kf_1,h\rangle$.
From~\eqref{equation: 5.1}, we now have $\langle f_1,g\rangle=-\langle f_1,Kh\rangle$
$(f_1\in{\mathcal H}_1)$.
Thus $g+Kh\in {\mathcal H}_1^\perp=[u]$, and $g=-Kh+au$ for some scalar $a$.
Thus $g\in{\mathscr{D}}_3$, and ${\mathscr{D}}(D_3^*)\subseteq {\mathscr{D}}_3$.
\end{proof}

\subsection{The classic representation}\label{s5.3}

Given the discussion and results to this point, what are we to understand by a~``representation of the
Heisenberg relation'', $QP-PQ=i\hbar I$?
Having proved that this representation cannot be achieved with
f\/inite matrices in place of~$Q$ and $P$ and $I$, nor even with bounded operators on a~Hilbert space, nor
elements~$Q$,~$P$,~$I$ in a~complex Banach algebra, we begin to examine the possibility that this representation
can be ef\/fected with unbounded operators for~$Q$ and~$P$.
It is ``rumored'', loosely, that~$Q$, which is associated with the physical observable ``position'' on~${\mathbb{R}}$, and~$P$, which is associated with the (conjugate) ``momentum'' observable, will provide
such a~representation.
The observable~$Q$ is modeled, nicely, by the self-adjoint operator, multiplication by $x$ on
$L_2({\mathbb{R}})$, with domain those~$f$ in~$L_2({\mathbb{R}})$ such that $xf$ is in $L_2({\mathbb{R}})$.
The observable~$P$ is modeled by $i\frac{{\rm d}}{{\rm d}t}$, dif\/ferentiation on some appropriate domain of
dif\/ferentiable functions with derivatives in $L_2({\mathbb{R}})$.
But $QP-PQ$ certainly can't equal $i\hbar I$, since its domain is contained in
${\mathscr{D}}(Q)\cap{\mathscr{D}}(P)$, which is not~${\mathcal H}$.
The domain of~$P$ must be chosen so that $P$ is self-adjoint and ${\mathscr{D}}(QP-PQ)$ is dense in~${\mathcal H}$ and $QP-PQ$ agrees with~$i\hbar I$ on this dense domain.
In particular, $QP-PQ\subseteq i\hbar I$.
Since $i\hbar I$ is bounded, it is closed, and $QP-PQ$ is closable with closure~$i\hbar I$.
We cannot insist that, with the chosen domains for~$Q$ and~$P$, $QP-PQ$ be skew-adjoint, for then it would
be closed, bounded, and densely def\/ined, hence, everywhere def\/ined.
In the end, we shall mean by ``a representation of the Heisenberg relation $QP-PQ=i\hbar I$ on the Hilbert
space ${\mathcal H}$'' a~choice of self-adjoint operators $Q$ and $P$ on ${\mathcal H}$ such that $QP-PQ$
has closure $i\hbar I$.

As mentioned above, the classic way~\cite{VNIV} to represent the Heisenberg relation $QP-PQ=i\hbar I$ with unbounded
self-adjoint operators $Q$ and $P$ on a~Hilbert space ${\mathcal H}$ is to realize ${\mathcal H}$ as
$L_2({\mathbb{R}})$, the space of square-integrable, complex-valued functions on~${\mathbb{R}}$ and~$Q$ and~$P$ as, respectively, the operator~$Q$ corresponding to multiplication by~$x$, the identity transform on
${\mathbb{R}}$, and the operator~$P$ corresponding to $i\frac{{\rm d}}{{\rm d}t}$, where~$\frac{{\rm d}}{{\rm d}t}$ denotes
dif\/ferentiation, each of~$Q$ and~$P$ with a~suitable domain in~$L_2({\mathbb{R}})$.
The domain of~$Q$ consists of those~$f$ in~$L_2({\mathbb{R}})$ such that~$xf$ is in~$L_2({\mathbb{R}})$.
The operator~$\frac{{\rm d}}{{\rm d}t}$ is intended to be dif\/ferentiation on~$L_2({\mathbb{R}})$, where that
dif\/ferentiation makes sense~-- certainly, on every dif\/ferentiable functions with derivative in~$L_2({\mathbb{R}})$.
However, specifying a~dense domain, precisely, including such functions, on which ``dif\/ferentiation'' is
a~self-adjoint operator is not so simple.
A step function, a~function on~${\mathbb{R}}$ that is constant on each connected component of an open dense
subset of~${\mathbb{R}}$ (those components being open intervals) has a~derivative almost everywhere (at all
but the set of endpoints of the intervals~-- a~countable set), and that derivative is~0.
The set of such step functions in~$L_2({\mathbb{R}})$ is dense in~$L_2({\mathbb{R}})$, as is their linear
span.
To include that linear span in a~proposed domain for our dif\/ferentiation operator condemns any closed
operator extending our dif\/ferentiation operator to be the everywhere-def\/ined operator~0.
Of course, that is not what we are aiming for.
Another problem that we face in this discussion is that of ``mixing'' measure theory with dif\/ferentiation.
We speak, loosely, of elements of our Hilbert space~$L_2({\mathbb{R}})$ as ``functions''.
We have learned to work quickly and accurately with the mathematical convenience that this looseness
provides us avoiding such pitfalls as taking the union of ``too many'' sets of measure~0 in the process.
The elements of~$L_2({\mathbb{R}})$ are, in fact, equivalence classes of functions dif\/fering from one
another on sets of measure~0.
On the other hand, dif\/ferentiation is a~process that focuses on points, each point being a~set of
Lebesgue measure zero.
When we speak of the~$L_2$-norm of a~\textit{function} in~$L_2({\mathbb{R}})$ it doesn't matter which
function in the \textit{class} in question we work with; they all have the same norm.
It is not the same with dif\/ferentiability.
Not each function in the class of an everywhere dif\/ferentiable function is everywhere dif\/ferentiable.
There are functions in such classes that are nowhere dif\/ferentiable, indeed, nowhere continuous (at each
point of dif\/ferentiability a~function is continuous).
The measure class of each function on~${\mathbb{R}}$ contains a~function that is nowhere continuous.
To see this, choose two disjoint, countable, everywhere-dense subsets, for example, the rationals~$\mathbb{Q}$ in ${\mathbb{R}}$ and $\mathbb{Q}+\sqrt{2}$.
With~$f$ a~given function on~${\mathbb{R}}$, the function $g$ that agrees with~$f$, except on~$\mathbb{Q}$
where it takes the value~0 and on~$\mathbb{Q}+\sqrt{2}$ where it takes the value~1 is in the measure class
of~$f$ and is continuous nowhere (since each non-null open set in~${\mathbb{R}}$ contains a~point at which~$g$ takes the value~0 and a~point at which it takes the value~1).
These are some of the problems that arise in dealing with an appropriate domain for~$\frac{{\rm d}}{{\rm d}t}$.

There is an elegant way to approach the problem of f\/inding precisely the self-adjoint operator and its
domain that we are seeking.
That approach is through the use of ``Stone's theorem''~\cite{S} (from the very beginning of the theory of unitary
representations of inf\/inite groups).
We start with a~clear statement of the theorem.
Particular attention should be paid to the description of the domain of the generator~$iH$ in this
statement.
\begin{Theorem}[Stone's theorem] If~$H$ is a~$($possibly unbounded$)$ self-adjoint operator on the Hilbert space
${\mathcal H}$, then $t\to\exp itH$ is a~one-parameter unitary group on~${\mathcal H}$.
Conversely, if~$t\to U_t$ is a~one-parameter unitary group on~${\mathcal H}$, there is a~$($possibly
unbounded$)$ self-adjoint operator~$H$ on ${\mathcal H}$ such that $U_t=\exp itH$ for each real~$t$.
The domain of~$H$ consists of precisely those vectors~$x$ in~${\mathcal H}$ for which $t^{-1}(U_tx-x)$
tends to a~limit as~$t$ tends to~$0$, in which case this limit is~$iHx$.
\end{Theorem}

The relevance of Stone's theorem emerges from the basic case of the one-parameter unitary group $t\to U_t$
on $L_2({\mathbb{R}})$, where $(U_tf)(s)=f(s+t)$.
That is, $U_t$ is ``translation by $t$''.
In this case, $U_t=\exp itH$, with $H$ a~self-adjoint operator acting on $L_2({\mathbb{R}})$.
The domain of~$H$ consists of those $f$ in $L_2({\mathbb{R}})$ such that $t^{-1}(U_tf-f)$ tends to a~limit
$g$ in $L_2({\mathbb{R}})$, as $t$ tends to 0, in which case, $iHf=g$.
We treat $\frac{{\rm d}}{{\rm d}t}$ as the inf\/initesimal generator of this one-parameter unitary group.
An easy measure-theoretic argument shows that this one-parameter unitary group is strong-operator
continuous on ${\mathcal H}$.
That is, $U_tf\to U_{t'}f$, in the norm topology of~${\mathcal H}$, as $t\to t'$, for each $f$ in
${\mathcal H}$, or what amounts to the same thing, since $t\to U_t$ is a~one-parameter group, if~$U_{t''}f=U_{t-t'}f\to f$, when $(t-t')=t''\to 0$ for each $f$ in $L_2({\mathbb{R}})$.
From Stone's theorem, there is a~skew-adjoint (unbounded) operator ($iH$) we denote by $\frac{{\rm d}}{{\rm d}t}$ on
${\mathcal H}$ such that $U_t=\exp t\frac{{\rm d}}{{\rm d}t}$ for each real~$t$.
The domain of~$\frac{{\rm d}}{{\rm d}t}$ consists of those~$f$ in $L_2({\mathbb{R}})$ such that $t^{-1}(U_tf-f)$ tends
to some~$g$ in~$L_2({\mathbb{R}})$ as $t$ tends to~0, in which case $g=\frac{{\rm d}}{{\rm d}t}f$.

Now, let us make some observations to see how Stone's theorem works in our situation.
Our aim, at this point, is to study just which functions are and are not in the domain of~$\frac{{\rm d}}{{\rm d}t}$.
(This study will make clear how apt the notation $\frac{{\rm d}}{{\rm d}t}$ is for the inf\/initesimal generator of the
group of real translations of~${\mathbb{R}}$.) To begin with, Stone's theorem requires us to study the
convergence behavior of~$t^{-1}(U_tf-f)$ as $t$ tends to 0.
This requirement is to study the convergence behavior in the Hilbert space metric (in the ``mean of order
2'', in the terminology of classical analysis), but there is no harm in examining how $t^{-1}(U_tf-f)$
varies pointwise with $t$ at points $s$ in ${\mathbb{R}}$.
For this, note that
\begin{gather*}
\big(t^{-1}(U_tf-f)\big)(s)=\frac{f(s+t)-f(s)}{t}\to f'(s)
\qquad
\text{as}
\qquad
t\to0,
\end{gather*}
which suggests $f'$ as the limit of~$t^{-1}(U_tf-f)$ when $f$ is dif\/ferentiable with $f'$ in
$L_2({\mathbb{R}})$ (and motivates the use of the notation ``$\frac{{\rm d}}{{\rm d}t}$'' for the inf\/initesimal
generator of~$t\to U_t$).
However, the ``instructions'' of Stone's theorem tell us to f\/ind $g$ in $L_2({\mathbb{R}})$ such that
\begin{gather*}
\int\left|\frac{f(s+t)-f(s)}{t}-g(s)\right|^2{\rm d}\mu(s)\to0
\end{gather*}
as $t\to 0$, where $\mu$ is Lebesgue measure on ${\mathbb{R}}$.
Our f\/irst observation is that if~$f$ fails to have a~derivative at some point $s_0$ in ${\mathbb{R}}$
\textit{in an essential way}, then $f$ is not in the domain of~$\frac{{\rm d}}{{\rm d}t}$.
This may be surprising, at f\/irst, for the behavior of a~function at a~point rarely has (Lebesgue)
measure-theoretic consequences.
In the present circumstances, we shall see that the ``local'' nature of dif\/ferentiation can result in
exclusion from the domain of an unbounded dif\/ferentiation operator because of non-dif\/ferentiability at
a~single point.

We begin with a~def\/inition of ``jump in a~function'' that is suitable for our measure-theoretic situation.
\begin{Definition}
We say that $f$ has \textit{jump} $a$ $({\ge}0)$ for width $\delta$ $({>}0)$ at $s_0$ in ${\mathbb{R}}$ when
$\inf\{f(s)\}$ with $s$ in one of the intervals $[s_0-\delta,s_0)$ or $(s_0,s_0+\delta]$ is
$a+\sup\{f(s)\}$ with $s$ in the other of those intervals.
\end{Definition}

Typically, one speaks of a~``jump discontinuity'' when $\lim\limits_{s\to s^-_0}f(s)$ and $\lim\limits_{s\to s^+_0}f(s)$
exist and are distinct.
In the strictly measure-theoretic situation with which we are concerned, the concept of ``jump'', as just
def\/ined, seems more appropriate.

\begin{Remark}\label{remark: jump}
If~$f$ has a~jump $a$ for width $\delta$ at some point $s_0$ in ${\mathbb{R}}$, then $U_{s_0}f$ has
a~jump $a$ for width $\delta$ at 0, and $bU_{s_0}f$ has jump $ba$ for width $\delta$ at 0 when $0<b$.
Letting $f_r$ be the function whose value at $s$ is $f(rs)$, one has that~$f_r$ has a~jump $a$ at
$r^{-1}s_0$ for width $r^{-1}\delta$.
Thus $a^{-1}(U_{s_0}f)_\delta$ has jump~1 at~0 for width~1.
\end{Remark}

\begin{Theorem}[\protect{cf.~\cite[Theorem 4.6]{Z}}] If~$f$ has a~positive jump, then $f\notin{\mathscr{D}}\big(\frac{{\rm d}}{{\rm d}t}\big)$.
\end{Theorem}
\begin{proof}
We shall show that $\|t^{-1}(U_tf-f)\|$ is unbounded for $t$ in each open interval in ${\mathbb{R}}$
containing 0.
Of course, this is so if and only if~$\|t^{-1}bU_s(U_tf-f)\|$ is unbounded for each given positive $b$ and
$U_s$.
Thus, from Remark~\ref{remark: jump}, it will suf\/f\/ice to show that $\|t^{-1}(U_tf-f)\|$ is unbounded
when $f$ has jump 1 at 0.
Noting that $\|g_r\|=r^{-1}\|g\|$ for $g$ in $L_2({\mathbb{R}})$, that $(g+h)_r=g_r+h_r$, and that
($U_tf)_r=U_{r^{-1}t}f_r=U_{t'}f_r$, where $t'=r^{-1}t\to0$ as $t\to 0$, we have that
\begin{gather*}
r^{-1}t^{-1}\|U_tf-f\|=t^{-1}\|(U_tf-f)_r\|=t^{-1}\|(U_tf)_r-f_r\|
\\
\hphantom{r^{-1}t^{-1}\|U_tf-f\|}{}  =t^{-1}\|U_{r^{-1}t}f_r-f_r\|=r^{-1}t'^{-1}\|U_{t'}f_r-f_r\|.
\end{gather*}
Thus $\|t^{-1}(U_tf-f)\|=\|t'^{-1}(U_{t'}f_r-f_r)\|$.
It follows that $\|t^{-1}(U_tf-f)\|$ is bounded for~$t$ near~0 if and only if~$\|t'^{-1}(U_{t'}f_r-f_r)\|$
is.
This holds for each positive $r$, in particular, when $r$ is $\delta$, where $f$ has jump~1 at~0 for width
$\delta$.
Since $f_\delta$ has jump~1 at~0 for width~1 $({=}\delta^{-1}\delta$)), from Remark~\ref{remark: jump}, it
will suf\/f\/ice to show that $\|t^{-1}(U_tf-f)\|$ is unbounded for~$t$ near~0, when~$f$ has jump~1 at~0
for width~1.
We shall do this by f\/inding a~sequence $t_2,t_3,\ldots$ of positive numbers~$t_j$ tending to~0 such that
$\|t_j^{-1}(U_{t_j}f-f)\|\to\infty$ as $j\to\infty$.
We assume that $f$ has jump~1 at~0 for width~1.
In this case, $|f(s')-f(s'')|\ge1$ when $s'\in[-1,0)$ and $s''\in(0,1]$.
Thus, when $t_n=\frac1{n-1}$,
\begin{gather*}
\|t_n^{-1}(U_{t_n}f-f)\|^2=\int_{\mathbb{R}}\left|t_n^{-1}(U_{t_n}f-f)\right|^2(s){\rm d}\mu(s)
\\
\phantom{\|t_n^{-1}(U_{t_n}f-f)\|^2}
\ge\int_{[-\frac{1}{n},0)}\left|(n-1)(f(s+t_n)-f(s))\right|^2{\rm d}\mu(s)
\ge\tfrac1n(n-1)^2=n-2+\tfrac{1}{n}.
\end{gather*}
It follows that $\|(n-1)(U_{(n-1)^{-1}}f-f)\|\to\infty$ as $n\to\infty$.
Hence $t^{-1}(U_tf-f)$ has no limit in~$L_2({\mathbb{R}})$ as $t\to0$ and $f\notin{\mathscr{D}}(\frac{{\rm d}}{{\rm d}t})$.
\end{proof}

\begin{Theorem}[\protect{cf.~\cite[Theorem~4.7]{Z}}]\label{thm: roundoff} If~$f_1$ is a~continuously differentiable function on
${\mathbb{R}}$ such that $f_1$ and $f_1'$ are in $L_2({\mathbb{R}})$,
then $f_1\in{\mathscr{D}}(\frac{{\rm d}}{{\rm d}t})$; and $\frac{{\rm d}}{{\rm d}t}(f_1)=f_1'$.
\end{Theorem}

\begin{proof}
We prove, f\/irst, that if~$f$, in $L_2({\mathbb{R}})$, vanishes outside some interval $[-n,n]$, with $n$
a~positive integer, and $f$ is continuously dif\/ferentiable on ${\mathbb{R}}$ with derivative $f'$ in
$L_2({\mathbb{R}})$, then $f\in{\mathscr{D}}(\frac{{\rm d}}{{\rm d}t})$ and $\frac{{\rm d}}{{\rm d}t}(f)=f'$.

From Stone's theorem, we must show that $\|t^{-1}(U_tf-f)-f'\|_2\to0$ as $t\to0$.
Now,
\begin{gather*}
\big\|t^{-1}(U_tf-f)-f'\big\|_2^2=\int_{[-n,n]}\left|\big[t^{-1}(U_tf-f)-f'\big](s)\right|^2{\rm d}\mu(s)
\\
\phantom{\big\|t^{-1}(U_tf-f)-f'\big\|_2^2}
=\int_{[-n,n]}\left|\frac{f(s+t)-f(s)}t-f'(s)\right|^2{\rm d}\mu(s).
\end{gather*}
Note that $t^{-1}(U_tf-f)-f'$ tends to 0 (pointwise) everywhere on ${\mathbb{R}}$ as $t$ tends to 0.
Of course, $t^{-1}(U_tf-f)$ and $f'$ vanish outside of~$[-(n+1),n+1]$ when $|t|<1$.
Since $f$ is dif\/ferentiable, it is continuous and bounded on $[-(n+1),n+1]$.
By assumption, $f'$ is continuous, hence bounded on $[-(n+1),n+1]$ (on ${\mathbb{R}}$).
Say, $|f'(s)|\le M$, for each $s$.
From the mean value theorem, for~$s$ in $[-n,n]$,
\begin{gather*}
\left|t^{-1}(U_tf-f)(s)\right|=\left|\frac{f(s+t)-f(s)}t\right|=|f'(s')|\le M,
\end{gather*}
for some $s'$ in the interval with endpoints $s$ and $s+t$.
Thus $|t^{-1}(U_tf-f)|$ is bounded by $M$, on $[-n,n]$ for all $t$ in $(-1,1)$.
At the same time, $t^{-1}(U_tf-f)$ tends to $f'$ everywhere (that is, pointwise) on $[-n,n]$.
From Egorof\/f's theorem, $t^{-1}(U_tf-f)$ tends \textit{almost uniformly} to $f'$ on $[-n,n]$ as $t$ tends
to 0.
Hence, given a~positive $\varepsilon$, there is a~subset $S$ of~$[-n,n]$ of measure less than
$\varepsilon/8M^2$ such that $t^{-1}(U_tf-f)$ converges uniformly to $f'$ on $[-n,n]\setminus S$.

We show, now, that $t^{-1}(U_tf-f)$ converges to $f'$ in $L_2({\mathbb{R}})$.
With $\varepsilon$ and $S$ chosen as in the preceding paragraph, by uniform convergence on $[n,-n]\setminus
S$, we f\/ind a~positive $\delta$ such that for $0<|t|<\delta$, and $s$ in $[-n,n]\setminus S$,
$\big|t^{-1}(f(s+t)-f(s))-f'(s)\big|^2<\varepsilon/4n$.
Hence, when $0<|t|<\delta$,
\begin{gather*}
\big\|t^{-1}(U_tf-f)-f'\big\|_2^2
\\
\qquad=\int_{[-n,n]\setminus S}\left|\frac{f(s+t)-f(s)}t-f'(s)\right|^2{\rm d}\mu(s)
+\int_S\left|\frac{f(s+t)-f(s)}t-f'(s)\right|^2{\rm d}\mu(s)
\\
\qquad
\le2n\frac\varepsilon{4n}+4M^2\frac\varepsilon{8M^2}=\varepsilon.
\end{gather*}
The desired convergence of~$t^{-1}(U_tf-f)$ to $f'$ in $L_2({\mathbb{R}})$ follows from this.

With $f_1$ as in the statement of this theorem, suppose that we can f\/ind $f$ as in the preceding
discussion (that is, vanishing outside a~f\/inite interval) such that $\|f_1-f\|_2$ and $\|f_1'-f'\|_2$ are
less than a~preassigned positive $\varepsilon$.
Then $(f_1,f_1')$ is in the closure of the graph of~$\frac{{\rm d}}{{\rm d}t}$, since each $(f,f')$ is in that closure
from what we have proved.
But $\frac{{\rm d}}{{\rm d}t}$ is skew-adjoint (from Stone's theorem); hence, $\frac{{\rm d}}{{\rm d}t}$ is closed.
Thus, if we can ef\/fect the described approximation of~$f_1$ and~$f_1'$ by~$f$ and~$f'$, it will follow
that $f_1\in{\mathscr{D}}(\frac{{\rm d}}{{\rm d}t})$ and~$\frac{{\rm d}}{{\rm d}t}(f_1)=f_1'$.

Since $f_1$ and $f_1'$ are continuous and in $L_2({\mathbb{R}})$, the same is true for
$|f_1|+|f_1^-|+|f_1'|+|f_1'^-|$, where $g^-(s)=g(-s)$ for each $s$ in ${\mathbb{R}}$ and each
complex-valued function $g$ on ${\mathbb{R}}$.
(Note, for this, that $s\to-s$ is a~Lebesgue-measure-preserving homeomorphism of~${\mathbb{R}}$ onto
${\mathbb{R}}$.) It follows that, for each positive integer $n$, there is a~real $s_n$ such that $n<s_n$ and
\begin{gather*}
|f_1(s_n)|+|f_1(-s_n)|+|f_1'(s_n)|+|f_1'(-s_n)|<\frac1n.
\end{gather*}
(Otherwise, $|f_1(s)|+|f_1(-s)|+|f_1'(s)|+|f_1'(-s)|\ge\frac1n$, for each $s$ outside of~$[-n,n]$,
contradicting the fact that $|f_1|+|f_1^-|+|f_1'|+|f_1'^-|\in L_2({\mathbb{R}})$.) We can choose $s_n$ such
that $s_{n-1}<s_n$.
Since $n<s_n$, we have that $s_n\to\infty$ as $n\to\infty$, and
\begin{gather*}
\int_{[-s_n,s_n]}|h(s)|^2{\rm d}\mu(s)\to\|h\|_2^2,
\qquad
n\to\infty,
\end{gather*}
for each $h$ in $L_2({\mathbb{R}})$.
Thus $\|h-h^{(n)}\|_2\to0$ as $n\to\infty$, where $h^{(n)}$ is the function that agrees with $h$ on
$[-s_n,s_n]$ and is 0 outside this interval.
With $\varepsilon$ (${<}1$) positive, there is an $n_0$ such that, if~$n>n_0$, then each of~$\|f_1-f_1^{(n)}\|_2$, $\|f_1^--{f_1^-}^{(n)}\|_2$, $\|f'_1-{f'_1}^{(n)}\|_2$, and
$\|{f'_1}^--{f'_1}^{-(n)}\|_2$ is less than $\frac\varepsilon 2$.
At the same time, we may choose $n_0$ large enough so that $\frac1n<\frac\varepsilon4$ when $n>n_0$.
For such an $n$, a~``suitably modif\/ied'' $f_1^{(n)}$ will serve as the desired $f$ for our approximation.
In the paragraphs that follow, we describe that modif\/ication.

Our aim is to extend $f_1^{(n)}$ to ${\mathbb{R}}$ from $[-s_n,s_n]$ so that the extension $f$ remains
continuously dif\/ferentiable with $f$ and $f'$ vanishing outside some f\/inite interval and so that the
projected approximations $\|f_1-f\|_2<\varepsilon$ and $\|f'_1-f'\|_2<\varepsilon$ are realized.
In ef\/fect, we want $\|f_1^{(n)}-f\|_2$ and $\|{f'_1}^{(n)}-f'\|_2$ to be less than $\frac\varepsilon 2$.
Combined, then, with our earlier choice of~$n_0$ such that, for $n>n_0$,
$\|f_1-f_1^{(n)}\|_2<\frac\varepsilon 2$ and $\|f'_1-{f'_1}^{(n)}\|_2<\frac\varepsilon 2$, we have the
desired approximation.

To construct $f$, we add to $f_1^{(n)}$ a~function $g$ continuous and continuously dif\/ferentiable on
$(-\infty,-s_n]\cup[s_n,\infty)$ such that $g(s_n)=f_1(s_n)$, $g'(s_n)=f_1'(s_n)$, $g(-s_n)=f_1(-s_n)$,
$g'(-s_n)=f_1'(-s_n)$, $g$ vanishes on $(-\infty,-s_n-1]\cup[s_n+1,\infty)$, and $\|g\|_2<\frac\varepsilon
2$, $\|g'\|_2<\frac\varepsilon 2$.
With $f$ so def\/ined, $\|f_1^{(n)}-f\|_2=\|g\|_2<\frac\varepsilon 2$ and
$\|{f'_1}^{(n)}-f'\|_2=\|g'\|_2<\frac\varepsilon 2$, as desired.
We describe the construction of~$g$ on $[s_n,\infty)$.
The construction of~$g$ on $(-\infty,-s_n]$ follows the same pattern.
We present the construction of~$g$ geometrically~-- with reference to the graphs of the functions involved.
The graphs are described in an $XY$ plane, where ${\mathbb{R}}$ is identif\/ied with the $X$-axis.
By choice of~$s_n$ and $n$ ($>n_0$), $|f_1(s_n))|<\frac\varepsilon4$, and $|f_1'(s_n)|<\frac\varepsilon4$.

Translating $s_n$ to the origin, we see that our task is to construct a~function $h$ on $[0,1]$
continuously dif\/ferentiable, 0 on $[\frac12,1]$, with given initial data $h(0)$, $h'(0)$ satisfying
$|h(0)|<\frac\varepsilon 4$, $|h'(0)|<\frac\varepsilon 4$ such that $\|h\|_2<\frac\varepsilon2$ and
$\|h'\|_2<\frac\varepsilon2$.
If~$h(0)=h'(0)=0$, then~$h$, with $h(x)=0$, for each~$x$ in $[0,1]$, will serve as our $h$.
If~$h'(0)\ne0$, we def\/ine $h$, f\/irst, on $[0,x_0]$, where $x_0=\frac12h(0)h'(0)$ and
$(y_0=)h(x_0)=\tfrac12h(0)\big[1+(1+h'(0)^2)^\frac12\big]$.
The restriction of~$h$ to $[0,x_0]$ has as its graph the (``upper, smaller'') arch of the circle with center
$(x_0,\frac12h(0))$ and radius $\frac12h(0)(1+h'(0)^2)^\frac12$ (tangent to the line with slope~$h'(0)$ at
$(0,h(0))$).
Note that $h(0)<y_0<2h(0)<\frac\varepsilon 2$ and that the circle described has a~horizontal tangent at
$(x_0,y_0)$; that is, $h'(x_0)=0$, as~$h$ has been def\/ined.

We complete the def\/inition of~$h$ by adjoining to the graph of~$h$ over $[0,x_0]$ the graph of~$\frac12y_0[\cos((\frac12-x_0)^{-1}\pi(x-x_0))+1]$ over $[x_0,\frac12]$.
Note that this graph passes through $(x_0,y_0)$ and $(\frac{1}{2},0)$.
Finally, we def\/ine $h(x)$ to be 0 when $x\in[\frac12,1]$.
As constructed, $h$ is continuously dif\/ferentiable on $[0,1]$.
Since $|h(x)|\le2|h((0)|<\frac\varepsilon 2$ for $x$ in $[0,\frac12]$ and $h$ vanishes on $[\frac12,1]$,
$\|h\|_2<\frac\varepsilon 2$.
\end{proof}

We may ask whether the converse statement to the preceding theorem holds as well.
Does a~function class in ${\mathscr{D}}(\frac{{\rm d}}{{\rm d}t})$ necessarily contain a~continuously dif\/ferentiable
function with derivative in $L_2({\mathbb{R}})$? As it turns out, there are more functions, not as well
behaved as continuously dif\/ferentiable functions, in the domain of~$\frac{{\rm d}}{{\rm d}t}$.
We shall give a~complete description of that domain in Theorem~\ref{thm: domain}.

Our notation and terminology has a~somewhat ``schizophrenic'' character to it~-- much in the style of the
way mathematics treats certain topics.
In the present instance, we use the notation `$L_2({\mathbb{R}})$' to denote, both, the collection (linear
space) of measurable functions $f$ such that $|f|^2$ is Lebesgue integrable on ${\mathbb{R}}$ and the
Hilbert space of (measure-theoretic) equivalence classes of such functions equipped with the usual Hilbert
space structure associated with $L_2$ spaces.
In most circumstances, there is no danger of serious confusion or misinterpretation.
In our present discussion of the domain of~$\frac{{\rm d}}{{\rm d}t}$, these dangers loom large.
We note, earlier in this section, that each measure-theoretic equivalence class of functions contains
a~function that is continuous at no point of~${\mathbb{R}}$.
It can make no sense to attempt to characterize special elements~$x$ of~$L_2({\mathbb{R}})$ by the
``smoothness'' properties of \textit{all} the functions in the equivalence class denoted by~`$x$' (their
continuity, dif\/ferentiability, and so forth).
Despite this, our next theorem describes the domain given to us by the generator, which we are denoting by
`$\frac{{\rm d}}{{\rm d}t}$', of the one-parameter unitary group $t\rightarrow U_t$ of translations of the equivalence
classes of functions in $L_2({\mathbb{R}})$ (to other such classes) in terms of smoothness properties.
However, these smoothness properties will be those of a~\textit{single} element in the class as we shall
see.
We note, f\/irst, that if an equivalence class contains a~continuous function~$f$ on~${\mathbb{R}}$, then~$f$ is the unique such function in the class.
This is immediate from the fact that $f-g$ vanishes nowhere on some non-null, open interval when~$f$ and~$g$ are distinct continuous functions, whence~$f$ and~$g$ dif\/fer on a~set of positive Lebesgue measure
and lie in dif\/ferent measure classes.

The unique continuous function in each measure class of some family of measure classes allows us to
distinguish subsets of this family by smoothness properties of that continuous function in the class.
In the case of the one-parameter unitary group induced by translations on~${\mathbb{R}}$, corresponding to
an element $x$ in the domain of the Stone generator $\frac{{\rm d}}{{\rm d}t}$, the measure class~$x$ contains
a~continuous function (hence, as noted, a~unique such function), and this function must be absolutely
continuous, in $L_2({\mathbb{R}})$, of course, with derivative almost everywhere on ${\mathbb{R}}$ in
$L_2({\mathbb{R}})$.
Moreover, an absolutely continuous function in $L_2({\mathbb{R}})$ with derivative almost everywhere in
$L_2({\mathbb{R}})$ has measure class an element of the Hilbert space on which the unitary group
(corresponding to the translations of~${\mathbb{R}}$) acts that lies in the domain of~$\frac{{\rm d}}{{\rm d}t}$.
So, this absolute-continuity smoothness, together with the noted $L_2$ restrictions, characterizes the
domain of~$\frac{{\rm d}}{{\rm d}t}$.
It is dangerously misleading to speak of the domain of~$\frac{{\rm d}}{{\rm d}t}$ as ``consisting of absolutely
continuous functions in $L_2$ with almost everywhere derivatives in $L_2$''; it consists of the measure
classes of such functions and each such class contains, as noted, a~function which is nowhere continuous.

We undertake, now, the proof of the theorem that describes the domain of~$\frac{{\rm d}}{{\rm d}t}$, the generator of~$t\rightarrow U_t$, the one-parameter unitary group corresponding to translations of~$L_2({\mathbb{R}})$ $({=}{\mathcal H})$.
(Compare \cite[Theorem~4.8]{Z}, where a~sketch of the proof is given. See, also,~\cite{HP}.) The following results in real
analysis will be useful to us~\cite{G,SS}.

\begin{Lemma}\label{lemma: 1}
Suppose that $f \in L_1 (\mathbb{R})$.
Let $F(x)=\int_0^xf(s){\rm d}s$.
Then $F$ is differentiable almost everywhere, and the derivative is equal to $f$ almost everywhere.
\end{Lemma}

\begin{Lemma}\label{lemma: 2}
If~$f \in L_2 (\mathbb{R})$, then $\frac{1}{t}\int_x^{x+t}f(s){\rm d}s \rightarrow f(x)$ in $L_2$ norm, as $t
\rightarrow0$.
\end{Lemma}

\begin{Theorem}\label{thm: domain}
The domain of~$\frac{{\rm d}}{{\rm d}t}$ is the linear subspace of measure classes in ${\mathcal
H}$ $({=}L_2({\mathbb{R}}))$ corresponding to absolutely continuous functions on ${\mathbb{R}}$ whose
almost-everywhere derivatives lie in $L_2({\mathbb{R}})$.
\end{Theorem}

\begin{proof}
Suppose $x\in {\mathscr{D}}(\frac{{\rm d}}{{\rm d}t})$.
Then, from Stone's theorem, there is a~vector $y$ in ${\mathcal H}$ such that, with $f$ in the measure
class~$x$ and~$g$ in the class~$y$,
\begin{gather*}
\left\|\frac{1}{t}(U_tx-x)-y\right\|_2^2
=\int_{\mathbb{R}}\left|\frac{1}{t}\big[f(s+t)-f(s)\big]-g(s)\right|^2{\rm d}s\rightarrow0,
\end{gather*}
as $|t|\to 0^+$.
With $a$ and $b$ in ${\mathbb{R}}$,
$\int_a^b\big|\tfrac{1}{t}\big[f(s+t)-f(s)\big]-g(s)\big|^2{\rm d}s\rightarrow 0$, and
\begin{gather*}
0\leq\left|\int_a^b\left(\frac{1}{t}\big[f(s+t)-f(s)\big]-g(s)\right){\rm d}s\right|\leq\int_a^b\left|\frac{1}
{t}\big[f(s+t)-f(s)\big]-g(s)\right|\cdot1{\rm d}s
\\
\hphantom{0}{}
\leq\left(\int_a^b\left|\frac{1}{t}\big[f(s+t)-f(s)\big]-g(s)\right|^2{\rm d}s\right)^{1/2}
\left(\int_a^b1{\rm d}s\right)^{1/2}\rightarrow0,
\qquad
|t|\to0^+.
\end{gather*}
Thus $\int_a^b\tfrac{1}{t}\big[f(s+t)-f(s)\big]{\rm d}s\rightarrow \int_a^bg(s){\rm d}s$ as $|t|\to 0^+$.
However,
\begin{gather*}
\int_a^b\frac{1}{t}\big[f(s+t)-f(s)\big]{\rm d}s=\frac{1}{t}\int_{a+t}^{b+t}f(s){\rm d}s-\frac{1}{t}
\int_a^bf(s){\rm d}s
\\
\phantom{\int_a^b\frac{1}{t}\big[f(s+t)-f(s)\big]{\rm d}s}
=\frac{1}{t}\int_b^{b+t}f(s){\rm d}s+\frac{1}{t}\int_{a+t}^bf(s){\rm d}s-\frac{1}{t}\int_a^bf(s){\rm d}s
\\
\phantom{\int_a^b\frac{1}{t}\big[f(s+t)-f(s)\big]{\rm d}s}
=\frac{1}{t}\int_b^{b+t}f(s){\rm d}s-\frac{1}{t}\int_a^{a+t}f(s){\rm d}s.
\end{gather*}
Now, from Lemma~\ref{lemma: 1}, $\tfrac{1}{t}\int_a^{a+t}f(s){\rm d}s$ and $\tfrac{1}{t}\int_b^{b+t}f(s){\rm d}s$ tend
to $f(a)$ and $f(b)$, respectively, as $|t|\to 0^+$, for almost every $a$ and $b$.
Choose $a$ for which this limit is valid.
Then, with this choice of~$a$, for almost all $b$, as noted, $\int_a^b\tfrac{1}{t}\big[f(s+t)-f(s)\big]{\rm d}s$
tends, as $|t|\to 0^+$, to $f(b)-f(a)$ and to $\int_a^bg(s){\rm d}s$.
Hence, for almost all $b$,
\begin{gather*}
f(b)=f(a)+\int_a^bg(s){\rm d}s.
\end{gather*}
Since $g\in L_2({\mathbb{R}})$, $g\in L_1([c,d])$, for each f\/inite interval $[c,d]$, and
$f(t)=f(a)+\int_a^tg(s){\rm d}s$ for almost all $t$, $f$ is in the measure class of~$h$ where
$h(t)=f(a)+\int_a^tg(s){\rm d}s$ for all real $t$.
Moreover, $h$ is absolutely continuous with almost-everywhere derivative $g$ in $L_2({\mathbb{R}})$.

Suppose, now, that $x$ in ${\mathcal H}$ $({=}L_2({\mathbb{R}}))$ contains an absolutely continuous function~$f$
with almost-everywhere derivative $g$ in $L_2({\mathbb{R}})$.
Let $y$ be the measure class of~$g$.
With this notation, $\frac{1}{t}[f(s+t)-f(s)]$ tends to $g(s)$ for almost every $s$ as $|t|\to 0^+$.
Now $f(s)=\int_0^sg(r){\rm d}r+f(0)$, so that
\begin{gather*}
\frac{1}{t}\big[f(s+t)-f(s)\big]
=\frac{1}{t}\left[\int_0^{s+t}g(r){\rm d}r-\int_0^sg(r){\rm d}r\right]=\frac{1}{t}\int_s^{s+t}g(r){\rm d}r,
\end{gather*}
and $\frac{1}{t}\int_s^{s+t}g(r){\rm d}r$ $({=}g_t(s))$ tends to $g$ in $L_2$ norm as $|t|\to 0^+$ (Lemma~\ref{lemma: 2}).
\end{proof}

We now describe a~core, for $\frac{{\rm d}}{{\rm d}t}$, that is particularly useful for computations.

\begin{Theorem}[\protect{cf.~\cite[Theorem 4.9]{Z}}] The family ${\mathscr{D}}_0$ of functions in $L_2 ({\mathbb{R}})$
that vanish outside a~finite interval and are continuously differentiable with derivatives in $L_2
({\mathbb{R}})$ determines a~core for the generator $\frac{{\rm d}}{{\rm d}t}$ of the one-parameter, translation,
unitary group on $L_2 ({\mathbb{R}})$.
\end{Theorem}
\begin{proof}
Suppose $f$ is the (unique) continuous function in a~measure class $\{f\}$ in ${\mathscr{D}}(\frac{{\rm d}}{{\rm d}t})$.
Suppose, moreover, $f$ is continuously dif\/ferentiable with derivative in $L_2 ({\mathbb{R}})$.
For any $\varepsilon>0$, there is a~positive integer $N$ ($N\geq 1$) such that
\begin{gather*}
\big\|f-f_{[-N,N]}\big\|_2<\frac{\varepsilon}{2}
\qquad
\text{and}
\qquad
\big\|f'-f_{[-N,N]}'\big\|_2<\frac{\varepsilon}{2},
\end{gather*}
where $f_{[-N,N]}$ denotes the function on ${\mathbb{R}}$ that agrees with $f$ on $[-N,N]$ and is $0$
outside $[-N,N]$.

From Theorem~\ref{thm: domain}, $f$ is absolutely continuous on ${\mathbb{R}}$; hence $f_{[-N,N]}$ is
absolutely continuous on $[-N,N]$.
Thus, $f_{[-N,N]}$ is dif\/ferentiable almost everywhere on $[-N,N]$ with derivative $f'_{[-N,N]}$ (in
$L_2([-N,N])$) and
\begin{gather*}
f_{[-N,N]}(x)=\int_{-N}^xf'_{[-N,N]}(s){\rm d}s+f_{[-N,N]}(-N),
\qquad
x\in[-N,N],
\end{gather*}
from the absolute continuity of~$f_{[-N,N]}$ on $[-N,N]$.

We approximate $f'_{[-N,N]}$ by a~continuous function $g'_N$ on $[-N,N]$ so that $\|f'_{[-N,N]}-g_N'\|_2<
\varepsilon/8N$.
Now, comparing the indef\/inite integrals,
\begin{gather*}
f_{[-N,N]}(x)=\int_{-N}^xf'_{[-N,N]}(s){\rm d}s+f_{[-N,N]}(-N)
\end{gather*}
and
\begin{gather*}
g_N(x)=\int_{-N}^xg'_N(s){\rm d}s+f_{[-N,N]}(-N),
\end{gather*}
we have
\begin{gather*}
|f_{[-N,N]}(x)-g_N(x)|=\left|\int_{-N}^x\big[f'_{[-N,N]}(s)-g_N'(s)\big]{\rm d}s\right|
\\
\hphantom{|f_{[-N,N]}(x)-g_N(x)|}{}
 \leq\left(\int_{-N}^N\left|f'_{[-N,N]}(s)-g_N'(s)\right|^2{\rm d}s\right)^{\frac{1}{2}}\left(\int_{-N}
^N|1|^2{\rm d}s\right)^{\frac{1}{2}}<\frac{\varepsilon}{8N}\sqrt{2N}.
\end{gather*}
Hence
\begin{gather*}
\|f_{[-N,N]}-g_N\|_2=\left(\int_{-N}^N\left|f_{[-N,N]}(x)-g_N(x)\right|^2{\rm d}x\right)^{\frac{1}{2}}
<\frac{\varepsilon}{8N}2N=\frac{\varepsilon}{4}.
\end{gather*}

Using the technique in the proof of Theorem~\ref{thm: roundoff}, we extend $g_N$ to ${\mathbb{R}}$ from
$[-N,N]$ so that the extension $g$ remains continuously dif\/ferentiable with $g$ and $g'$ vanishing
outside some f\/inite interval and
\begin{gather*}
\|g_N-g\|_2<\frac{\varepsilon}{4}
\qquad
\text{and}
\qquad
\|g_N'-g'\|_2<\frac{\varepsilon}{4}.
\end{gather*}
Then
\begin{gather*}
\|f_{[-N,N]}-g\|_2\leq\|f_{[-N,N]}-g_N\|_2+\|g_N-g\|_2<\frac{\varepsilon}{4}+\frac{\varepsilon}{4}
=\frac{\varepsilon}{2},
\\[1.11mm]
\|f_{[-N,N]}'-g'\|_2\leq\|f'_{[-N,N]}-g_N'\|_2+\|g_N'-g'\|_2<\frac{\varepsilon}{8N}+\frac{\varepsilon}{4}
<\frac{\varepsilon}{2}.
\end{gather*}
Finally,
\begin{gather*}
\|f-g\|_2\leq\big\|f-f_{[-N,N]}\big\|_2+\big\|f_{[-N,N]}-g\big\|_2<\frac{\varepsilon}{2}+\frac{\varepsilon}{2}
=\varepsilon,
\\[1.11mm]
\|f'-g'\|_2\leq\big\|f'-f'_{[-N,N]}\big\|_2+\big\|f'_{[-N,N]}-g'\big\|_2<\frac{\varepsilon}{2}+\frac{\varepsilon}{2}
=\varepsilon.
\end{gather*}
Thus, if~$(\{f\},\{f'\})\in \mathscr{G}(\frac{{\rm d}}{{\rm d}t})$, it can be approximated as closely as we wish by
$(\{g\},\{g'\})$ with $g\in {\mathscr{D}}_0$.
It follows that ${\mathscr{D}}_0$ is a~core for $\frac{{\rm d}}{{\rm d}t}$.
\end{proof}

In the classic representation of the Heisenberg relation, $QP-PQ=i\hbar I$, the operator $Q$ corresponds to
multiplication by $x$, the identity transform on ${\mathbb{R}}$.
The domain of~$Q$ consists of measure classes of functions $f$ in $L_2 ({\mathbb{R}})$ such that $xf$ is in
$L_2 ({\mathbb{R}})$.
Elementary measure-theoretic considerations establish that ${\mathscr{D}}_0$ is also a~core for $Q$.
Moreover, ${\mathscr{D}}_0\subseteq {\mathscr{D}}(QP)\cap{\mathscr{D}}(PQ)$, that is, ${\mathscr{D}}_0$ is
contained in the domain of~$QP-PQ$.
A calculation, similar to the one at the end of Example~\ref{example: L2}, shows that
\begin{gather*}
\big[QP-PQ\big]\big|{\mathscr{D}}_0=-iI\big|{\mathscr{D}}_0.
\end{gather*}
Moreover, for any $\{f\}\in{\mathscr{D}}$ ($={\mathscr{D}}(QP-PQ)$, the domain of~$QP-PQ$), with $f$ the
unique continuous function in the measure class $\{f\}$, for all $t$ at which $f$ is dif\/ferentiable,
\begin{gather*}
t\left(i\frac{{\rm d}}{{\rm d}t}f\right)(t)-\left(i\frac{{\rm d}}{{\rm d}t}\right)(tf(t))=itf'(t)-\big(if(t)+itf'(t)\big)=-if(t).
\end{gather*}
Thus
\begin{gather*}
\big[QP-PQ\big]\big|{\mathscr{D}}=-iI\big|{\mathscr{D}}.
\end{gather*}

As noted, the family of continuously dif\/ferentiable functions on ${\mathbb{R}}$ vanishing outside
f\/inite intervals constitutes a~very useful core for $\frac{{\rm d}}{{\rm d}t}$ for computing purposes.
It may be made even more useful, for these purposes, by introducing a~class of polynomials associated with
an $f$ in this core, the \textit{Bernstein polynomials}~\cite{B}, $B_n(f)$ $(n=1,2,\dots)$, which have remarkable
approximation properties.
We shall show that $\{B_n(f)\}$ tends uniformly to $f$ and $\{B'_n(f)\}$, the derivatives of~$\{B_n(f)\}$
(not $\{B_n(f')\}$, in general!), tends uniformly to $f'$.
Thus the set of Bernstein polynomials $B_n(f)$ with $f$ in the core we are studying, while not a~linear
space, hence not a~core for~$\frac{{\rm d}}{{\rm d}t}$, generates a~subset $\{(B_n(f),B'_n(f))\}$ of~${\mathscr{G}}(\frac{{\rm d}}{{\rm d}t})$ that is dense in ${\mathscr{G}}(\frac{{\rm d}}{{\rm d}t})$.
Having found $f$ continuously dif\/ferentiable and vanishing outside $[-N,N]$ for some positive $N$, we use
the mapping $\varphi$ on $[0,1]$ to $[-N,N]$ def\/ined by $\varphi(x)=2Nx-N$, for each~$x$ in $[0,1]$, to
transform $[0,1]$ onto $[-N,N]$.
Then $f\circ\varphi$ vanishes outside $[0,1]$ and is continuously dif\/ferentiable on~${\mathbb{R}}$.
We def\/ine~$B_n(f)$ as~$B_n(f\circ\varphi)\circ\varphi^{-1}$, where $B_n(h)$ for a~function $h$ def\/ined
on $[0,1]$ is as described in the following def\/inition.
(What follows, through the proof of Theorem~\ref{thm: Bernstein}, is our elaboration and completion of
a~few remarks of the great twentieth century classical analyst and leading expert on trigonometric
polynomials, Antoni Zygmund, during a~course of lectures on ``Approximation of functions'' at the
University of Chicago in 1948.)

\begin{Definition}
With $f$ a~real-valued function def\/ined and bounded on the interval $[0,1]$, let~$B_n(f)$ be the
polynomial on $[0,1]$ that assigns to $x$ the value
\begin{gather*}
\sum_{k=0}^n\binom{n}{k}x^k(1-x)^{n-k}f\left(\frac{k}{n}\right).
\end{gather*}
{}$B_n(f)$ is the $n$th Bernstein polynomial for~$f$.
\end{Definition}

The following identities will be useful to us in the proof of Theorem~\ref{thm: Bernstein}.
\begin{gather}
B_n(1)=\sum_{k=0}^n\binom{n}{k}x^k(1-x)^{n-k}=1,
\nonumber
\\
B_n(x)=\sum_{k=0}^n\binom{n}{k}\frac{k}{n}x^k(1-x)^{n-k}=x,
\label{equation: 3}
\\
B_n\big(x^2\big)=\sum_{k=0}^n\binom{n}{k}\frac{k^2}{n^2}x^k(1-x)^{n-k}=\frac{(n-1)x^2}{n}+\frac{x}{n},
\label{equation: 4}
\\
B_n\big(x^3\big)=\sum_{k=0}^n\binom{n}{k}\frac{k^3}{n^3}x^k(1-x)^{n-k}=\frac{(n-1)(n-2)x^3}{n^2}
+\frac{3(n-1)x^2}{n^2}+\frac{x}{n^2},
\label{equation: 5}
\\
B_n\big(x^4\big)=\sum_{k=0}^n\binom{n}{k}\frac{k^4}{n^4}x^k(1-x)^{n-k}=\frac{(n-1)(n-2)(n-3)x^4}{n^3}
\nonumber
\\
\phantom{B_n\big(x^4\big)=}
{}+\frac{6(n-1)(n-2)x^3}{n^3}+\frac{7(n-1)x^2}{n^3}+\frac{x}{n^3},
\label{equation: 6}
\\
\sum_{k=0}^n\binom{n}{k}\left(\frac{k}{n}-x\right)^2x^k(1-x)^{n-k}=x(1-x)\frac{1}{n},
\label{equation: 7}
\\
\sum_{k=0}^n\binom{n}{k}\left(\frac{k}{n}-x\right)^4x^k(1-x)^{n-k}=x(1-x)\frac{(3n-6)x(1-x)+1}{n^3}.
\label{equation: 8}
\end{gather}

To prove these identities, f\/irst, from the binomial theorem,
\begin{gather*}
B_n(1)=\sum_{k=0}^n\binom{n}{k}x^k(1-x)^{n-k}=[x+(1-x)]^n=1.
\end{gather*}
Note that
\begin{gather*}
\frac{{\rm d}}{{\rm d}p}\left(\sum_{k=0}^n\binom{n}{k}p^kq^{n-k}\right)=\frac{{\rm d}}{{\rm d}p}\big((p+q)^n\big)=n(p+q)^{n-1}.
\end{gather*}
Thus
\begin{gather*}
\sum_{k=0}^n\binom{n}{k}\frac{k}{n}p^kq^{n-k}=(p+q)^{n-1}p.
\end{gather*}
Replacing $p$ by $x$ and $q$ by $1-x$ in the above expression, we have identity~\eqref{equation: 3}.
Now, dif\/ferentiating this expression with respect to $p$ three more times and each time multiplying both
sides of the result by $\frac{p}{n}$, we have the following
\begin{gather*}
\sum_{k=0}^n\binom{n}{k}\frac{k^2}{n^2}p^kq^{n-k}=\frac{(n-1)(p+q)^{n-2}}{n}p^2+\frac{(p+q)^{n-1}}{n}p,
\\[1.11mm]
\sum_{k=0}^n\binom{n}{k}\frac{k^3}{n^3}p^kq^{n-k}=\frac{(n-1)(n-2)(p+q)^{n-3}}{n^2}
p^3+\frac{3(n-1)(p+q)^{n-2}}{n^2}p^2+\frac{(p+q)^{n-1}}{n^2}p,
\\[1.11mm]
\sum_{k=0}^n\binom{n}{k}\frac{k^4}{n^4}p^kq^{n-k}=\frac{(n-1)(n-2)(n-3)(p+q)^{n-4}}{n^3}
p^4+\frac{6(n-1)(n-2)(p+q)^{n-3}}{n^3}p^3
\\
\phantom{\sum_{k=0}^n\binom{n}{k}\frac{k^4}{n^4}p^kq^{n-k}=}
{}+\frac{7(n-1)(p+q)^{n-2}}{n^3}p^2+\frac{(p+q)^{n-1}}{n^3}p.
\end{gather*}
Replacing $p$ by $x$ and $q$ by $1-x$ in the above three identities, we obtain the
identities~\eqref{equation: 4},~\eqref{equation: 5} and~\eqref{equation: 6}.
It follows that
\begin{gather*}
\sum_{k=0}^n\binom{n}{k}\left(\frac{k}{n}-x\right)^2x^k(1-x)^{n-k}
=\left[\frac{(n-1)x^2}{n}+\frac{x}{n}\right]-2x^2+x^2=x(1-x)\frac{1}{n},
\end{gather*}
and
\begin{gather*}
\sum_{k=0}^n\binom{n}{k}\left(\frac{k}{n}-x\right)^4x^k(1-x)^{n-k}
\\
\qquad=\left[\frac{(n-1)(n-2)(n-3)x^4}{n^3}+\frac{6(n-1)(n-2)x^3}{n^3}+\frac{7(n-1)x^2}{n^3}+\frac{x}{n^3}\right]
\\
\qquad\phantom{=}
{}-4x\left[\frac{(n-1)(n-2)x^3}{n^2}+\frac{3(n-1)x^2}{n^2}+\frac{x}{n^2}\right]
+6x^2\left[\frac{(n-1)x^2}{n}+\frac{x}{n}\right]-4x^4+x^4
\\
\qquad=x(1-x)\frac{(3n-6)x(1-x)+1}{n^3}.
\end{gather*}
\begin{Theorem}\label{thm: Bernstein}
Let $f$ be a~real-valued function defined, and bounded by $M$ on the interval~$[0,1]$.
For each point~$x$ of continuity of~$f$, $B_n(f)(x)\rightarrow f(x)$ as $n\rightarrow \infty$.
If~$f$ is continuous on~$[0,1]$, then the Bernstein polynomial $B_n(f)$ tends uniformly to~$f$ as
$n\rightarrow \infty$.
With $x$ a~point of differentiability of~$f$, $B_n'(f)(x)\rightarrow f'(x)$ as $n\rightarrow \infty$.
If~$f$ is continuously differentiable on~$[0,1]$, then $B_n'(f)$ tends to $f'$ uniformly as $n\rightarrow
\infty$.
\end{Theorem}

\begin{proof}
From
\begin{gather*}
B_n(f)(x)-f(x)=\sum_{k=0}^n\binom{n}{k}x^k(1-x)^{n-k}f\left(\frac{k}{n}\right)-f(x)\sum_{k=0}
^n\binom{n}{k}x^k(1-x)^{n-k}
\\
\phantom{B_n(f)(x)-f(x)}
=\sum_{k=0}^n\binom{n}{k}x^k(1-x)^{n-k}\left[f\left(\frac{k}{n}\right)-f(x)\right],
\end{gather*}
it follows that, for each $x$ in $[0,1]$,
\begin{gather*}
|B_n(f)(x)-f(x)|\leq\sum_{k=0}^n\binom{n}{k}x^k(1-x)^{n-k}\left|f\left(\frac{k}{n}\right)-f(x)\right|.
\end{gather*}

To estimate this last sum, we separate the terms into two sums $\sum'$ and $\sum''$, those where
$|\frac{k}{n}-x|$ is less than a~given positive $\delta$ and the remaining terms, those for which
$\delta\leq|\frac{k}{n}-x|$.
Suppose that $x$ is a~point of continuity of~$f$.
Then for any $\varepsilon>0$, there is a~positive $\delta$ such that $|f(x')-f(x)|<\frac{\varepsilon}{2}$
when $|x'-x|<\delta$.
For the f\/irst sum,
\begin{gather*}
\sideset{}{'}\sum\binom{n}{k}x^k(1-x)^{n-k}\left|f\left(\frac{k}{n}\right)-f(x)\right|\\
\qquad{}
<\sideset{}{'}\sum\binom{n}{k}x^k(1-x)^{n-k}\frac{\varepsilon}{2}
\leq\frac{\varepsilon}{2}\sum_{k=0}^n\binom{n}{k}x^k(1-x)^{n-k}=\frac{\varepsilon}{2}.
\end{gather*}
For the remaining terms, we have $\delta^2\leq|\frac{k}{n}-x|^2$,
\begin{gather*}
\delta^2\sideset{}{''}\sum\binom{n}{k}x^k(1-x)^{n-k}\left|f\left(\frac{k}{n}
\right)-f(x)\right|
\\
\qquad
\leq\sideset{}{''}\sum\binom{n}{k}\left(\frac{k}{n}-x\right)^2x^k(1-x)^{n-k}
\left|f\left(\frac{k}{n}\right)-f(x)\right|
\\
\qquad\leq\sideset{}{''}\sum\binom{n}{k}\left(\frac{k}{n}-x\right)^2x^k(1-x)^{n-k}2M
\\
\qquad{} \overset{\text{from~\eqref{equation: 7}}}{\leq} 2M\sum_{k=0}^n\binom{n}{k}\left(\frac{k}{n}-x\right)^2x^k(1-x)^{n-k}
=2M\frac{x(1-x)}{n}\leq\frac{2M}{n}.
\end{gather*}
Thus
\begin{gather*}
\sideset{}{''}\sum\binom{n}{k}x^k(1-x)^{n-k}\left|f\left(\frac{k}{n}\right)-f(x)\right|\leq\frac{2M}{\delta^2n}
.
\end{gather*}
For this $\delta$, we can choose $n_0$ large enough so that, when $n\geq n_0$,
$\frac{2M}{\delta^2n}<\frac{\varepsilon}{2}$.
For such an $n$ and the given $x$
\begin{gather*}
|B_n(f)(x)-f(x)|\leq\sideset{}{'}\sum+\sideset{}{''}\sum<\frac{\varepsilon}{2}+\frac{\varepsilon}{2}
=\varepsilon.
\end{gather*}
Hence $B_n(f)(x)\rightarrow f(x)$ as $n\rightarrow \infty$ for each point $x$ of continuity of the function
$f$.
If~$f$ is continuous at each point of~$[0,1]$, then it is uniformly continuous on $[0,1]$, and for this
given~$\varepsilon$, we can choose $\delta$ so that $|f(x')-f(x)|<\frac{\varepsilon}{2}$ for each pair of
points~$x'$ and~$x$ in $[0,1]$ such that $|x'-x|<\delta$.
From the preceding argument, with $n_0$ chosen for this $\delta$, and when $n\geq n_0$,
$|B_n(f)(x)-f(x)|<\varepsilon$ for each $x$ in $[0,1]$.
Thus $\|B_n(f)-f\|\leq\varepsilon$, and $B_n(f)$ tends uniformly to~$f$ as $n\rightarrow \infty$.

Now, with $x$ in $[0,1]$,
\begin{gather*}
B_n'(f)=\frac{{\rm d}}{{\rm d}x}\bigg(\sum_{k=0}^n\binom{n}{k}x^k(1-x)^{n-k}f\left(\frac{k}{n}\right)\bigg)
\\
\phantom{B_n'(f)}
=\sum_{k=1}^n\binom{n}{k}kx^{k-1}(1-x)^{n-k}f\left(\frac{k}{n}\right)
\\
\phantom{B_n'(f)=}
{}-\sum_{k=0}^{n-1}\binom{n}{k}(n-k)x^k(1-x)^{n-k-1}f\left(\frac{k}{n}\right)+\big[nx^{n-1}f(1)-n(1-x)^{n-1}f(0)\big]
\\
\phantom{B_n'(f)}
=\sum_{k=0}^n\binom{n}{k}\big[k(1-x)-(n-k)x\big]x^{k-1}(1-x)^{n-k-1}f\left(\frac{k}{n}\right)
\\
\phantom{B_n'(f)}
=n\sum_{k=0}^n\binom{n}{k}\left(\frac{k}{n}-x\right)x^{k-1}(1-x)^{n-k-1}f\left(\frac{k}{n}\right).
\end{gather*}
(Note that $\big(\frac{k}{n}-x\big)x^{k-1}=-1$ when $k=0$ and $\big(\frac{k}{n}-x\big)(1-x)^{n-k-1}=1$ when
$k=n$.) Also,
\begin{gather*}
0=f(x)\frac{{\rm d}}{{\rm d}x}(1)=f(x)\frac{{\rm d}}{{\rm d}x}\bigg(\sum_{k=0}^n\binom{n}{k}x^k(1-x)^{n-k}\bigg)
\\
\phantom{0=f(x)\frac{{\rm d}}{{\rm d}x}(1)}
=f(x)n\sum_{k=0}^n\binom{n}{k}\left(\frac{k}{n}-x\right)x^{k-1}(1-x)^{n-k-1}.
\end{gather*}
Thus
\begin{gather*}
B_n'(f)(x)
=n\sum_{k=0}^n\binom{n}{k}\left(\frac{k}{n}-x\right)x^{k-1}(1-x)^{n-k-1}\left[f\left(\frac{k}{n}\right)-f(x)\right]
\end{gather*}
for all $x$ in $[0,1]$.

Suppose that $x$ is a~point of dif\/ferentiability of~$f$.
Let a~positive $\varepsilon$ be given.
We write
\begin{gather*}
\frac{f\big(\frac{k}{n}\big)-f(x)}{\frac{k}{n}-x}=f'(x)+\xi_k.
\end{gather*}
From the assumption of dif\/ferentiability of~$f$ at $x$, there is a~positive $\delta$ such that, when
$0<|x'-x|<\delta$, $|\frac{f(x')-f(x)}{x'-x}-f'(x)|<\frac{\varepsilon}{2}$.
Thus, when $0<|\frac{k}{n}-x|<\delta$,
\begin{gather*}
|\xi_k|=\left|\frac{f\big(\frac{k}{n}\big)-f(x)}{\frac{k}{n}-x}-f'(x)\right|<\frac{\varepsilon}{2}.
\end{gather*}
If~$\frac{k}{n}$ happens to be $x$ for some $k$, we def\/ine $\xi_k$ to be $0$ for that $k$ and note that
the inequality just stated, when $|\frac{k}{n}-x|>0$, remains valid when $\frac{k}{n}=x$.
It follows that
\begin{gather*}
B_n'(f)(x)
=n\sum_{k=0}^n\binom{n}{k}\left(\frac{k}{n}-x\right)x^{k-1}(1-x)^{n-k-1}\left[f\left(\frac{k}
{n}\right)-f(x)\right]
\\
\phantom{B_n'(f)(x)}
=n\sum_{k=0}^n\binom{n}{k}\left(\frac{k}{n}-x\right)x^{k-1}(1-x)^{n-k-1}\left[\left(\frac{k}{n}
-x\right)f'(x)+\left(\frac{k}{n}-x\right)\xi_k\right]
\\
\phantom{B_n'(f)(x)}
=f'(x)n\sum_{k=0}^n\binom{n}{k}\left(\frac{k}{n}-x\right)^2x^{k-1}(1-x)^{n-k-1}
\\
\phantom{B_n'(f)(x)=}
{}+n\sum_{k=0}^n\binom{n}{k}\left(\frac{k}{n}-x\right)^2x^{k-1}(1-x)^{n-k-1}\xi_k
\\
\phantom{B_n'(f)(x)}
=f'(x)+n\sum_{k=0}^n\binom{n}{k}\left(\frac{k}{n}-x\right)^2x^{k-1}(1-x)^{n-k-1}\xi_k.
\end{gather*}
For the last equality we made use of~\eqref{equation: 7}.
We estimate this last sum by separating it, again, into the two sums $\sum'$ and $\sum''$, those with the
$k$ for which $|\frac{k}{n}-x|<\delta$ and those for which $\delta\leq|\frac{k}{n}-x|$, respectively.
For the f\/irst sum, we have
\begin{gather*}
\left|n\sideset{}{'}\sum\right|
\leq n\sideset{}{'}\sum\binom{n}{k}\left(\frac{k}{n}-x\right)^2x^{k-1}(1-x)^{n-k-1}|\xi_k|
\\
\phantom{\left|n\sideset{}{'}\sum\right|}
<n\sideset{}{'}\sum\binom{n}{k}\left(\frac{k}{n}-x\right)^2x^{k-1}(1-x)^{n-k-1}\frac{\varepsilon}{2}
\\
\phantom{\left|n\sideset{}{'}\sum\right|}
\leq\frac{\varepsilon}{2}n\sum_{k=0}^n\binom{n}{k}\left(\frac{k}{n}-x\right)^2x^{k-1}(1-x)^{n-k-1}=\frac{\varepsilon}{2}
\end{gather*}
from~\eqref{equation: 7} and the choice of~$\delta$ (that is, the dif\/ferentiability of~$f$ at $x$).
For the second sum, we have that $\delta\leq|\frac{k}{n}-x|$ so that
\begin{gather*}
|\xi_k|\leq\left|\frac{f\big(\frac{k}{n}\big)-f(x)}{\frac{k}{n}-x}\right|+|f'(x)|\leq\frac{2M}{\delta}+|f'(x)|
\end{gather*}
and ($\delta^2\leq|\frac{k}{n}-x|^2$)  
\begin{gather*}
\delta^2\left|n\sideset{}{''}\sum\right|
\leq\delta^2n\sideset{}{''}\sum\binom{n}{k}\left(\frac{k}{n}-x\right)^2x^{k-1}(1-x)^{n-k-1}|\xi_k|
\\
\phantom{\delta^2\left|n\sideset{}{''}\sum\right|}
\leq n\sideset{}{''}\sum\binom{n}{k}\left(\frac{k}{n}-x\right)^4x^{k-1}(1-x)^{n-k-1}|\xi_k|
\\
\phantom{\delta^2\left|n\sideset{}{''}\sum\right|}
\leq n\sideset{}{''}\sum\binom{n}{k}\left(\frac{k}{n}-x\right)^4x^{k-1}(1-x)^{n-k-1}\left(\frac{2M}
{\delta}+|f'(x)|\right)
\\
\phantom{\delta^2\left|n\sideset{}{''}\sum\right|}
\leq n\sum_{k=0}^n\binom{n}{k}\left(\frac{k}{n}-x\right)^4x^{k-1}(1-x)^{n-k-1}\left(\frac{2M}{\delta}
+|f'(x)|\right)
\\
\phantom{\delta^2\left|n\sideset{}{''}\sum\right|}
\overset{\text{from~\eqref{equation: 8}}}{=} n\frac{(3n-6)x(1-x)+1}{n^3}\left(\frac{2M}{\delta}+|f'(x)|\right)
\\
\phantom{\delta^2\left|n\sideset{}{''}\sum\right|}
\leq n\frac{3}{n^2}\left(\frac{2M}{\delta}+|f'(x)|\right)=\frac{6M+3\delta|f'(x)|}{n\delta}.
\end{gather*}
Thus
\begin{gather*}
\left|n\sideset{}{''}\sum\right|\leq\frac{6M+3\delta|f'(x)|}{n\delta^3}.
\end{gather*}
For this $\delta$, we can choose $n_0$ large enough so that, when $n\geq n_0$,
\begin{gather*}
\frac{6M+3\delta|f'(x)|}{n\delta^3}<\frac{\varepsilon}{2}.
\end{gather*}
For such $n$ and the given $x$
\begin{gather*}
\left|B_n'(f)(x)-f'(x)\right|\leq\left|n\sideset{}{'}\sum\right|+\left|n\sideset{}{''}\sum\right|<\frac{\varepsilon}
{2}+\frac{\varepsilon}{2}=\varepsilon.
\end{gather*}
Hence $B_n'(f)(x)\rightarrow f'(x)$ as $n\rightarrow \infty$ for each point $x$ of dif\/ferentiability of
the function $f$.

We show, now, that if~$f$ is continuously dif\/ferentiable on $[0,1]$, then the sequence $\{B_n'(f)\}$
tends to $f'$ uniformly.
We intercept the proof for pointwise convergence at each point of dif\/ferentiability of~$f$ at the formula:
\begin{gather*}
B_n'(f)(x)=f'(x)+n\sum_{k=0}^n\binom{n}{k}\left(\frac{k}{n}-x\right)^2x^{k-1}(1-x)^{n-k-1}\xi_k.
\end{gather*}
Assuming that $f$ is everywhere dif\/ferentiable on $[0,1]$ and $f'$ is continuous on $[0,1]$, let $M'$ be
$\sup\{|f'(x)|: x\in [0,1]\}$.
Choose $\delta$ positive and such that $|f'(x')-f'(x)|<\frac{\varepsilon}{2}$ when $|x'-x|<\delta$.
Now, for any given $x$ in $[0,1]$, recall that we had def\/ined
\begin{gather*}
\xi_k=\frac{f\big(\frac{k}{n}\big)-f(x)}{\frac{k}{n}-x}-f'(x)
\quad
\text{when}
\quad
\frac{k}{n}\neq x,
\qquad
\text{and}
\qquad
\xi_k=0
\quad
\text{when}
\quad
\frac{k}{n}=x.
\end{gather*}
From the dif\/ferentiability of~$f$ on $[0,1]$, the mean value theorem applies, and
\begin{gather*}
f\left(\frac{k}{n}\right)-f(x)=f'(x_k)\left(\frac{k}{n}-x\right),
\end{gather*}
where $x_k$ is in the open interval with endpoints $\frac{k}{n}$ and $x$, when $\frac{k}{n}\neq x$.
In case $\frac{k}{n}=x$, we may choose $f'(x_k)$ as we wish, and we choose $x$ as $x_k$.
With these choices, $\xi_k=f'(x_k)-f'(x)$.
Our formula becomes
\begin{gather*}
B_n'(f)(x)-f'(x)=n\sum_{k=0}^n\binom{n}{k}\left(\frac{k}{n}-x\right)^2x^{k-1}(1-x)^{n-k-1}
\big(f'(x_k)-f'(x)\big).
\end{gather*}
In this case when we estimate the sum in the right-hand side of this equality by separating it into the two
parts $\sum'$ and $\sum''$ exactly as we did before (for approximation of the derivatives at the single
point $x$ of dif\/ferentiability), except that in this case, $|\xi_k|$ is replaced by $|f'(x_k)-f'(x)|$ and~$\delta$ has been chosen by means of the uniform continuity of~$f'$ on $[0,1]$ such that
$|f'(x_k)-f'(x)|<\frac{\varepsilon}{2}$ when $|x_k-x|<\delta$, as is the case when $|\frac{k}{n}-x|<\delta$.
For the f\/irst sum $\sum'$, the sum over those $k$ such that $|\frac{k}{n}-x|<\delta$,
\begin{gather*}
n\sideset{}{'}\sum\binom{n}{k}\left(\frac{k}{n}-x\right)^2x^{k-1}(1-x)^{n-k-1}
|f'(x_k)-f'(x)|
\\
\qquad{}
\overset{\text{from~\eqref{equation: 7}}}{<}\frac{\varepsilon}{2}n\sum_{k=0}^n\binom{n}{k}\left(\frac{k}{n}-x\right)^2x^{k-1}(1-x)^{n-k-1}
=\frac{\varepsilon}{2}.
\end{gather*}
For the second sum $\sum''$, the sum over those $k$ such that $\delta\leq|\frac{k}{n}-x|$, again, we have
$\delta^2\leq|\frac{k}{n}-x|^2$.
This time, $|f'(x_k)-f'(x)|\leq2M'$ (and we really don't care that $x_k$ may be very close to $x$ as long
as $|\frac{k}{n}-x|\geq\delta$ in this part of the estimate),
\begin{gather*}
\delta^2n\sideset{}{''}\sum\binom{n}{k}\left(\frac{k}{n}-x\right)^2x^{k-1}(1-x)^{n-k-1}
|f'(x_k)-f'(x)|
\\
\qquad\leq n\sideset{}{''}\sum\binom{n}{k}\left(\frac{k}{n}-x\right)^4x^{k-1}(1-x)^{n-k-1}2M'
\\
\qquad\leq2M'n\sum_{k=0}^n\binom{n}{k}\left(\frac{k}{n}-x\right)^4x^{k-1}(1-x)^{n-k-1}
\\
\qquad
\overset{\text{from~\eqref{equation: 8}}}{=}2M'n\frac{(3n-6)x(1-x)+1}{n^3}
\leq\frac{6M'}{n}.
\end{gather*}
Again, for this $\delta$, we can choose $n_0$ large enough so that, when $n>n_0$
\begin{gather*}
n\sideset{}{''}\sum\binom{n}{k}\left(\frac{k}{n}-x\right)^2x^{k-1}(1-x)^{n-k-1}|f'(x_k)-f'(x)|\leq\frac{6M'}
{\delta^2n}<\frac{\varepsilon}{2},
\end{gather*}
and
\begin{gather*}
\left|B_n'(f)(x)-f'(x)\right|
\leq n\sideset{}{'}\sum\binom{n}{k}\left(\frac{k}{n}-x\right)^2x^{k-1}(1-x)^{n-k-1}|f'(x_k)-f'(x)|
\\
\qquad{}
+n\sideset{}{''}\sum\binom{n}{k}\left(\frac{k}{n}-x\right)^2x^{k-1}(1-x)^{n-k-1}|f'(x_k)-f'(x)|
<\frac{\varepsilon}{2}+\frac{\varepsilon}{2}=\varepsilon
\end{gather*}
for each $x$ in $[0,1]$.
Thus $\|B_n'(f)-f'\|\leq\varepsilon$, and $\{B_n'(f)\}$ tends to $f'$ uniformly.
\end{proof}

\section{Murray--von Neumann algebras}\label{s6}  

\subsection{Finite von Neumann algebras}

Let ${\mathcal H}$ be a~Hilbert space.
Two projections $E$ and $F$ are said to be \textit{orthogonal} if~$EF=0$.
If the range of~$F$ is contained in the range of~$E$ (equivalently, $EF=F$), we say that $F$ is
a~subprojection of~$E$ and write $F\leq E$.
Let ${\mathcal R}$ be a~von Neumann algebra acting on ${\mathcal H}$.
Suppose that $E$ and $F$ are nonzero projections in ${\mathcal R}$.
We say $E$ is a~\textit{minimal projection} in ${\mathcal R}$ if~$F\leq E$ implies $F=E$.
Murray and von Neumann conceived the idea of comparing the ``sizes'' of projections in a~von Neumann algebra
in the following way: $E$ and $F$ are said to be \textit{equivalent} (\textit{modulo} or \textit{relative
to} ${\mathcal R}$), written $E\sim F$, when $V^*V=E$ and $VV^*=F$ for some $V$ in ${\mathcal R}$.
(Such an operator $V$ is called a~\textit{partial isometry} with \textit{initial} projection $E$ and
\textit{final} projection $F$.) We write $E \precsim F$ when $E \sim F_0$ and $F_0 \leq F$ and $E\prec F$
when $E$ is, in addition, not equivalent to $F$.
It is apparent that $\sim$ is an equivalence relation on the projections in ${\mathcal R}$.
In addition, $\precsim$ is a~{\it partial ordering} of the equivalence classes of projections in ${\mathcal
R}$, and it is a~non-trivial and crucially important fact that this partial ordering is a~{\it total
ordering} when ${\mathcal R}$ is a~factor (Factors are von Neumann algebras whose centers consist of scalar
multiples of the identity operator).
Murray and von Neumann also def\/ine inf\/inite and f\/inite projections in this framework modeled on the
set-theoretic approach.
The projection $E$ in ${\mathcal R}$ is {\it infinite} (relative to ${\mathcal R}$) when $E\sim F < E$, and
{\it finite} otherwise.
We say that the von Neumann algebra ${\mathcal R}$ is {\it finite} when the identity operator~$I$ is
f\/inite.

\begin{Proposition}
\label{prop: IEIF}
Suppose that $E$ and $F$ are projections in a~finite von Neumann algebra ${\mathcal R}$.
If~$E\sim F$, then $I-E\sim I-F$.
\end{Proposition}

\begin{proof}
Suppose $I-E$ and $I-F$ are not equivalent.
Then there is a~central projection $P$ such that either $P(I-E)\prec P(I-F)$ or $P(I-F)\prec P(I-E)$.
Suppose $P(I-E)\sim G<P(I-F)$.
Then, since $PE\sim PF$, $P=P(I-E)+PE\sim G+PF<P(I-F)+PF=P$, contrary to the assumption that ${\mathcal R}$
is f\/inite.
The symmetric argument applies if~$P(I-F)\prec P(I-E)$.
Thus $I-E\sim I-F$.
\end{proof}
\begin{Proposition}\label{prop: Kaplansky}
For any projections $E$ and $F$ in a~finite von Neumann algebra ${\mathcal R}$,
\begin{gather*}
\Delta(E\vee F)+\Delta(E\wedge F)=\Delta(E)+\Delta(F),
\end{gather*}
where $\Delta$ is the center-valued dimension function on ${\mathcal R}$.
\end{Proposition}
\begin{proof}
Since $E\vee F-F\sim E-E\wedge F$ (Kaplansky formula), we have
\begin{gather*}
\Delta(E\vee F)-\Delta(F)=\Delta(E\vee F-F)=\Delta(E-E\wedge F)=\Delta(E)-\Delta(E\wedge F).
\end{gather*}
Thus $\Delta(E\vee F)+\Delta(E\wedge F)=\Delta(E)+\Delta(F)$.
\end{proof}
\begin{Proposition}
\label{prop: sequence}
Suppose that $E$, $F$, and $G$ are projections in a~f\/inite von Neumann algebra~${\mathcal R}$, and $E$
and~$F$ are the $($strong-operator$)$ limits of increasing nets $\{E_a\}$ and $\{F_a\}$, respectively, of
projections in ${\mathcal R}$ $($the index set being the same$)$.
Then
\begin{enumerate}\itemsep=0pt
\item[$(i)$] $\{E_a\vee G\}$ is strong-operator convergent to $E\vee G$;

\item[$(ii)$] $\{E_a\wedge G\}$ is strong-operator convergent to $E\wedge G$;

\item[$(iii)$] $\{E_a\wedge F_a\}$ is strong-operator convergent to $E\wedge F$.
\end{enumerate}
\end{Proposition}

\begin{proof}
$(i)$ Since the net $\{E_a\vee G\}$ is increasing and bounded above by $E\vee G$, it converges to
a~projection $P$ in ${\mathcal R}$, and $P\leq E\vee G$.
For each index $a$, $E_a\leq E_a\vee G\leq P$, so $\bigvee E_a\leq P$; that is $E\leq P$.
Also, $G\leq E_a\vee G\leq P$; so $E\vee G\leq P$.
Thus $P=E\vee G$.

$(ii)$ Since the net $\{E_a\wedge G\}$ is increasing and bounded above by $E\wedge G$, it converges
to a~projection $P$ in ${\mathcal R}$, and $P\leq E\wedge G$.
Recall that the center-valued dimension function $\Delta$ on ${\mathcal R}$ is weak-operator continuous on
the set of all projections on ${\mathcal R}$; together with Proposition~\ref{prop: Kaplansky},
\begin{gather*}
\Delta(P)=\lim\Delta(E_a\wedge G)=\lim[\Delta(E_a)+\Delta(G)-\Delta(E_a\vee G)]
\\
\phantom{\Delta(P)}
=\Delta(E)+\Delta(G)-\Delta(E\vee G)=\Delta(E\wedge G).
\end{gather*}
Since $E\wedge G-P$ is a~projection in ${\mathcal R}$ and $\Delta(E\wedge G-P)=0$, it follows that
$P=E\wedge G$.

$(iii)$ The net $\{E_a\wedge F_a\}$ is increasing and therefore has a~projection $P$ as
a~strong-operator limit and least upper bound.
Since $E_a\wedge F_a \leq E\wedge F$ for each~$a$, $P\leq E\wedge F$.
With $a'$ f\/ixed the net $\{E_a\wedge F_{a'}\}$ has strong-operator limit $E\wedge F_{a'}$ from $(ii)$.
Since $E_a\wedge F_{a'}\leq E_a\wedge F_a$ when $a'\leq a$, $E\wedge F_{a'}\leq P$ for each $a'$.
Again, from $(ii)$, $\{E\wedge F_a\}$ has $E\wedge F$ as its strong-operator limit.
Thus $E\wedge F\leq P$.
Hence $P=E\wedge F$.
\end{proof}

\begin{Proposition}
\label{prop: projection range}
Let $E$ be a~projection in a~finite von Neumann algebra ${\mathcal R}$ acting on a~Hilbert space
${\mathcal H}$.
With $T$ in ${\mathcal R}$, let $F$ be the projection with range $\{x:Tx\in E({\mathcal H})\}$.
Then $F\in {\mathcal R}$ and $E\precsim F$.
\end{Proposition}
\begin{proof}
With $A'$ in ${\mathcal R}'$ and $Tx$ in $E({\mathcal H})$, $TA'x=A'Tx\in E({\mathcal H})$ since $A'E=EA'$.
Thus $F({\mathcal H})$ is stable under ${\mathcal R}'$, and $F\in {\mathcal R}''$ $({=}{\mathcal R})$.

Note that $Tx\in E({\mathcal H})$ if and only if~$(I-E)Tx=0$.
Thus $F({\mathcal H})$ is the null space of~$(I-E)T$ (that is, $F=N[(I-E)T]$).
Then $I-F=I-N[(I-E)T]=R[T^*(I-E)]\sim R[(I-E)T]\leq I-E$.
If~$E\precnsim F$, then there is a~central projection $P$ in ${\mathcal R}$ such that $PF\prec PE$.
At the same time, $P(I-F)\precsim P(I-E)$ so that $P(I-F)\sim E_0\leq P(I-E)$.
Thus $P=PF+P(I-F)\prec PE+E_0\leq PE+P(I-E)=P$.
This is contrary to the assumption that ${\mathcal R}$ is f\/inite.
It follows that $E\precsim F$.
\end{proof}

\subsection{The algebra of af\/f\/iliated operators}

Recall (Def\/inition~\ref{df affiliation}) that a~closed densely def\/ined operator $T$ is af\/f\/iliated
with a~von Neumann algebra ${\mathcal R}$ acting on a~Hilbert space ${\mathcal H}$ when $TU'=U'T$ for each
unitary operator $U'$ in ${\mathcal R}'$ (the commutant of~${\mathcal R}$).
\begin{Proposition}
\label{prop: RTRT*}
If~$T$ is affiliated with a~von Neumann algebra ${\mathcal R}$, then
\begin{enumerate}\itemsep=0pt
\item[$(i)$] $R(T)$ and $N(T)$ are in ${\mathcal R}$;

\item[$(ii)$] $R(T^*)=R(T^*T)=R((T^*T)^{1/2})$;

\item[$(iii)$] $R(T)\sim R(T^*)$ relative to ${\mathcal R}$.
\end{enumerate}
\end{Proposition}
\begin{proof}
$(i)$ From Proposition~\ref{prop: RN}, $x\in N(T)({\mathcal H})$ if and only if~$x\in{\mathscr{D}}(T)$ and $Tx=0$.
If~$U'$ is a~unitary operator in ${\mathcal R}'$, then $U'x\in {\mathscr{D}}(T)$ when
$x\in{\mathscr{D}}(T)$ and $TU'x=U'Tx$.
Thus $TU'x=0$ when $x\in N(T)({\mathcal H})$, and $N(T)({\mathcal H})$ is stable under each unitary
operator in ${\mathcal R}'$.
Hence, $N(T)\in {\mathcal R}$.
From Proposition~\ref{prop: RN}, $R(T)\in {\mathcal R}$.

$(ii)$ We show that $N((T^*T)^{1/2})=N(T^*T)$.
If~$x\in N((T^*T)^{1/2})({\mathcal H})$, then $x\in {\mathscr{D}}((T^*T)^{1/2})$ and $(T^*T)^{1/2}x=0$.
Thus $x\in {\mathscr{D}}(T^*T)$, $T^*Tx=(T^*T)^{1/2}(T^*T)^{1/2}x=0$, and $x\in N(T^*T)({\mathcal H})$.

If~$x\in N(T^*T)({\mathcal H})$, then $x\in{\mathscr{D}}(T^*T)$ and $T^*Tx=0$.
Thus $x\in {\mathscr{D}}((T^*T)^{1/2})$,
$0=\langle T^*Tx,x\rangle=\langle(T^*T)^{1/2}(T^*T)^{1/2}x,x\rangle=\|(T^*T)^{1/2}x\|^2$, and $x\in
N((T^*T)^{1/2})({\mathcal H})$.
It follows that $N((T^*T)^{1/2})=N(T^*T)$.
From Proposition~\ref{prop: RN}, $R(T^*)=R(T^*T)=R((T^*T)^{1/2})$.

$(iii)$ From Theorem~\ref{thm: polar}, $T=V(T^*T)^{1/2}$, where $V$ is a~partial isometry in
${\mathcal R}$ with initial projection $R((T^*T)^{1/2})$ and f\/inal projection~$R(T)$.
From $(ii)$, $R(T^*)=R((T^*T)^{1/2})$.
Thus $R(T)$ and $R(T^*)$ are equivalent in~${\mathcal R}$.
\end{proof}

Throughout the rest of this section, ${\mathcal R}$ denotes a~f\/inite von Neumann algebra acting on
a~Hilbert space ${\mathcal H}$, and ${\mathscr{A}}_{\rm f}({\mathcal R})$ denotes the family of operators
af\/f\/iliated with~${\mathcal R}$.
We shall show that ${\mathscr{A}}_{\rm f}({\mathcal R})$ is a~$*$ algebra (cf.~\cite{KIII,MVI}).
The hypothesis that ${\mathcal R}$ is f\/inite is crucial for the results that follow.

\begin{Proposition}
\label{prop: symmetric}
If~$S$ is a~symmetric operator affiliated with~${\mathcal R}$, then $S$ is self-adjoint.
\end{Proposition}

\begin{proof}
Since $S\in {\mathscr{A}}_{\rm f}({\mathcal R})$, $(S+iI)\in {\mathscr{A}}_{\rm f}({\mathcal R})$.
It follows that
\begin{gather*}
R(S+iI)\overset{\text{Proposition~\ref{prop: RTRT*}}}{\sim} R\big((S+iI)^*\big),
\\
I-R(S+iI) \overset{\text{Proposition~\ref{prop: IEIF}}}{\sim} I-R\big((S+iI)^*\big),
\\
I-R(S+iI)\overset{\text{Proposition~\ref{prop: RN}}}{=} N\big((S+iI)^*\big)\sim N(S+iI)=I-R((S+iI)^*).
\end{gather*}
If~$x$ is in the range of~$N(S+iI)$, then $x\in {\mathscr{D}}(S+iI)$ $({=}{\mathscr{D}}(S))$ and $Sx+ix=0$.
Since $S\subseteq S^*$, $x\in {\mathscr{D}}(S^*)$ and $Sx=S^*x$, so that
$\langle Sx,x\rangle=\langle x,S^*x\rangle=\langle x,Sx\rangle=\overline{\langle Sx,x\rangle}$ and
$0=\langle Sx+ix,x\rangle=\langle Sx,x\rangle+i\langle x,x\rangle$.
Thus $\langle x,x\rangle=0$ and $x=0$.
Hence $N(S+iI)=0$ and $N((S+iI)^*)=0$.
Similarly, $N((S-iI)^*)=0$.
From Proposition~\ref{prop:2.1.6}, $S$ is self-adjoint (for $(S\pm iI)^*=S^*\mp iI$).
\end{proof}
\begin{Proposition}\label{prop: equal}
If operators $A$ and $B$ are affiliated with~${\mathcal R}$ and $A\subseteq B$, then $A=B$.
\end{Proposition}
\begin{proof}
Let $VH$ be the polar decomposition of~$B$.
Since $A\subseteq B$,
\begin{gather*}
V^*A\subseteq V^*B=V^*VH=H=H^*\subseteq(V^*A)^*.
\end{gather*}
Thus $V^*A$ is symmetric.
If fact, $V^*A$ is af\/f\/iliated with~${\mathcal R}$.
To see this, f\/irst, $V^*A$ is densely def\/ined since ${\mathscr{D}}(V^*A)={\mathscr{D}}(A)$.
Now, suppose $\{x_n\}$ is a~sequence of vectors in ${\mathscr{D}}(V^*A)$ such that $x_n \rightarrow x$ and
$V^*Ax_n \rightarrow y$.
As $V^*$ is isometric on the range of~$A$, $\|Ax_n-Ax_m\|=\|V^*Ax_n-V^*Ax_m\|\rightarrow 0$ as
$m,n\rightarrow 0$, so that $\{Ax_n\}$ converges to some vector $z$ and $V^*Ax_n \rightarrow V^*z=y$.
But since $A$ is closed, $x\in{\mathscr{D}}(A)$ and $Ax=z$.
Thus $y=V^*z=V^*Ax$, and $V^*A$ is closed.
If~$U'$ is a~unitary operator in ${\mathcal R}'$, then $U'^*AU'=A$ so that $U'^*V^*AU'=V^*U'^*AU'=V^*A$
(since $V^*\in {\mathcal R}$).
Thus $V^*A\in{\mathscr{A}}_{\rm f}({\mathcal R})$.

From Proposition~\ref{prop: symmetric}, $V^*A$ is self-adjoint.
Since $V^*A$ is contained in $H$ and self-adjoint operators are maximal symmetric (Remark~\ref{remark:
maximal symmetric}), $V^*A=H$.
Hence $A=R(B)A=VV^*A=VH=B$.
\end{proof}
\begin{Proposition}\label{prop: hat df}
If operators $S$ and $T$ are affiliated with~${\mathcal R}$, then:
\begin{enumerate}\itemsep=0pt
\item[$(i)$] $S+T$ is densely defined, preclosed and has a~unique closed extension $S\;\hat{+}\;T$
affiliated with~${\mathcal R}$;

\item[$(ii)$] $ST$ is densely defined, preclosed and has a~unique closed extension $S\;\hat{\cdot}\;T$
affiliated with~${\mathcal R}$.
\end{enumerate}
\end{Proposition}
\begin{proof}
Let $VH$ and $WK$ be the polar decompositions of~$S$ and $T$, respectively, and let $E_n$ and~$F_n$ be the
spectral projections for $H$ and $K$, respectively, corresponding to the interval $[-n,n]$ for each
positive integer $n$.

$(i)$ From the spectral theorem, $\{E_n\}$ and $\{F_n\}$ are increasing sequences of projections
with strong-operator limit~$I$.
From Proposition~\ref{prop: sequence}, $\{E_n \wedge F_n\}$ is an increasing sequence with strong-operator
limit~$I$.
Thus $\bigcup_{n=1}^\infty(E_n \wedge F_n)({\mathcal H})$ is dense in ${\mathcal H}$.
If~$x\in (E_n \wedge F_n)({\mathcal H})$, then $x\in {\mathscr{D}}(H)\cap{\mathscr{D}}(K)$.
Hence $x\in {\mathscr{D}}(S+T)$.
It follows that $S+T$ is densely def\/ined.

Since $S$ and $T$ are af\/f\/iliated with~${\mathcal R}$, $S^*$ and $T^*$ are af\/f\/iliated with~${\mathcal R}$.
From what we just proved, $S^*+T^*$ is densely def\/ined.
Since $S^*+T^*\subseteq (S+T)^*$, ${\mathscr{D}}((S+T)^*)$ is dense in ${\mathcal H}$.
From Theorem~\ref{thm: preclosed}, $S+T$ is preclosed.
The closure $S\;\hat{+}\;T$ of~$S+T$ is the smallest closed extension of~$S+T$.
If~$U'$ is a~unitary operator in ${\mathcal R}'$ and $x\in {\mathscr{D}}(S+T)$, then $x\in
{\mathscr{D}}(S)$, $x\in {\mathscr{D}}(T)$, $U'x\in {\mathscr{D}}(S)$, $U'x\in {\mathscr{D}}(T)$ (recall
that a~unitary operator transforms the domain of each af\/f\/iliated operator onto the domain itself), and
\begin{gather*}
(S+T)U'x=SU'x+TU'x=U'Sx+U'Tx=U'(S+T)x.
\end{gather*}
From Remark~\ref{remark: check core}, $S\;\hat{+}\;T\in {\mathscr{A}}_{\rm f}({\mathcal R})$ since
${\mathscr{D}}(S+T)$ is a~core for $S\;\hat{+}\;T$.
If~$A$ is a~closed extension of~$(S+T)$ and $A\in {\mathscr{A}}_{\rm f}({\mathcal R})$, then
$S\;\hat{+}\;T\subseteq A$ and, from Proposition~\ref{prop: equal}, $S\;\hat{+}\;T=A$.
Therefore, $S\;\hat{+}\;T$ is the only closed extension of~$S+T$ af\/f\/iliated with~${\mathcal R}$.

$(ii)$ By choice of~$F_n$, $KF_n$ is a~bounded, everywhere-def\/ined, self-adjoint operator in~${\mathcal R}$.
Let $T_n=TF_n$.
Then $T_n$ $({=}TF_n=WKF_n)$ is a~bounded, everywhere-def\/ined, operator in~${\mathcal R}$.
From Proposition~\ref{prop: projection range}, the projection $M_n$ with range $\{x: T_nx\in E_n({\mathcal
H})\}$ is in ${\mathcal R}$ and $E_n\precsim M_n$.
Since~$\{E_n\}$ is an increasing sequence of projections with strong-operator limit~$I$,
$\Delta(E_n)=\tau(E_n)\uparrow \tau(I)=I$ in the strong-operator topology, where~$\Delta$ is the
center-valued dimension function and~$\tau$ is the center-valued trace on ${\mathcal R}$.
Since $\{M_n\}$ is an increasing sequence and $\tau (E_n)\leq\tau(M_n)$, $\tau (M_n)\uparrow I$.
Hence $\{M_n\}$ has strong-operator limit~$I$.
From Proposition~\ref{prop: sequence}, $\{G_n\}=\{F_n\wedge M_n\}$ is an increasing sequence with
strong-operator limit~$I$.
It follows that $\bigcup_{n=1}^\infty G_n({\mathcal H})$ is dense in ${\mathcal H}$.
If~$x\in G_n({\mathcal H})$, then $T_nx\in E_n({\mathcal H})$ so that $T_nx\in
{\mathscr{D}}(H)={\mathscr{D}}(S)$.
At the same time, $x\in F_n({\mathcal H})$ so that $x\in {\mathscr{D}}(K)={\mathscr{D}}(T)$ and
$Tx=TF_nx=T_nx$.
Thus $x\in {\mathscr{D}}(ST)$.
It follows that $ST$ is densely def\/ined.

Now, $T^*S^*$ is densely def\/ined since $S^*$ and $T^*$ are in ${\mathscr{A}}_{\rm f}({\mathcal R})$.
Note that $T^*S^*\subseteq(ST)^*$, thus $(ST)^*$ is densely def\/ined.
From Theorem~\ref{thm: preclosed}, $ST$ is preclosed.
The closure $S\;\hat{\cdot}\;T$ of~$ST$ is the smallest closed extension of~$ST$.
If~$U'$ is a~unitary operator in ${\mathcal R}'$ and $x\in {\mathscr{D}}(ST)$, then $x\in
{\mathscr{D}}(T)$, $Tx\in {\mathscr{D}}(S)$, $U'x\in {\mathscr{D}}(T)$, $TU'x=U'Tx\in {\mathscr{D}}(S)$, and
\begin{gather*}
STU'x=SU'Tx=U'STx.
\end{gather*}
As with $S\;\hat{+}\;T$ in $(i)$, $S\;\hat{\cdot}\;T\in {\mathscr{A}}_{\rm f}({\mathcal R})$ and
$S\;\hat{\cdot}\;T$ is the only closed extension of~$ST$ af\/f\/iliated with~${\mathcal R}$.
\end{proof}

\begin{Proposition}
If operators $A$, $B$ and $C$ are affiliated with~${\mathcal R}$, then
\begin{gather*}
(A\;\hat{+}\;B)\;\hat{+}\;C=A\;\hat{+}\;(B\;\hat{+}\;C),
\end{gather*}
that is, the associative law holds under the addition $\hat{+}$ described in Proposition~{\rm \ref{prop: hat df}}.
\end{Proposition}
\begin{proof}
First, we note that $(A\;\hat{+}\;B)\;\hat{+}\;C$ and $A\;\hat{+}\;(B\;\hat{+}\;C)$ are closed extensions
of~$(A+B)+C$ and $A+(B+C)$, respectively.
Hence, both are preclosed.
In addition, both are densely def\/ined with domain
${\mathscr{D}}$ $({=}{\mathscr{D}}(A)\cap{\mathscr{D}}(B)\cap{\mathscr{D}}(C))$ since $A$, $B$, and $C$ are
af\/f\/iliated with~${\mathcal R}$ (see Proposition~\ref{prop: hat df}$(i)$).

Note, also, that $(A+B)+C=A+ (B+C)$ (on ${\mathscr{D}}$).
Hence, $\overline{(A+B)+C}=\overline{A+ (B+C)}$.
Moreover, ${\mathscr{D}}$ is a~core for both.
Once we note that these closed, densely def\/ined (equal) operators are af\/f\/iliated with~${\mathcal R}$,
their closed extensions $(A\;\hat{+}\;B)\;\hat{+}\;C$ and $A\;\hat{+}\;(B\;\hat{+}\;C)$ are equal to each
of them (Proposition~\ref{prop: equal}), hence, to each other, which is what we wish to prove.
To establish this af\/f\/iliation, let $U'$ be a~unitary operator in ${\mathcal R}'$ and $x$ a~vector in
${\mathscr{D}}$.
As $A$, $B$, and $C$ are in ${\mathscr{A}}_{\rm f}({\mathcal R})$, $AU'x=U'Ax$, $BU'x=U'Bx$, and
$CU'x=U'Cx$.
Thus $((A+B)+C)U'x=U'((A+B)+C)x$, and from Remark~\ref{remark: check core},
$\overline{(A+B)+C}\in{\mathscr{A}}_{\rm f}({\mathcal R})$.
Thus, $\overline{A+(B+C)}\in{\mathscr{A}}_{\rm f}({\mathcal R})$ and
$(A\;\hat{+}\;B)\;\hat{+}\;C=\overline{(A+B)+C}=\overline{A+(B+C)}=A\;\hat{+}\;(B\;\hat{+}\;C)$.
\end{proof}

For the addition operation described in Proposition~\ref{prop: hat df}, one can also show that
$A\;\hat{+}\;B=B\;\hat{+}\;A$, $A\;\hat{+}\;0=A$, and $A\;\hat{+}\;(-A)=0$, for $A,B$ af\/f\/iliated with~${\mathcal R}$.
(Note, in $A\;\hat{+}\;(-A)=0$, f\/irst, the 0-operator is def\/ined on a~dense domain ${\mathscr{D}}(A)$.
Since it is bounded, it has a~unique extension to the 0-operator on the whole Hilbert space ${\mathcal H}$.)
\begin{Proposition}
\label{prop: associative law}
If operators $A$, $B$ and $C$ are affiliated with~${\mathcal R}$, then
\begin{gather*}
(A\;\hat{\cdot}\;B)\;\hat{\cdot}\;C=A\;\hat{\cdot}\;(B\;\hat{\cdot}\;C),
\end{gather*}
that is, the associative law holds under the multiplication $\hat{\cdot}$ described in
Proposition~{\rm \ref{prop: hat df}}.
\end{Proposition}
\begin{proof}
First, we note that
\begin{gather*}
(A\cdot B)\cdot C\subseteq(A\;\hat{\cdot}\;B)\;\hat{\cdot}\;C
\qquad
\text{and}
\qquad
A\cdot(B\cdot C)\subseteq A\;\hat{\cdot}\;(B\;\hat{\cdot}\;C),
\end{gather*}
where ``$\cdot$'' is the usual multiplication of operators, hence, $(A\cdot B)\cdot C$ and $A\cdot (B\cdot
C)$ are preclosed.
Note, also, that $(A\cdot B)\cdot C=A\cdot (B\cdot C)$ on ${\mathscr{D}}={\mathscr{D}}((A\cdot B)\cdot
C)$ $({=}{\mathscr{D}}(A\cdot (B\cdot C)))$.
We shall show that the operator $A\cdot (B\cdot C)$ is densely def\/ined and its closure, $\overline{A\cdot
(B\cdot C)}$ $({=}\overline{(A\cdot B)\cdot C})$, is af\/f\/iliated with~${\mathcal R}$.
Then from Proposition~\ref{prop: equal},
$(A\;\hat{\cdot}\;B)\;\hat{\cdot}\;C=A\;\hat{\cdot}\;(B\;\hat{\cdot}\;C)$.

Let $V_1H_1$, $V_2H_2$ and $V_3H_3$ be the polar decompositions of~$A$, $B$ and $C$, respectively.
Let $E_n$, $F_n$ and $G_n$ be the spectral projections for $H_1$, $H_2$ and $H_3$, respectively,
corresponding to the interval $[-n,n]$ for each positive integer $n$.
We note that the operator $CG_n$ $({=}V_3H_3G_n)$, denoted by $C_n$, is a~bounded, everywhere-def\/ined operator.
Let $J_n$ be the projection with range $G_n({\mathcal H})\cap \{x: C_nx\in F_n({\mathcal H})\}$.
As in the proof of Proposition~\ref{prop: hat df}$(ii)$, $\{J_n\}$ is an increasing sequence with
strong-operator limit $I$.
Thus $\bigcup _{n=1}^\infty J_n({\mathcal H})$ is dense in ${\mathcal H}$.
If~$x\in J_n({\mathcal H})$, then $C_nx\in F_n({\mathcal H})$ so that $C_nx\in
{\mathscr{D}}(H_2)={\mathscr{D}}(B)$.
At the same time, $x\in G_n({\mathcal H})$ so that $x\in {\mathscr{D}}(H_3)={\mathscr{D}}(C)$ and
$Cx=CG_nx=C_nx$.
Thus $x\in{\mathscr{D}}(BC)$.
Let $B_n=(BC)J_n$.
By our def\/inition of~$J_n$, $J_n\leq G_n$ so that $CJ_n$ is a~bounded, everywhere-def\/ined operator in
${\mathcal R}$ and
\begin{gather*}
CJ_n({\mathcal H})=(CG_n)J_n({\mathcal H})=C_nJ_n({\mathcal H})\subseteq F_n({\mathcal H}).
\end{gather*}
It follows that $B_n=(BC)J_n=B(CJ_n)$ is a~bounded, everywhere-def\/ined operator in ${\mathcal R}$.
Let~$K_n$ be the projection with range $J_n({\mathcal H})\cap \{x: B_nx\in E_n({\mathcal H})\}$.
Similarly, $\{K_n\}$ is an increasing sequence with strong-operator limit $I$.
Thus $\bigcup _{n=1}^\infty K_n({\mathcal H})$ is dense in ${\mathcal H}$.
If~$x\in K_n({\mathcal H})$, then $B_nx\in E_n({\mathcal H})$ so that $B_nx\in
{\mathscr{D}}(H_1)={\mathscr{D}}(A)$.
At the same time, $x\in J_n({\mathcal H})$ so that $x\in {\mathscr{D}}(BC)$ and $BCx=BCJ_nx=B_nx$.
Thus $x\in {\mathscr{D}}(A\cdot (B\cdot C))$.
It follows that $A\cdot (B\cdot C)$  $({=}(A\cdot B)\cdot C)$ is densely def\/ined.

Now, we show that the closure $\overline{A\cdot (B\cdot C)}$ is af\/f\/iliated with~${\mathcal R}$, which
completes the proof.
If~$U'$ is a~unitary operator in ${\mathcal R}'$ and $x\in {\mathscr{D}}$ $({=}{\mathscr{D}}(A\cdot (B\cdot
C)))$, since $A$, $B$, and $C$ are af\/f\/iliated with~${\mathcal R}$, we have $A\cdot (B\cdot C)\cdot
U'x=A\cdot U'\cdot (B\cdot C)x=U'\cdot A\cdot (B\cdot C)x$.
From Remark~\ref{remark: check core}, $\overline{A\cdot (B\cdot C)}$ is af\/f\/iliated with~${\mathcal R}$
since ${\mathscr{D}}$ is a~core for $\overline{A\cdot (B\cdot C)}$.
\end{proof}
\begin{Proposition}
If operators $A$, $B$ and $C$ are affiliated with~${\mathcal R}$, then
\begin{gather*}
(A\;\hat{+}\;B)\;\hat{\cdot}\;C=A\;\hat{\cdot}\;C\;\hat{+}\;(B\;\hat{\cdot}\;C)
\qquad
\text{and}
\qquad
C\;\hat{\cdot}\;(A\;\hat{+}\;B)=C\;\hat{\cdot}\;A\;\hat{+}\;(C\;\hat{\cdot}\;B),
\end{gather*}
that is, the distributive laws hold under the addition $\hat{+}$ and multiplication $\hat{\cdot}$ described
in Proposition~{\rm \ref{prop: hat df}}.
\end{Proposition}

\begin{proof}
First, we note the following
\begin{gather*}
(A+B)C\subseteq(A\;\hat{+}\;B)\;\hat{\cdot}\;C,
\qquad
AC+BC\subseteq A\;\hat{\cdot}\;C\;\hat{+}\;(B\;\hat{\cdot}\;C),
\\[1.11mm]
C(A+B)\subseteq C\;\hat{\cdot}\;(A\;\hat{+}\;B),
\qquad
CA+AB\subseteq C\;\hat{\cdot}\;A\;\hat{+}\;(C\;\hat{\cdot}\;B),
\end{gather*}
and
\begin{gather*}
(A+B)C=AC+BC,
\qquad
CA+CB\subseteq C(A+B).
\end{gather*}
Hence, $(A+B)C$ and $CA+CB$ are preclosed.
We shall show that $(A+B)C$ and $CA+CB$ are densely def\/ined and their closures are af\/f\/iliated with~${\mathcal R}$.
Then, again, using Proposition~\ref{prop: equal}, we obtain
$(A\;\hat{+}\;B)\;\hat{\cdot}\;C=A\;\hat{\cdot}\;C\;\hat{+}\;(B\;\hat{\cdot}\;C)$ and
$C\;\hat{\cdot}\;(A\;\hat{+}\;B)=C\;\hat{\cdot}\;A\;\hat{+}\;(C\;\hat{\cdot}\;B)$.

We def\/ine $V_1H_1$, $V_2H_2$, $V_3H_3$ and $E_n$, $F_n$, $G_n$ as in the proof of Proposition~\ref{prop:
associative law}.
By choice of~$G_n$, the operator $C_n=CG_n=V_3H_3G_n$ is a~bounded and everywhere-def\/ined.
Let~$J_n$ be the projection on the range $G_n({\mathcal H})\cap\{x: C_nx \in (E_n\wedge F_n)({\mathcal
H})\}$.
Then $\bigcup_{n=1}^\infty J_n({\mathcal H})$ is dense in~${\mathcal H}$ since~$\{J_n\}$ is an increasing
sequence with strong-operator limit $I$.
If~$x\in J_n({\mathcal H})$, then $C_nx\in (E_n\wedge F_n)({\mathcal H})$ so that $C_nx\in
{\mathscr{D}}(A+B)$.
At the same time, $x\in G_n({\mathcal H})$ so that $x\in {\mathscr{D}}(H_3)={\mathscr{D}}(C)$ and
$Cx=CG_nx=C_nx$.
Thus $x\in {\mathscr{D}}((A+B)C)$.
It follows that $(A+B)C$ is densely def\/ined.

Let $A_n=AE_n$ and $B_n=BF_n$.
Then $A_n$ and $B_n$ are bounded, everywhere-def\/ined operators in ${\mathcal R}$.
Let $K_n$ be the projection on the range
\begin{gather*}
E_n({\mathcal H})\cap\{x:A_nx\in G_n({\mathcal H})\}\cap F_n({\mathcal H})\cap\{x:B_nx\in G_n({\mathcal H}
)\}.
\end{gather*}
Again, $\{K_n\}$ is an increasing sequence with strong-operator limit $I$ so that $\bigcup_{n=1}^\infty
K_n({\mathcal H})$ is dense in ${\mathcal H}$.
If~$x\in K_n({\mathcal H})$, then $A_nx\in G_n({\mathcal H})$ and $B_nx\in G_n({\mathcal H})$ so that
$A_nx\in {\mathscr{D}}(C)$ and $B_nx\in {\mathscr{D}}(C)$.
At the same time, $x\in E_n({\mathcal H})$ and $x\in F_n({\mathcal H})$ so that $x\in {\mathscr{D}}(A)$,
$x\in {\mathscr{D}}(B)$ and $Ax=AE_nx=A_nx$, $Bx=BF_nx=B_nx$.
Thus $x\in {\mathscr{D}}(CA+CB)$.
It follows that $CA+CB$ is densely def\/ined.

If~$U'$ is a~unitary operator in ${\mathcal R}'$, for $x\in {\mathscr{D}}((A+B)C)$,
\begin{gather*}
(A+B)CU'x=(A+B)U'Cx=AU'Cx+BU'Cx=U'ACx+U'BCx
\\
\phantom{(A+B)CU'x}
=U'(ACx+BCx)=U'(A+B)Cx
\end{gather*}
and for $x\in {\mathscr{D}}(CA+CB)$,
\begin{gather*}
(CA+CB)U'x=CAU'x+CBU'x=CU'Ax+CU'Bx
\\
\phantom{(CA+CB)U'x}
=U'CAx+U'CBx=U'(CA+CB)x.
\end{gather*}
From Remark~\ref{remark: check core}, $\overline{(A+B)C}$ and $\overline{CA+CB}$ are af\/f\/iliated with~${\mathcal R}$.
\end{proof}
\begin{Proposition}
If operators $A$ and $B$ are affiliated with~${\mathcal R}$, then
\begin{gather*}
(aA\;\hat{+}\;bB)^*=\bar{a}A^*\;\hat{+}\;\bar{b}B^*
\qquad
\text{and}
\qquad
(A\;\hat{\cdot}\;B)^*=B^*\;\hat{\cdot}\;A^*,
\qquad
a,b\in\mathbb{C},
\end{gather*}
where $*$ is the usual adjoint operation on operators $($possibly unbounded$)$.
\end{Proposition}

\begin{proof}
From Proposition~\ref{prop: hat df}, $aA+bB$ and $AB$ are densely def\/ined and preclosed with closures
$aA\;\hat{+}\;bB$ and $A\;\hat{\cdot}\;B$ (af\/f\/iliated with~${\mathcal R}$), respectively.
Then from Theorem~\ref{thm: preclosed},
\begin{gather}
(aA+bB)^*=(aA\;\hat{+}\;bB)^*,
\qquad
(AB)^*=(A\;\hat{\cdot}\;B)^*.
\label{3.2}
\end{gather}
At the same time,
\begin{gather*}
\bar{a}A^*+\bar{b}B^*\subseteq(aA+bB)^*,
\qquad
B^*A^*\subseteq(AB)^*;
\end{gather*}
and both $(aA+bB)^*$ and $(AB)^*$ are closed (Remark~\ref{selfadjoint}).
We also have $\bar{a}A^*\;\hat{+}\;\bar{b}B^*$ and $B^*\;\hat{\cdot}\;A^*$ as the closures (smallest closed
extensions) of~$\bar{a}A^*+\bar{b}B^*$ and $B^*A^*$, respectively.
It follows that
\begin{gather}
\bar{a}A^*+\bar{b}B^*\subseteq\bar{a}A^*\;\hat{+}\;\bar{b}B^*\subseteq(aA+bB)^*,
\qquad
B^*A^*\subseteq B^*\;\hat{\cdot}\;A^*\subseteq(AB)^*.
\label{3.4}
\end{gather}
Now,~\eqref{3.2} together with~\eqref{3.4},
\begin{gather*}
\bar{a}A^*\;\hat{+}\;\bar{b}B^*\subseteq(aA\;\hat{+}\;bB)^*,
\qquad
B^*\;\hat{\cdot}\;A^*\subseteq(A\;\hat{\cdot}\;B)^*.
\end{gather*}
Since $\bar{a}A^*\;\hat{+}\;\bar{b}B^*$, $(aA\;\hat{+}\;bB)^*$, $B^*\;\hat{\cdot}\;A^*$ and
$(A\;\hat{\cdot}\;B)^*$ are all af\/f\/iliated with~${\mathcal R}$, from Proposition~\ref{prop: equal},
$\bar{a}A^*\;\hat{+}\;\bar{b}B^*=(aA\;\hat{+}\;bB)^*$ and $B^*\;\hat{\cdot}\;A^*=(A\;\hat{\cdot}\;B)^*$.
\end{proof}
\begin{Theorem}
The family ${\mathscr{A}}_{\rm f}({\mathcal R})$ is a~$*$ algebra $($with unit $I)$ when provided with the
operations~$\hat{+}$ $($addition$)$ and~$\hat{\cdot}$ $($multiplication$)$.
\end{Theorem}
\begin{Definition}
 We call ${\mathscr{A}}_{\rm f}({\mathcal R})$, the $*$ algebra of operators af\/f\/iliated with
a~f\/inite von Neumann algebra ${\mathcal R}$, {\it the Murray--von Neumann algebra} associated with~${\mathcal R}$.
\end{Definition}

\section{The Heisenberg--von Neumann puzzle}

The Heisenberg--von Neumann puzzle asks whether there is a~representation of the Heisenberg commutation
relation in terms of unbounded operators af\/f\/iliated with a~factor of Type~II$_1$.
Recall that factors are von Neumann algebras whose centers consist of scalar multiples of the identity
operator~$I$.
A von Neumann algebra is said to be f\/inite when the identity operator $I$ is f\/inite.
Factors without minimal projections in which $I$ is f\/inite are said to be of ``Type~II$_1$''.
So, factors of Type II$_1$ are f\/inite von Neumann algebras.
As noted in Section~\ref{s6}, the operators af\/f\/iliated with a~f\/inite von Neumann algebra ${\mathcal
R}$ have special properties and they form an algebra~${\mathscr{A}}_{\rm f}({\mathcal R})$ (the Murray--von
Neumann algebra associated with~${\mathcal R}$).
Von Neumann had great respect for his physicist colleagues and the uncanny accuracy of their results in
experiments at the subatomic level.
In ef\/fect, the physicists worked with unbounded operators, but in a~loose way.
If taken at face value, many of their mathematical assertions were demonstrably incorrect.
When the algebra ${\mathscr{A}}_{\rm f}({{\mathcal M}})$, with ${\mathcal M}$ a~factor of Type~II$_1$,
appeared, von Neumann hoped that it would provide a~framework for the formal computations the physicists
made with the unbounded ope\-ra\-tors.
As it turned out, in more advanced areas of modern physics, factors of Type~II$_1$ do not suf\/f\/ice, by
themselves, for the mathematical framework needed.
It remains a~tantalizing question, nonetheless, whether the most fundamental relation of quantum mechanics,
the Heisenberg relation, can be realized with self-adjoint operators in some~${\mathscr{A}}_{\rm
f}({{\mathcal M}})$.
\begin{Lemma}\label{lemma: TB}
Suppose that $T$ is a~closed operator on the Hilbert space ${\mathcal H}$ and $B\in {\mathcal B}({\mathcal H})$.
Then the operator $TB$ is closed.
\end{Lemma}

\begin{proof}
Suppose $(x_n,y_n)\in \mathscr{G}(TB)$ and $x_n\rightarrow x$, $y_n=TBx_n\rightarrow y$.
We show that $(x,y)\in \mathscr{G}(TB)$.
By assumption, $Bx_n\in \mathscr{D}(T)$.
Since $B$ is bounded (hence, continuous), $Bx_n\rightarrow Bx$.
Since $T$ is closed and $TBx_n=y_n\rightarrow y$, we have that $(Bx,y)\in \mathscr{G}(T)$, so that $Bx\in
\mathscr{D}(T)$ and $TBx=y$.
Hence $(x,y)\in \mathscr{G}(TB)$ and $TB$ is closed.
\end{proof}

\begin{Remark}
{\rm With $T$ and $B$ as in the preceding lemma, the operator $BT$ is not necessarily closed in general,
even not preclosed (closable).}
\end{Remark}

Consider the following example.
Let $\{y_1,y_2, y_3, \dotsc\}$ be an orthonormal basis for a~Hilbert space~${\mathcal H}$, and let
\begin{gather*}
{\mathscr{D}}=\left\{x\in{\mathcal H}:\sum_{n=1}^\infty n^4|\langle x,y_n\rangle|^2<\infty\right\},
\qquad
z=\sum_{n=1}^\infty n^{-1}y_n.
\end{gather*}
Def\/ine $B$ in ${\mathcal B}({\mathcal H})$ by $Bx=\langle x,z\rangle z$; and def\/ine mapping $T$ with
domain ${\mathscr{D}}$ by
\begin{gather*}
Tx=\sum_{n=1}^\infty n^2\langle x,y_n\rangle y_n.
\end{gather*}
Note that $T$ is a~closed densely def\/ined operator.
First, ${\mathscr{D}}$ certainly contains the subma\-ni\-fold of all f\/inite linear combinations of the basis
elements $y_1,y_2, y_3, \dotsc$, from which ${\mathscr{D}}$ is dense in ${\mathcal H}$.
Suppose $\{u_m\}$ is a~sequence in ${\mathscr{D}}$ tending to $u$ and $\{Tu_m\}$ converges to $v$.
For $y_{n'}\in \{y_1,y_2, y_3, \dotsc\}$, $\langle Tu_m, y_{n'}\rangle=\langle \sum\limits_{n=1}^\infty
n^2\langle u_m,y_n\rangle y_n, y_{n'}\rangle =n^2\langle u_m,y_{n'}\rangle \rightarrow n^2\langle
u,y_{n'}\rangle$.
But $\langle Tu_m, y_{n'}\rangle\rightarrow \langle v, y_{n'}\rangle$, so that $\langle v,
y_{n'}\rangle=n^2\langle u,y_{n'}\rangle$; and $\sum\limits_{n=1}^\infty|n^2\langle
u,y_n\rangle|^2=\sum\limits_{n=1}^\infty|\langle v,y_n\rangle|^2=\|v\|^2<\infty$.
Thus $u\in {\mathscr{D}}$ and $Tu=\sum\limits_{n=1}^\infty n^2\langle u,y_n\rangle
y_n=\sum\limits_{n=1}^\infty \langle v,y_n\rangle y_n=v$, so that $\mathscr{G}(T)$ is closed.
However, $BT$ is not preclosed.
If~$u_m=m^{-1}y_m$, then $u_m\rightarrow 0$, but
\begin{gather*}
BTu_m=\langle Tu_m,z\rangle z=\left\langle\sum_{n=1}^\infty n^2\langle u_m,y_n\rangle y_n,\sum_{n=1}
^\infty n^{-1}y_n\right\rangle z   
\\
\phantom{BTu_m}
=\left\langle\sum_{n=1}^\infty n^2\langle m^{-1}y_m,y_n\rangle y_n,\sum_{n=1}^\infty n^{-1}y_n\right\rangle z
=\langle my_m,m^{-1}y_m\rangle z=z\neq0.  
\end{gather*}
Hence $BT$ is not preclosed.
(Recall that an operator $S$ is preclosed, i.e.\
$\mathscr{G}(S)^-$ is a~graph of a~linear transformation, if and only if convergence of the sequence
$\{x_n\}$ in $\mathscr{D}(S)$ to~$0$ and $\{Sx_n\}$ to~$z$ implies that~$z=0$.)

\begin{Lemma}\label{special case}
If~${\mathcal R}$ is a~finite von Neumann algebra, $P$ is a~self-adjoint operator affiliated
with~${\mathcal R}$, and $A$ is an operator in ${\mathcal R}$, such that
$P\;\hat{\cdot}\;A\;\hat{-}\;A\;\hat{\cdot}\;P$ is a~bounded operator $B$, necessarily, affiliated
with~${\mathcal R}$ and, hence, in ${\mathcal R}$, then $\tau(B)$, where $\tau$ is the center-valued trace on
${\mathcal R}$, is $0$.
In particular, $B$ is not of the form $aI$ with $a$ some non-zero scalar in this case.
\end{Lemma}
\begin{proof}
Let $E_n$ be the spectral projection for $P$ corresponding to the interval $[-n,n]$ for each positive
integer $n$.
Then $PE_n$ is an everywhere def\/ined bounded self-adjoint operator as is $E_nPE_n$, and $E_nPE_n=PE_n$.
Note, for this, that $E_nP\subseteq PE_n$, so, $E_nP$ is bounded and its closure $E_n\;\hat{\cdot}\;P=PE_n$.
From the (algebraic) properties, established in Section~\ref{s6}, of the Murray--von Neumann algebra
${\mathscr{A}}_{\rm f}({\mathcal R})$ (of operators af\/f\/iliated with~${\mathcal R}$),
\begin{gather*}
E_n\;\hat{\cdot}\;(P\;\hat{\cdot}\;A)\;\hat{\cdot}\;E_n\;\hat{-}\;E_n\;\hat{\cdot}\;(A\;\hat{\cdot}
\;P)\;\hat{\cdot}\;E_n=E_nBE_n;
\end{gather*}
and from Lemma~\ref{lemma: TB},
\begin{gather*}
E_n\;\hat{\cdot}\;(P\;\hat{\cdot}\;A)E_n\;\hat{-}\;E_n\;\hat{\cdot}\;(A\;\hat{\cdot}\;P)E_n=E_nBE_n.
\end{gather*}
(Since $P\;\hat{\cdot}\;A$ and $A\;\hat{\cdot}\;P$ are closed and $E_n$ is bounded,
$(P\;\hat{\cdot}\;A)E_n$ and $(A\;\hat{\cdot}\;P)E_n$ are closed.
Hence they are equal to their closures $(P\;\hat{\cdot}\;A)\;\hat{\cdot}\;E_n$ and
$(A\;\hat{\cdot}\;P)\;\hat{\cdot}\;E_n$, respectively.) Now, since $E_n$, $A$ and
$E_n\;\hat{\cdot}\;P=PE_n=E_nPE_n$ are all bounded,
\begin{gather*}
E_n\;\hat{\cdot}\;(P\;\hat{\cdot}\;A)E_n=(E_n\;\hat{\cdot}\;P)\;\hat{\cdot}\;AE_n=E_nPE_nAE_n=E_nPE_nE_nAE_n
\end{gather*}
and
\begin{gather*}
E_n\;\hat{\cdot}\;(A\;\hat{\cdot}\;P)E_n=(E_n\;\hat{\cdot}\;A)\;\hat{\cdot}
\;(PE_n)=E_nAE_nPE_n=E_nAE_nE_nPE_n.
\end{gather*}
Thus
\begin{gather}
E_nPE_nE_nAE_n-E_nAE_nE_nPE_n=E_nBE_n.
\label{*}
\end{gather}
Since $E_nPE_n$ and $E_nAE_n$ are bounded and in ${\mathcal R}$, the left-hand side of~\eqref{*} is
a~commutator in~${\mathcal R}$.
Hence $\tau(E_nBE_n)=0$.
As $\|E_nBE_n\|\leq\|B\|$, for each $n$, and $E_n\uparrow I$ in the strong-operator topology, $E_nBE_n$ is
strong (hence, weak)-operator convergent to $B$.
From Theorem 8.2.8 of~\cite{KRII}, $\tau$~is ultraweakly continuous on~${\mathcal R}$.
Thus $0=\tau(E_nBE_n)\rightarrow \tau(B)$.
\end{proof}

\begin{Theorem}\label{thm main}
If~${\mathcal R}$ is a~finite von Neumann algebra, $P$ and $Q$ are self-adjoint operators affi\-liated
with~${\mathcal R}$, and $P\;\hat{\cdot}\;Q\;\hat{-}\;Q\;\hat{\cdot}\;P$ is a~bounded operator~$B$, then
$\tau(B)$, where $\tau$ is the center-valued trace on ${\mathcal R}$, is~$0$.
In particular, $P\;\hat{\cdot}\;Q\;\hat{-}\;Q\;\hat{\cdot}\;P$ is not of the form~$aI$ for some non-zero
scalar~$a$.
\end{Theorem}
\begin{proof}
Since $P\;\hat{\cdot}\;Q\;\hat{-}\;Q\;\hat{\cdot}\;P$ is af\/f\/iliated with~${\mathcal R}$, it is, by
def\/inition, closed on its dense domain.
We are given that $B$ is bounded on this domain.
Hence $B$ is everywhere def\/ined.
With $E_n$ as in Lemma~\ref{special case}, we argue as in Lemma~\ref{special case}, with $Q$ in place of~$A$, to conclude that
\begin{gather*}
E_n\;\hat{\cdot}\;(P\;\hat{\cdot}\;Q)\;\hat{\cdot}\;E_n\;\hat{-}\;E_n\;\hat{\cdot}\;(Q\;\hat{\cdot}
\;P)\;\hat{\cdot}\;E_n=E_nBE_n.
\end{gather*}
In this case,
\begin{gather*}
E_n\;\hat{\cdot}\;(P\;\hat{\cdot}\;Q)\;\hat{\cdot}\;E_n=(E_n\;\hat{\cdot}\;P)\;\hat{\cdot}
\;(Q\;\hat{\cdot}\;E_n)=E_nPE_n\;\hat{\cdot}\;(Q\;\hat{\cdot}\;E_n)
\\
\phantom{E_n\;\hat{\cdot}\;(P\;\hat{\cdot}\;Q)\;\hat{\cdot}\;E_n}
\overset{\text{Lemma~\ref{lemma: TB}}}{=} E_nPE_nE_n\;\hat{\cdot}\;QE_n
=E_nPE_n\;\hat{\cdot}\;(E_n\;\hat{\cdot}\;QE_n),
\end{gather*}
and
\begin{gather*}
E_n\;\hat{\cdot}\;(Q\;\hat{\cdot}\;P)\;\hat{\cdot}\;E_n=E_n\;\hat{\cdot}\;Q\;\hat{\cdot}
\;PE_n
\overset{\text{Lemma~\ref{lemma: TB}}}{=} E_n\;\hat{\cdot}\;QPE_n
\\
\phantom{E_n\;\hat{\cdot}\;(P\;\hat{\cdot}\;Q)\;\hat{\cdot}\;E_n}
=E_n\;\hat{\cdot}\;QE_nPE_n=E_n\;\hat{\cdot}\;QE_nE_nPE_n=(E_n\;\hat{\cdot}\;QE_n)\;\hat{\cdot}\;E_nPE_n.
\end{gather*}
Thus
\begin{gather*}
(E_nPE_n)\;\hat{\cdot}\;(E_n\;\hat{\cdot}\;QE_n)\;\hat{-}\;(E_n\;\hat{\cdot}\;QE_n)\;\hat{\cdot}
\;(E_nPE_n)=E_nBE_n.
\end{gather*}
Since $E_nPE_n$ and $E_nBE_n$ are bounded operators in ${\mathcal R}$, Lemma~\ref{special case} applies,
and $\tau(E_nBE_n)=0$.
Again, $E_n\uparrow I$ and $\tau(B)=0$.
It follows that $B$ cannot be $aI$ with $a\neq0$.
\end{proof}
\begin{Corollary}
The Heisenberg relation, $QP-PQ=i\hbar I$, cannot be satisfied with self-adjoint operators~$Q$ and~$P$ in
the algebra of operators affiliated with a~finite von Neumann algebra, in particular, with a~factor
of Type~{\rm II}$_1$.
\end{Corollary}

\begin{Corollary}
Let ${\mathcal R}$ be a~finite von Neumann algebra with the center-valued trace $\tau$.
If~$T\in {\mathscr{A}}_{\rm f}({\mathcal R})$, $A\in {\mathscr{A}}_{\rm f}({\mathcal R})$, $A=A^*$, and
$T\;\hat{\cdot}\;A\;\hat{-}\;A\;\hat{\cdot}\;T=B\in {\mathcal R}$, then $\tau(B)=0$.
\end{Corollary}

\begin{proof}
First, we note that in the $*$ algebra ${\mathscr{A}}_{\rm f}({\mathcal R})$, $T$ can be uniquely written as
$T_1\;\hat{+}\;iT_2$ with $T_1$ $({=}(T\;\hat{+}\;T^*)/2)$ and $T_2$ $({=}(T\;\hat{-}\;T^*)/2i)$ self-adjoint
operators in ${\mathscr{A}}_{\rm f}({\mathcal R})$.
Making use of this and the algebraic properties of~${\mathscr{A}}_{\rm f}({\mathcal R})$, we obtain
\begin{gather*}
B_1+iB_2=B=T\;\hat{\cdot}\;A\;\hat{-}\;A\;\hat{\cdot}\;T
=(T_1\;\hat{+}\;iT_2)\;\hat{\cdot}\;A\;\hat{-}\;A\;\hat{\cdot}\;(T_1\;\hat{+}\;iT_2)
\\
\phantom{B_1+iB_2}
=(T_1\;\hat{\cdot}\;A\;\hat{-}\;A\;\hat{\cdot}\;T_1)\;\hat{-}\;i(A\;\hat{\cdot}\;T_2\;\hat{-}
\;T_2\;\hat{\cdot}\;A).
\end{gather*}
Hence
\begin{gather*}
B_1=-i(A\;\hat{\cdot}\;T_2\;\hat{-}\;T_2\;\hat{\cdot}\;A)
\qquad
\text{and}
\qquad
iB_2=T_1\;\hat{\cdot}\;A\;\hat{-}\;A\;\hat{\cdot}\;T_1.
\end{gather*}
Since $T_1$, $T_2$, and $A$ are self-adjoint operators in ${\mathscr{A}}_{\rm f}({\mathcal R})$, and $B_1$
and $B_2$ are bounded operators in ${\mathcal R}$, from Theorem~\ref{thm main}, $\tau(B_1)=\tau(B_2)=0$.
Thus, $\tau(B)=0$.
\end{proof}
\begin{Corollary}
Let ${\mathcal R}$ be a~finite von Neumann algebra with the center-valued trace $\tau$.
If~$A\in{\mathscr{A}}_{\rm f}({\mathcal R})$ and
$A^*\;\hat{\cdot}\;A\;\hat{-}\;A\;\hat{\cdot}\;A^*=B\in{\mathcal R}$, then $\tau(B)=0$.
\end{Corollary}
\begin{proof}
Write $A=A_1\;\hat{+}\;iA_2$ with $A_1$ and $A_2$ self-adjoint operators in ${\mathscr{A}}_{\rm
f}({\mathcal R})$.
Then
\begin{gather*}
A^*\;\hat{\cdot}\;A\;\hat{-}\;A\;\hat{\cdot}\;A^*=(A_1\;\hat{-}\;iA_2)\;\hat{\cdot}\;(A_1\;\hat{+}
\;iA_2)\;\hat{-}\;(A_1\;\hat{+}\;iA_2)\;\hat{\cdot}\;(A_1\;\hat{-}\;iA_2)
\\
\phantom{A^*\;\hat{\cdot}\;A\;\hat{-}\;A\;\hat{\cdot}\;A^*}
=2i(A_1\;\hat{\cdot}\;A_2\;\hat{-}\;A_2\;\hat{\cdot}\;A_1)=B\in{\mathcal R}.
\end{gather*}
From Theorem~\ref{thm main}, $\tau(B)=0$.
\end{proof}

The question whether or not the Heisenberg relation can be realized with unbounded ope\-ra\-tors (not
necessarily self-adjoint) af\/f\/iliated with a~f\/inite von Neumann algebra remains open.
The authors have obtained some new results.
It is a~work in progress for us.
During our work, we conjectured the following:

{\it Let ${\mathcal R}$ be a~finite von Neumann algebra.
If~$p$ is a~non-commutative polynomial in $n$ variables with the property that, whenever the variables are
replaced by operators in~${\mathcal R}$ the resulting operator in~${\mathcal R}$ has trace~$0$, then,
whenever replacing the variables in~$p$ by operators in ${\mathscr{A}}_{\rm f}({\mathcal R})$ produces
a~bounded operator, necessarily in~${\mathcal R}$, that operator has trace~$0$.}

Recently, in a~joint work with Andreas Thom, we found a~counter example to the conjecture.
That work will appear elsewhere.

\pdfbookmark[1]{References}{ref}
\LastPageEnding

\end{document}